\documentclass[13pt,reqno,a4paper,oneside,11pt]{amsart}%
\usepackage{amsfonts}
\usepackage{amsfonts}
\usepackage{amsfonts}
\usepackage{amsfonts}
\usepackage{amsfonts}
\usepackage{amsfonts}
\usepackage{amsfonts}
\usepackage{amsfonts}
\usepackage{amsfonts}
\usepackage{amsfonts}
\usepackage{amsfonts}
\usepackage{amsfonts}
\usepackage{amsfonts}
\usepackage{amsfonts}
\usepackage{amsfonts}
\usepackage{amsfonts}
\usepackage{amsfonts}
\usepackage{amsfonts}
\usepackage{amsfonts}
\usepackage{amsfonts}
\usepackage{mathrsfs}
\usepackage{mathrsfs}
\usepackage{amsfonts}
\usepackage{amssymb}
\usepackage{amsmath}
\usepackage{amsthm}
\usepackage{graphicx}
\usepackage{color}
\usepackage{extarrows}
\usepackage{cite}
\usepackage[cp1251]{inputenc}
\providecommand{\U}[1]{\protect\rule{.1in}{.1in}}
\newtheorem{theorem}{Theorem}[section]
\newtheorem{corollary}[theorem]{Corollary}
\newtheorem{lemma}[theorem]{Lemma}
\newtheorem{proposition}[theorem]{Proposition}

\newtheorem{definition}{Definition}[section]
\newtheorem{remark}[theorem]{Remark}
\newtheorem{algorithm}[theorem]{Algorithm}
\newtheorem{example}[theorem]{Example}

\theoremstyle{definition}
\theoremstyle{remark}
\numberwithin{equation}{section}

\ifx\pdfoutput\relax\let\pdfoutput=\undefined\fi
\newcount\msipdfoutput
\ifx\pdfoutput\undefined\else
\ifcase\pdfoutput\else
\ifx\paperwidth\undefined\else
\ifdim\paperheight=0pt\relax\else\pdfpageheight\paperheight\fi
\ifdim\paperwidth=0pt\relax\else\pdfpagewidth\paperwidth\fi
\fi\fi\fi
\begin{document}
\pagestyle{myheadings}
\begin{center}
{\Large \textbf{A new generalization  of a system of  two-sided coupled Sylvester-like quaternion tensor equations}}\footnote{This research was supported by
the grants from the National Natural
Science Foundation of China (11971294).
\par

{} First author: Qing-Wen Wang (Q.W. Wang), wqw@t.shu.edu.cn.
\par  Second author: Mahmoud Saad Mehany (M.S. Mehany), mahmoud2006@shu.edu.cn.}

\bigskip

{ \textbf{Qing-Wen Wang$^{a,b}$, Mahmoud Saad Mehany$^{a,c}$}}

{\small
\vspace{0.5cm}

$a.$ Department of Mathematics, Shanghai University, Shanghai 200444, P. R. China}\\
{\small
\vspace{0.1cm}
$b.$ Collaborative Innovation Center for the Marine Artificial Intelligence, Shanghai 200444, P. R. China}\\
{\small
\vspace{0.1cm}
$c.$ Department of Mathematics, Ain Shams University, Cairo, 11566, A.R. Egypt }
\end{center}
\vspace{0.2cm}
\begin{quotation}
\noindent\textbf{Abstract:}
 This study establishes consistency conditions and a general solution for a coupled system that consists of five two-sided Sylvester-like tensor equations in ten quaternion variables throughout the Einstein tensor product. Certain specific cases are thus established. In a direct application, we investigate certain necessary and sufficient conditions for the existence of an $\eta$-Hermitian solution to  five coupled two-sided Sylvester-like quaternion tensor equations. Finally, we present an algorithm and a numerical example to validate the main result.

\vspace{2mm}
\end{quotation}
\section{\textbf{Introduction}}
We introduce certain notations and definitions for convenience. Consider $I_{1},$ ..., $I_{M}$ to be positive integers for the positive integer $M$. An $M$ order tensor $D$ with entry $D_{i_{1}..i_{M}}$ $(1\leq i_{j}\leq I_{j},\  i=1,...,M)$ is a multidimensional array with the subscripts $i_{1},i_{2},..i_{M}$ \cite{FrontMah,Qii,Qi3,c,QII,QII1,QII2,QII3,a,b}. We utilize the notation that $I(M)$ represents $I_{1}\times I_{2}\times...\times I_{M}$. The quaternion concept was investigated by Hamilton in \cite{12h}, and quaternion algebra can be considered a non-commutative skew field. Let $\mathbb{R}$ and $\mathbb{C}$ be the fields of real numbers and complex numbers, respectively, and let $\mathbb{H}$ be the quaternion algebra
\begin{center}
$\mathbb{H}=\{d_{0}+d_{1}\mathbf{i}+d_{2}\mathbf{j}+d_{3}\mathbf{k} \  | \  \mathbf{i}^{2}=\mathbf{j}^{2}=\mathbf{k}^{2}=\mathbf{ijk}=-1, \  d_{0}, d_{1}, d_{2}, d_{3}\in \mathbb{R}\}.$
\end{center}
Let $\mathbb{H}^{I(M)}$ be the set of the order $M$ dimension $I(M)$ tensors over the quaternion algebra $\mathbb{H}$. Tensors are the natural expansions of vectors and matrices. Tensor equations and computations have applications in machine learning, signal processing, mechanics, physics, Markov processes, control theory, numerical analysis, partial differential equations, and engineering problems \cite{1,5}. Tensor decompositions, tensor eigenvalue, and non-negative tensors  \cite{c},\cite{a},\cite{e} have implementations in signal processing, color image processing \cite{13}, quantum mechanics \cite{15}, quaternion tensor computing \cite{16}, Iterative algorithms for solving some tensor equations \cite{Chu11,Chu22,Chu33,TAO1,TAO2,TAO3,Zheng1,Zheng2}.
Let $\mathcal{A} \in \mathbb{H}^{I(N)\times J(N)}$ and $\mathcal{B} \in \mathbb{H}^{J(N)\times K(M)}$, then the Einstein tensor product \cite{23} of tensors $\mathcal{A}$ and $\mathcal{B}$ is denoted by $\mathcal{A}*_{N}\mathcal{B} \in \mathbb{H}^{I(N)\times K(M)}$, where \[(\mathcal{A}*_{N}\mathcal{B})_{i_{1}..i_{N}k_{1}..k_{M}}=
 \sum_{j_{1}...j_{N}}a_{i_{1}...i_{N}j_{1}...j_{N}}b_{j_{1}...j_{N}k_{1}...k_{M}}.\] The operation $*_{N}$ is associative over the set of all quaternion tensors with qualified order.

Let ${\psi}$ be a nonstandard involution of the quaternion algebra  $\mathbb{H}$ (Definition 3.4.5 \cite{Rodman2014}). If $D \in \mathbb{H}^{m\times n}$, then $(D)_{\psi}$ is an $n\times m$ matrix over $\mathbb{H}$ obtained by applying  ${\psi}$ entrywise to the transpose of $D$. let $D$ be an $n\times n$ matrix over $\mathbb{H}$. $D$  is called a  ${\psi}$-Hermitian matrix if  $(D)_{\psi} = D$ (Definition 3.6.1 \cite{Rodman2014}). Took et al. \cite{Took4} introduced an example of a ${\psi}$-Hermitian matrix called an ${\eta}$-Hermitian matrix. For fixed ${\eta} \in \{\mathbf{i},\mathbf{j},\mathbf{k}\}$, A square matrix $A$  is called an ${\eta}$-Hermitian matrix if $D^{\eta^{*}} = D$, where $D^{\eta^{*}} = -{\eta}D^{*}{\eta}$. An ${\eta}$-Hermitian matrix has applications in linear modeling and statistical signal processing  \cite{Took1,Took2,Took3,Took4}. He \cite{18} gave a generalization of an  ${\eta}$-Hermitian matrix. A square quaternion tensor $\mathcal{D}$ is called an ${\eta}$-Hermitian tensor if $\mathcal{D}=\mathcal{D}^{\eta^{*}}$, where $\mathcal{D}^{\eta^{*}}=-\eta \mathcal{D}^{*}\eta$.
We \cite{18} investigated the  consistency conditions and the exact general solution formula for the following two-sided quaternion tensor equations:
\begin{small}
\begin{equation}
\label{1.1aa}
\begin{gathered}
\mathcal{A}_{1}*_{N}\mathcal{X}_{1}*_{M}\mathcal{B}_{1}+\mathcal{A}_{2}*_{N}\mathcal{X}_{2}*_{M}\mathcal{B}_{2}\\
+\mathcal{A}_{2}*_{N}(\mathcal{C}_{3}*_{N}\mathcal{X}_{3}*_{M}\mathcal{D}_{3}
+\mathcal{C}_{4}*_{N}\mathcal{X}_{4}*_{M}\mathcal{D}_{4})*_{M}\mathcal{B}_{1}=\mathcal{E}_{1}.
\end{gathered}
\end{equation}
\end{small}
where $\mathcal{A}_{i}$, $\mathcal{B}_{i}$, $\mathcal{C}_{j}$, $\mathcal{D}_{j}$  $(i=1,2,\ j=3,4),$ and $\mathcal{E}_{1}$  are given quaternion tensors. A tensor equation $\eqref{1.1aa}$  has applications in the discretization of higher-dimension linear partial differential equations, even including its generalizations \cite {19}. Recently, Wang et al. \cite{20} gave a proper extension to  the quaternion tensor equation  \eqref{1.1aa}. They established consistency conditions and a general solution to  the following  coupled two-sided Sylvester-type quaternion system of tensor equations in terms of the Moore–Penrose inverses for certain given tensors:
\begin{align}
\label{1.3aa}
\left\{
\begin{array}{rll}
&\mathcal{A}_{1}*_{N}\mathcal{X}_{1}*_{M}\mathcal{B}_{1}+\mathcal{A}_{2}*_{N}\mathcal{W}*_{M}\mathcal{B}_{2}=\mathcal{E}_{1}\\
&\mathcal{A}_{3}*_{N}\mathcal{Y}_{1}*_{M}\mathcal{B}_{3}+\mathcal{A}_{4}*_{N}\mathcal{W}*_{M}\mathcal{B}_{4}=\mathcal{E}_{2}.
\end{array}
 \right.
\end{align}
This is based on the various uses of quaternions, rank characterizations of some matrix expressions, matrix decompositions, the coupled Sylvester-like quaternion systems of matrix equations \cite{a1,a3,a4,a5,a6,a7,a8,a9,a10,a11,a12,a13,a14,a15,a16,a17,a18,a19,a21,21,Chu10,Chu20,Chu30}, and the theoretical studies  surrounding Sylvester-like quaternion tensor equations. This paper investigates the consistency of and general solution to the following coupled two-sided Sylvester-like quaternion system of tensor equations:
\begin{align}
\label{1.4aa}
\left\{
\begin{array}{rll}
\begin{gathered}
\mathcal{F}_{4}*_{N}\mathcal{Z}_{1}*_{M}\mathcal{G}_{4}=\mathcal{E}_{4},\\
\mathcal{A}_{i}*_{N}\mathcal{X}_{i}*_{M}\mathcal{B}_{i}+\mathcal{C}_{i}*_{N}\mathcal{Y}_{i}*_{M}\mathcal{D}_{i}\\
+\mathcal{C}_{i}*_{N}(\mathcal{F}_{i}*_{N}\mathcal{Z}_{i}*_{M}\mathcal{G}_{i}
+\mathcal{H}_{i}*_{N}\mathcal{Z}_{i+1}*_{M}\mathcal{J}_{i})*_{M}\mathcal{B}_{i}=\mathcal{E}_{i},\\
 \mathcal{H}_{4}*_{N}\mathcal{Z}_{4}*_{M}\mathcal{J}_{4}=\mathcal{E}_{5},
\end{gathered}
\end{array}
  \right.
\end{align}
$(i=\overline{1,3}),$ which gives us a proper generalization of both systems, \eqref{1.1aa} and \eqref{1.3aa}.
As a direct conclusion, we derive certain necessary and sufficient conditions for the consistency of:
\begin{align}
\label{1.5aa}
\left\{
\begin{array}{rll}
\begin{gathered}
\mathcal{F}_{4}*_{N}\mathcal{Z}_{1}*_{M}\mathcal{G}_{4}=\mathcal{E}_{4},\\
\mathcal{F}_{1}*_{N}\mathcal{Z}_{1}*_{M}\mathcal{G}_{1}
+\mathcal{H}_{1}*_{N}\mathcal{Z}_{2}*_{M}\mathcal{J}_{1}=\mathcal{E}_{1},\\
\mathcal{F}_{2}*_{N}\mathcal{Z}_{2}*_{M}\mathcal{G}_{2}
+\mathcal{H}_{2}*_{N}\mathcal{Z}_{3}*_{M}\mathcal{J}_{2}=\mathcal{E}_{2},\\
\mathcal{F}_{3}*_{N}\mathcal{Z}_{3}*_{M}\mathcal{G}_{3}
+\mathcal{H}_{3}*_{N}\mathcal{Z}_{4}*_{M}\mathcal{J}_{3}=\mathcal{E}_{3},\\
 \mathcal{H}_{4}*_{N}\mathcal{Z}_{4}*_{M}\mathcal{J}_{4}=\mathcal{E}_{5}.
\end{gathered}
\end{array}
  \right.
\end{align}
As an implementation of \eqref{1.5aa}, we obtain the consistency conditions for the existence of an  $\eta$-Hermitian solution to the following two-sided quaternion system of tensor equations:
\begin{align}
\label{1.6aa}
\left\{
\begin{array}{rll}
\begin{gathered}
\mathcal{F}_{4}*_{N}\mathcal{Z}_{1}*_{N}\mathcal{F}_{4}^{\eta^{*}}=\mathcal{E}_{4},\\
\mathcal{F}_{1}*_{N}\mathcal{Z}_{1}*_{N}\mathcal{F}_{1}^{\eta^{*}}
+\mathcal{H}_{1}*_{N}\mathcal{Z}_{2}*_{N}\mathcal{H}_{1}^{\eta^{*}}=\mathcal{E}_{1},\\
\mathcal{F}_{2}*_{N}\mathcal{Z}_{2}*_{N}\mathcal{F}_{2}^{\eta^{*}}
+\mathcal{H}_{2}*_{N}\mathcal{Z}_{3}*_{N}\mathcal{H}_{2}^{\eta^{*}}=\mathcal{E}_{2},\\
\mathcal{F}_{3}*_{N}\mathcal{Z}_{3}*_{N}\mathcal{F}_{3}^{\eta^{*}}
+\mathcal{H}_{3}*_{N}\mathcal{Z}_{4}*_{N}\mathcal{H}_{3}^{\eta^{*}}=\mathcal{E}_{3},\\
 \mathcal{H}_{4}*_{N}\mathcal{Z}_{4}*_{N}\mathcal{H}_{4}^{\eta^{*}}=\mathcal{E}_{5}.
\end{gathered}
\end{array}
  \right.
\end{align}
If we set $\mathcal{C}_{i}=\mathcal{B}_{i}=\mathcal{I}$ in \eqref{1.4aa} where $i=\overline{1,3},$ we obtain the following Sylvester-like quaternion system of tensor equations:
\begin{align}
\label{1.7aa}
\left\{
\begin{array}{rll}
\begin{gathered}
\mathcal{A}_{i}*_{N}\mathcal{X}_{i}+\mathcal{Y}_{i}*_{M}\mathcal{D}_{i}
+\mathcal{F}_{i}*_{N}\mathcal{Z}_{i}*_{M}\mathcal{G}_{i}
+\mathcal{H}_{i}*_{N}\mathcal{Z}_{i+1}*_{M}\mathcal{J}_{i}=\mathcal{E}_{i},\\
\mathcal{F}_{4}*_{N}\mathcal{Z}_{1}*_{M}\mathcal{G}_{4}=\mathcal{E}_{4},\
 \mathcal{H}_{4}*_{N}\mathcal{Z}_{4}*_{M}\mathcal{J}_{4}=\mathcal{E}_{5}\ .
\end{gathered}
\end{array}
  \right.
\end{align}
\par
The remainder of this manuscript is described as follows. The concept of an $\eta$-Hermitian quaternion tensor and the Moore–-Penrose inverse for a general tensor are reminiscent of Section 2. Section 3 expresses the general solution to the two-sided Sylvester-type quaternion system of tensor equations  \eqref{1.4aa} when the solvability conditions are applicable. In Section 4, we provide the necessary and sufficient conditions for the existence of a  $\eta$-Hermitian solution to a system \eqref{1.6aa} as a system \eqref{1.5aa} application. We briefly summarize the key results in Section 5.
\section{\textbf{Preliminaries}}
Throughout this paper tensors are considered quaternion tensors. A tensor $\mathcal{C}\in \mathbb{H}^{I(N)\times J(N)}$ is called an even-order tensor. An even-order tensor $\mathcal{C}\in \mathbb{H}^{I(N)\times I(N)}$ is called an even-order square tensor. Let $c \in \mathbb{H}$, then $\overline{c}$ stands for the conjugate of $c$. A quaternion tensor $\mathcal{C}^{*}$= $(\overline{c}_{j_{1}..j_{M}i_{1}..i_{N}}) \in \mathbb{H}^{J(M)\times I(N)}$ calls the conjugate transpose of the tensor $\mathcal{C}$= $(c_{i_{1}..i_{N}j_{1}..j_{M}}) \in \mathbb{H}^{I(N)\times J(M)}$. If $\mathcal{C}=\mathcal{C}^{*}$, then $\mathcal{C}$  is called  Hermitian tensor.
\begin{definition}\cite{23}
An even order square tensor $\mathcal{C} = (c_{i_{1}...i_{M}i_{1}...i_{M}})\in \mathbb{H}^{I(M)\times I(M)}$ is called a diagonal tensor if  $c_{i_{1}...i_{M}i_{1}...i_{M}}\neq 0$ and all its entries are zero. A diagonal tensor is said to be a unit tensor if $c_{i_{1}...i_{M}i_{1}...i_{M}}= 1$, which denotes by $\mathcal{I}$.
\end{definition}
\begin{definition}\cite{23}
Let $\mathcal{C} = (c_{i_{1}...i_{N}j_{1}...j_{M}})\in \mathbb{H}^{I(N)\times J(M)}$, $\mathcal{D} = (d_{i_{1}...i_{N}k_{1}...k_{M}})\in \mathbb{H}^{I(N)\times K(M)}$. The "row block tensor" of $\mathcal{C}$ and $\mathcal{D}$ is denoted by
\begin{align}
\label{2.1}
\begin{array}{rll}
&\begin{pmatrix}\mathcal{C}&\mathcal{D}\end{pmatrix} \in \mathbb{H}^{I(N)\times L(M)},
\end{array}
\end{align}
where $L_{s}=J_{s}+K_{s}$, $s=1, ...,M$ define as
\begin{align*}
&\begin{pmatrix}\mathcal{C}&\mathcal{D}\end{pmatrix}_{i_{1}...i_{N}l_{1}...l_{M}}=
\left\{
\begin{array}{rll}
&c_{i_{1}...i_{N}l_{1}...l_{M}},\ if\ i_{1}...i_{N} \in [I_{1}]\times...\times [I_{N}],\ l_{1}...l_{M} \in [J_{1}]\times...\times [J_{M}],\\
&d_{i_{1}...i_{N}l_{1}...l_{M}},\ if\ i_{1}...i_{N} \in [I_{1}]\times...\times [I_{N}],\ l_{1}...l_{M} \in \Gamma_{1}\times...\times \Gamma_{M},\\
&0,\ \ \ \ \ otherwise,
\end{array}
  \right.
\end{align*}
where $\Gamma_{s}=\{J_{s}+1, ..., J_{s}+K_{s}\}$, $s=1, ...,M$. For a given tensors $\mathcal{A} = (a_{j_{1}...j_{M}i_{1}...i_{N}})\in \mathbb{H}^{J(M)\times I(N)}$, $\mathcal{B} = (b_{k_{1}...k_{M}i_{1}...i_{N}})\in \mathbb{H}^{K(M)\times I(N)}$. The "column block tensor" of $\mathcal{A}$ and $\mathcal{B}$ is denoted by
\begin{align}
\label{2.2}
&\begin{pmatrix}\mathcal{A}\\ \mathcal{B}\end{pmatrix} \in \mathbb{H}^{L(M)\times I(N)},
\end{align}
where $L_{s}=J_{s}+K_{s}$, $s=1, ...,M$ define as
\begin{align*}
&\begin{pmatrix}\mathcal{A}\\ \mathcal{B}\end{pmatrix}_{l_{1}...l_{M}i_{1}...i_{N}}=
\left\{
\begin{array}{rll}
&a_{l_{1}...l_{M}i_{1}...i_{N}},\ if\ l_{1}...l_{M} \in [J_{1}]\times...\times [J_{M}],\ i_{1}...i_{N} \in [I_{1}]\times...\times [I_{N}],\\
&b_{l_{1}...l_{M}i_{1}...i_{N}},\ if\ l_{1}...l_{M} \in \Gamma_{1}\times...\times \Gamma_{M},\ i_{1}...i_{N} \in [I_{1}]\times...\times [I_{N}],\\
&0,\ \ \ \ \ otherwise,
\end{array}
  \right.
\end{align*}
where $\Gamma_{s}=\{J_{s}+1, ..., J_{s}+K_{s}\}$, $s=1, ...,M$.
\end{definition}
\begin{proposition}\cite{23} Let $\mathcal{A} \in \mathbb{H}^{I(P)\times K(N)}$ and $\mathcal{B} \in \mathbb{H}^{K(N)\times J(M)}$. Then
\begin{enumerate}
\item $(\mathcal{A}*_{N}\mathcal{B})^{*}=\mathcal{B}^{*}*_{N}\mathcal{A}^{*}$;
\item $\mathcal{I_{N}}*_{N}\mathcal{B}=\mathcal{B}$,\ $\mathcal{B}*_{M}\mathcal{I}_{M}=\mathcal{B}$,
where $\mathcal{I_{N}} \in \mathbb{H}^{K(N)\times K(N)}$ and $\mathcal{I_{M}} \in \mathbb{H}^{J(M)\times J(M)}$ are units.
\end{enumerate}
\end{proposition}
\begin{proposition}\label{lma 2.22}\cite{23} Consider the tensors $\begin{pmatrix}\mathcal{A}& \mathcal{B}\end{pmatrix}$ and $\begin{pmatrix}\mathcal{C}\\ \mathcal{D}\end{pmatrix}$ given in \eqref{2.1} and \eqref{2.2}. For a given quaternion tensor $\mathcal{G} \in \mathbb{H}^{I(N)\times I(N)}$, we have that
\begin{enumerate}
  \item $\mathcal{G}*_{N}\begin{pmatrix}\mathcal{A}& \mathcal{B}\end{pmatrix}=
\begin{pmatrix}\mathcal{G}*_{N}\mathcal{A}& \mathcal{G}*_{N}\mathcal{B}\end{pmatrix}\in \mathbb{H}^{I(N)\times L(M)},$
  \item $\begin{pmatrix}\mathcal{C}\\ \mathcal{D}\end{pmatrix}*_{N}\mathcal{G}=
\begin{pmatrix}\mathcal{C}*_{N}\mathcal{G}\\ \mathcal{D}*_{N}\mathcal{G}\end{pmatrix}\in \mathbb{H}^{L(M)\times I(N)},$
  \item $\begin{pmatrix}\mathcal{A}& \mathcal{B}\end{pmatrix}*_{M}\begin{pmatrix}\mathcal{C}\\ \mathcal{D}\end{pmatrix}=
\mathcal{A}*_{M}\mathcal{C}+\mathcal{B}*_{M}\mathcal{D}\in \mathbb{H}^{I(N)\times I(N)}.$
\end{enumerate}
\end{proposition}
\begin{definition}\cite{18}
For a given quaternion tensor $\mathcal{D} \in \mathbb{H}^{I(N)\times J(N)}$. The Moore-Penrose inverse of $\mathcal{D}$ is the unique quaternion tensor $\mathcal{X}\in \mathbb{H}^{J(N)\times I(N)}$ satisfies the following axioms:
\begin{enumerate}
 \item $\mathcal{D}*_{N}\mathcal{X}*_{N}\mathcal{D}=\mathcal{D}$,
 \item $\mathcal{X}*_{N}\mathcal{D}*_{N}\mathcal{X}=\mathcal{X}$,
 \item $(\mathcal{D}*_{N}\mathcal{X})^{*}=\mathcal{D}*_{N}\mathcal{X}$,
 \item $(\mathcal{X}*_{N}\mathcal{D})^{*}=\mathcal{X}*_{N}\mathcal{D}$.
 \end{enumerate}
which denotes by $\mathcal{D}^{\dagger}$. Furthermore, $\mathcal{R}_{\mathcal{D}}$ and $\mathcal{L}_{\mathcal{D}}$ denote the projections along $\mathcal{D}$.
\end{definition}
\begin{proposition}\cite{18}
Let $\mathcal{D} \in \mathbb{H}^{I(N)\times I(N)}$. Then
\begin{enumerate}
\item $\mathcal{L}_{\mathcal{D}}*_{N}\mathcal{D}^{\dagger}=\mathcal{D}*_{N}\mathcal{L}_{\mathcal{D}}=0,\
\mathcal{R}_{\mathcal{D}}*_{N}\mathcal{D}=\mathcal{D}^{\dagger}*_{N}\mathcal{R}_{\mathcal{D}}=0$,
\item $(\mathcal{D}^{*})^{\dagger}=(\mathcal{D}^{\dagger})^{*},\ (\mathcal{D}^{\eta^{*}})^{\dagger}=(\mathcal{D}^{\dagger})^{\eta^{*}},$
\item $(\mathcal{L}_{\mathcal{D}})^{\eta^{*}}=\mathcal{R}_{\mathcal{D}^{\eta^{*}}},\ (\mathcal{R}_{\mathcal{D}})^{\eta^{*}}=\mathcal{L}_{\mathcal{D}^{\eta^{*}}},$
\item $(\mathcal{D}^{*}*_{N}\mathcal{D})^{\dagger}=\mathcal{D}^{\dagger}*_{N}(\mathcal{D}^{*})^{\dagger},\
(\mathcal{D}*_{N}\mathcal{D}^{*})^{\dagger}=(\mathcal{D}^{*})^{\dagger}*_{N}\mathcal{D}^{\dagger}.$
\end{enumerate}
\end{proposition}
\begin{lemma}\label{lma 2.3}\cite{18} Let $\mathcal{A}_{1} \in \mathbb{H}^{I(N)\times J(N)}$,\ $\mathcal{A}_{2} \in \mathbb{H}^{I(N)\times G(N)}$,\ $\mathcal{B}_{1} \in \mathbb{H}^{K(M)\times L(M)}$,\ $\mathcal{B}_{2} \in \mathbb{H}^{H(M)\times L(M)}$,\ $\mathcal{C}_{3} \in \mathbb{H}^{G(N)\times Q(N)}$,\ $\mathcal{C}_{4} \in \mathbb{H}^{G(N)\times T(N)}$,\ $\mathcal{D}_{3} \in \mathbb{H}^{S(M)\times K(M)}$,\ $\mathcal{D}_{4} \in \mathbb{H}^{P(M)\times K(M)}$\ and $\mathcal{E}_{1} \in \mathbb{H}^{I(N)\times L(M)}$ be given. Set
\begin{small}
\begin{align*}
\begin{gathered}
\mathcal{M}_{1}=\mathcal{R}_{\mathcal{A}_{1}}*_{N}\mathcal{A}_{2},\ \mathcal{N}_{1}=\mathcal{B}_{2}*_{M}\mathcal{L}_{\mathcal{B}_{1}},\ \mathcal{S}_{1}=\mathcal{A}_{2}*_{N}\mathcal{L}_{\mathcal{M}_{1}},\
\mathcal{\widehat{A}}_{1}=\mathcal{M}_{1}*_{N}\mathcal{C}_{3},\\
\mathcal{\widehat{A}}_{2}=\mathcal{M}_{1}*_{N}\mathcal{C}_{4},\
\mathcal{\widehat{B}}_{1}=\mathcal{D}_{3}*_{M}\mathcal{B}_{1}*_{M}\mathcal{L}_{\mathcal{B}_{2}},\
\mathcal{\widehat{B}}_{2}=\mathcal{D}_{4}*_{M}\mathcal{B}_{1}*_{M}\mathcal{L}_{\mathcal{B}_{2}},\\
\mathcal{\widehat{M}}_{1}=\mathcal{R}_{\mathcal{\widehat{A}}_{1}}*_{N}\mathcal{\widehat{A}}_{2},\ \mathcal{\widehat{N}}_{1}=\mathcal{\widehat{B}}_{2}*_{M}\mathcal{L}_{\mathcal{\widehat{B}}_{1}},\  \mathcal{\widehat{S}}_{1}=\mathcal{\widehat{A}}_{2}*_{N}\mathcal{L}_{\mathcal{\widehat{M}}_{1}},\
\mathcal{\widehat{E}}_{1}=\mathcal{R}_{\mathcal{A}_{1}}*_{N}\mathcal{E}_{1}*_{M}\mathcal{L}_{\mathcal{B}_{2}},\\
\mathcal{\grave{E}}_{1}=\mathcal{E}_{1}-\mathcal{A}_{2}*_{N}(\mathcal{C}_{3}*_{N}\mathcal{X}_{3}*_{M}\mathcal{D}_{3}
+\mathcal{C}_{4}*_{N}\mathcal{W}*_{M}\mathcal{D}_{4})*_{M}\mathcal{B}_{1}.
\end{gathered}
\end{align*}
\end{small}
Then the following statements are equivalent:\\
(1) \eqref{1.1aa} is solvable.\\
(2)
\begin{small}
\begin{align*}
\begin{gathered}
\mathcal{R}_{\mathcal{M}_{1}}*_{N}\mathcal{R}_{\mathcal{A}_{1}}*_{N}\mathcal{E}_{1}=0,\ \mathcal{E}_{1}*_{M}\mathcal{L}_{\mathcal{B}_{1}}*_{M}\mathcal{L}_{\mathcal{N}_{1}}=0,\
\mathcal{R}_{\mathcal{A}_{2}}*_{N}\mathcal{E}_{1}*_{M}\mathcal{L}_{\mathcal{B}_{1}}=0,\\
\mathcal{R}_{\mathcal{\widehat{M}}_{1}}*_{N}\mathcal{R}_{\mathcal{\widehat{A}}_{1}}*_{N}\mathcal{\widehat{E}}_{1}=0,\ \mathcal{\widehat{E}}_{1}*_{M}\mathcal{L}_{\mathcal{\widehat{B}}_{1}}*_{M}\mathcal{L}_{\mathcal{\widehat{N}}_{1}}=0,\\ \mathcal{R}_{\mathcal{\widehat{A}}_{1}}*_{N}\mathcal{\widehat{E}}_{1}*_{M}\mathcal{L}_{\mathcal{\widehat{B}}_{2}}=0,\
\mathcal{R}_{\mathcal{\widehat{A}}_{2}}*_{N}\mathcal{\widehat{E}}_{1}*_{M}\mathcal{L}_{\mathcal{\widehat{B}}_{1}}=0.
\end{gathered}
\end{align*}
\end{small}
In that case, the general solution to  \eqref{1.1aa} can be expressed as follows:
\begin{small}
\begin{subequations}
\begin{align*}
&\begin{array}{l}
 \mathcal{X}_{1}=\mathcal{A}_{1}^{\dagger}*_{N}\mathcal{\grave{E}}_{1}*_{M}\mathcal{B}_{1}^{\dagger}
-\mathcal{A}_{1}^{\dagger}*_{N}\mathcal{A}_{2}*_{N}\mathcal{M}_{1}^{\dagger}*_{N}\mathcal{\grave{E}}_{1}*_{M}\mathcal{B}_{1}^{\dagger}
-\mathcal{A}_{1}^{\dagger}*_{N}\mathcal{S}_{1}*_{N}\mathcal{A}_{2}^{\dagger}\\
\ \ \ \ \ *_{N}\mathcal{\grave{E}}_{1}*_{M}\mathcal{N}_{1}^{\dagger}*_{M}\mathcal{B}_{2}*_{M}\mathcal{B}_{1}^{\dagger}-\mathcal{A}_{1}^{\dagger}*_{N}\mathcal{S}_{1}
*_{N}\mathcal{U}_{2}*_{M}\mathcal{R}_{\mathcal{N}_{1}}*_{M}\mathcal{B}_{2}*_{M}\mathcal{B}_{1}^{\dagger}
\\
\ \ \ \ \ +\mathcal{L}_{\mathcal{A}_{1}}*_{N}\mathcal{U}_{4}+\mathcal{U}_{5}*_{M}\mathcal{R}_{\mathcal{B}_{1}},
 \end{array}  \\
&\begin{array}{l}
 \mathcal{X}_{2}=\mathcal{M}_{1}^{\dagger}*_{N}\mathcal{\grave{E}}_{1}*_{M}\mathcal{B}_{2}^{\dagger}
 +\mathcal{S}_{1}^{\dagger}*_{N}\mathcal{S}_{1}*_{N}\mathcal{A}_{2}^{\dagger}
*_{N}\mathcal{\grave{E}}*_{M}\mathcal{N}_{1}^{\dagger}+\mathcal{L}_{\mathcal{M}_{1}}*_{N}\mathcal{L}_{\mathcal{S}_{1}}\\
\ \ \ \ \ \ *_{N}\mathcal{U}_{1}+\mathcal{L}_{\mathcal{M}_{1}}*_{N}\mathcal{U}_{2}*_{M}\mathcal{R}_{\mathcal{N}_{1}}+\mathcal{U}_{3}*_{M}\mathcal{R}_{\mathcal{B}_{2}},
 \end{array}\\
 &\begin{array}{l}
 \mathcal{X}_{3}=\mathcal{\widehat{A}}_{1}^{\dagger}*_{N}\mathcal{\widehat{E}}_{1}*_{M}\mathcal{\widehat{B}}_{1}^{\dagger}
-\mathcal{\widehat{A}}_{1}^{\dagger}*_{N}\mathcal{\widehat{A}}_{2}*_{N}\mathcal{\widehat{M}}_{1}^{\dagger}*_{N}\mathcal{\widehat{E}}_{1}
*_{M}\mathcal{\widehat{B}}_{1}^{\dagger}
-\mathcal{\widehat{A}}_{1}^{\dagger}*_{N}\mathcal{\widehat{S}}_{1}*_{N}\mathcal{\widehat{A}}_{2}^{\dagger}\\
\ \ \ \ \*_{N}\mathcal{\widehat{E}}_{1} *_{M}\mathcal{\widehat{N}}_{1}^{\dagger}*_{M}\mathcal{\widehat{B}}_{2}*_{M}\mathcal{\widehat{B}}_{1}^{\dagger}-\mathcal{\widehat{A}}_{1}^{\dagger}
*_{N}\mathcal{\widehat{S}}_{1}
*_{N}\mathcal{\widehat{U}}_{2}*_{M}\mathcal{R}_{\mathcal{\widehat{N}}_{1}}*_{M}\mathcal{\widehat{B}}_{2}*_{M}\mathcal{\widehat{B}}_{1}^{\dagger}
\\
\ \ \ \ \ +\mathcal{L}_{\mathcal{\widehat{A}}_{1}}*_{N}\mathcal{\widehat{U}}_{4}+\mathcal{\widehat{U}}_{5}*_{M}\mathcal{R}_{\mathcal{\widehat{B}}_{1}},
 \end{array}\\
&\begin{array}{l}
 \mathcal{X}_{4}=\mathcal{\widehat{M}}_{1}^{\dagger}*_{N}\mathcal{\widehat{E}}_{1}*_{M}\mathcal{\widehat{B}}_{2}^{\dagger}
 +\mathcal{\widehat{S}}_{1}^{\dagger}*_{N}\mathcal{\widehat{S}}_{1}*_{N}\mathcal{\widehat{A}}_{2}^{\dagger}
*_{N}\mathcal{\widehat{E}}*_{M}\mathcal{\widehat{N}}_{1}^{\dagger}+\mathcal{L}_{\mathcal{\widehat{M}}_{1}}
*_{N}\mathcal{L}_{\mathcal{\widehat{S}}_{1}}\\
\ \ \ \ \ \ *_{N}\mathcal{\widehat{U}}_{1}+\mathcal{L}_{\mathcal{\widehat{M}}_{1}}*_{N}\mathcal{\widehat{U}}_{2}*_{M}\mathcal{R}_{\mathcal{\widehat{N}}_{1}}
+\mathcal{\widehat{U}}_{3}*_{M}\mathcal{R}_{\mathcal{\widehat{B}}_{2}},
 \end{array}
\end{align*}
\end{subequations}
\end{small}
where $\mathcal{U}_{i},$ $\mathcal{\widehat{U}}_{i}$ $(i=\overline{1,5})$ are arbitrary  tensors with suitable orders.

In case of $\mathcal{A}_{2}=\mathcal{B}_{1}=\mathcal{I}$ and  $\mathcal{A}_{1}=\mathcal{B}_{2}=0$, we have the following special case of \eqref{1.1aa}
\begin{small}
\begin{align*}
\mathcal{C}_{3}*_{N}\mathcal{X}_{3}*_{M}\mathcal{D}_{3}
+\mathcal{C}_{4}*_{N}\mathcal{X}_{4}*_{M}\mathcal{D}_{4}=\mathcal{E}_{1},
\end{align*}
\end{small}
which is solvable if and only if
\begin{small}
\begin{align*}
\mathcal{R}_{\mathcal{\widehat{M}}_{1}}*_{N}\mathcal{R}_{\mathcal{C}_{3}}*_{N}\mathcal{E}_{1}=0,\ \mathcal{E}_{1}*_{M}\mathcal{L}_{\mathcal{D}_{3}}*_{M}\mathcal{L}_{\mathcal{\widehat{N}}_{1}}=0,\\ \mathcal{R}_{\mathcal{C}_{3}}*_{N}\mathcal{E}_{1}*_{M}\mathcal{L}_{D_{4}}=0,\
\mathcal{R}_{\mathcal{C}_{4}}*_{N}\mathcal{E}_{1}*_{M}\mathcal{L}_{\mathcal{D}_{3}}=0.
\end{align*}
\end{small}
In that case, the general solution can be expressed as follows:
\begin{small}
\begin{subequations}
\begin{align*}
 &\begin{array}{l}
 \mathcal{X}_{3}=\mathcal{C}_{3}^{\dagger}*_{N}\mathcal{E}_{1}*_{M}\mathcal{D}_{3}^{\dagger}
-\mathcal{C}_{3}^{\dagger}*_{N}\mathcal{C}_{4}*_{N}\mathcal{\widehat{M}}_{1}^{\dagger}*_{N}\mathcal{E}_{1}
*_{M}\mathcal{D}_{3}^{\dagger}
-\mathcal{C}_{3}^{\dagger}*_{N}\mathcal{\widehat{S}}_{1}*_{N}\mathcal{C}_{4}^{\dagger}\\
\ \ \ \ \*_{N}\mathcal{E}_{1} *_{M}\mathcal{\widehat{N}}_{1}^{\dagger}*_{M}\mathcal{D}_{4}*_{M}\mathcal{D}_{3}^{\dagger}-\mathcal{C}_{3}^{\dagger}
*_{N}\mathcal{\widehat{S}}_{1}
*_{N}\mathcal{\widehat{U}}_{2}*_{M}\mathcal{R}_{\mathcal{\widehat{N}}_{1}}*_{M}\mathcal{D}_{4}*_{M}\mathcal{D}_{3}^{\dagger}
\\
\ \ \ \ \ +\mathcal{L}_{\mathcal{C}_{3}}*_{N}\mathcal{\widehat{U}}_{4}+\mathcal{\widehat{U}}_{5}*_{M}\mathcal{R}_{\mathcal{D}_{3}},
 \end{array}\\
&\begin{array}{l}
 \mathcal{X}_{4}=\mathcal{\widehat{M}}_{1}^{\dagger}*_{N}\mathcal{E}_{1}*_{M}\mathcal{D}_{4}^{\dagger}
 +\mathcal{\widehat{S}}_{1}^{\dagger}*_{N}\mathcal{\widehat{S}}_{1}*_{N}\mathcal{C}_{4}^{\dagger}
*_{N}\mathcal{E}_{1}*_{M}\mathcal{\widehat{N}}_{1}^{\dagger}+\mathcal{L}_{\mathcal{\widehat{M}}_{1}}
*_{N}\mathcal{L}_{\mathcal{\widehat{S}}_{1}}\\
\ \ \ \ \ \ *_{N}\mathcal{\widehat{U}}_{1}+\mathcal{L}_{\mathcal{\widehat{M}}_{1}}*_{N}\mathcal{\widehat{U}}_{2}*_{M}\mathcal{R}_{\mathcal{\widehat{N}}_{1}}
+\mathcal{\widehat{U}}_{3}*_{M}\mathcal{R}_{\mathcal{D}_{4}},
 \end{array}
\end{align*}
\end{subequations}
\end{small}
\end{lemma}
\section{\textbf{The consistency conditions and the general Solution to $\mathbf{(1.4)}$}}
In the following Theorem, we provide consistency conditions and  general solution of a coupled Two-sided Sylvester-like quaternion system of tensor equations \eqref{1.4aa}.
\begin{theorem}\label{system 3.33AA}
Consider the quaternion system of tensor equations  \eqref{1.4aa}, where
\begin{align*}
&\mathcal{F}_{4} \in \mathbb{H}^{I(N)\times J(N)},\ \mathcal{G}_{4} \in \mathbb{H}^{L(M)\times K(M)},\
 \mathcal{H}_{4} \in \mathbb{H}^{I(N)\times Q(N)},\ \mathcal{J}_{4} \in \mathbb{H}^{S(M)\times K(M)},\\
&\mathcal{E}_{4} \in \mathbb{H}^{I(N)\times K(M)},\ \mathcal{E}_{5} \in \mathbb{H}^{I(N)\times K(M)},\
\mathcal{A}_{1} \in \mathbb{H}^{I(N)\times J(N)},\ \mathcal{A}_{2} \in \mathbb{H}^{I(N)\times Q(N)},\\
&\mathcal{A}_{3} \in \mathbb{H}^{I(N)\times P(N)},\ \mathcal{B}_{1} \in \mathbb{H}^{F(M)\times K(M)},\
\mathcal{B}_{2} \in \mathbb{H}^{G(M)\times K(M)},\ \mathcal{B}_{3} \in \mathbb{H}^{H(M)\times K(M)},\\
&\mathcal{C}_{1} \in \mathbb{H}^{I(N)\times A(N)},\ \mathcal{C}_{2} \in \mathbb{H}^{I(N)\times B(N)},\
\mathcal{C}_{3} \in \mathbb{H}^{I(N)\times C(N)},\ \mathcal{D}_{1} \in \mathbb{H}^{L(M)\times K(M)},\\
&\mathcal{D}_{2} \in \mathbb{H}^{L(M)\times K(M)},\ \mathcal{D}_{3} \in \mathbb{H}^{L(M)\times K(M)},\
\mathcal{F}_{1} \in \mathbb{H}^{A(N)\times J(N)},\ \mathcal{F}_{2} \in \mathbb{H}^{B(N)\times P(N)},\\
&\mathcal{F}_{3} \in \mathbb{H}^{C(N)\times J(N)},\ \mathcal{G}_{1} \in \mathbb{H}^{L(M)\times F(M)},\
\mathcal{G}_{2} \in \mathbb{H}^{Q(M)\times G(M)},\ \mathcal{G}_{3} \in \mathbb{H}^{L(M)\times H(M)},\\
&\mathcal{H}_{1} \in \mathbb{H}^{A(N)\times P(N)},\ \mathcal{H}_{2} \in \mathbb{H}^{B(N)\times J(N)},\
\mathcal{H}_{3} \in \mathbb{H}^{C(N)\times Q(N)},\ \mathcal{J}_{1} \in \mathbb{H}^{Q(M)\times F(M)},\\
&\mathcal{J}_{2} \in \mathbb{H}^{L(M)\times J(M)},\ \mathcal{J}_{3} \in \mathbb{H}^{S(M)\times H(M)},\
\mathcal{E}_{i} \in \mathbb{H}^{I(N)\times K(M)},\ (i=\overline{1,3}).
\end{align*} are given tensors over $\mathbb{H}$. Set
\begin{subequations}
\begin{align}
\label{abc11}
&\mathcal{\grave{E}}_{i}=\mathcal{E}_{i}
-\mathcal{C}_{i}*_{N}(\mathcal{F}_{i}*_{N}\mathcal{Z}_{i}*_{M}\mathcal{G}_{i}
-\mathcal{H}_{i}*_{N}\mathcal{Z}_{i+1}*_{M}\mathcal{J}_{i})*_{M}\mathcal{B}_{i},\\
\label{abc22}
&\mathcal{M}_{i}=\mathcal{R}_{\mathcal{A}_{i}}*_{N}\mathcal{C}_{i},\ \mathcal{N}_{i}=\mathcal{D}_{i}*_{M}\mathcal{L}_{\mathcal{B}_{i}},\ \mathcal{S}_{i}=\mathcal{C}_{i}*_{N}\mathcal{L}_{\mathcal{M}_{i}},\ \mathcal{\widehat{A}}_{i}=\mathcal{M}_{i}*_{N}\mathcal{F}_{i},\\
\label{abc33}
&\mathcal{\widehat{C}}_{i}=\mathcal{M}_{i}*_{N}\mathcal{H}_{i},\
\mathcal{\widehat{B}}_{i}=\mathcal{G}_{i}*_{M}\mathcal{B}_{i}*_{M}\mathcal{L}_{\mathcal{D}_{i}},\
\mathcal{\widehat{D}}_{i}=\mathcal{J}_{i}*_{M}\mathcal{B}_{i}*_{M}\mathcal{L}_{\mathcal{D}_{i}},\
\mathcal{\widehat{M}}_{i}=\mathcal{R}_{\mathcal{\widehat{A}}_{i}}*_{N}\mathcal{\widehat{C}}_{i},\\
\label{abc44}
&\mathcal{\widehat{N}}_{i}=\mathcal{\widehat{D}}_{i}*_{M}\mathcal{L}_{\mathcal{\widehat{B}}_{i}},\  \mathcal{\widehat{S}}_{i}=\mathcal{\widehat{C}}_{i}*_{N}\mathcal{L}_{\mathcal{\widehat{M}}_{i}},\
\mathcal{\widehat{E}}_{i}=\mathcal{R}_{\mathcal{A}_{i}}*_{N}\mathcal{E}_{i}*_{M}\mathcal{L}_{\mathcal{D}_{i}},\ (i=\overline{1,3}),\\
\label{abc55}
&\mathcal{A}_{11}=\begin{bmatrix}\mathcal{L}_{\mathcal{F}_{4}} & -\mathcal{L}_{\mathcal{\widehat{A}}_{1}}\end{bmatrix},\
\mathcal{D}_{11}=\begin{bmatrix}\mathcal{R}_{\mathcal{G}_{4}} \\ -\mathcal{R}_{\mathcal{\widehat{B}}_{1}}\end{bmatrix},\
\mathcal{\widehat{A}}_{11}=\mathcal{\widehat{A}}_{1}^{\dagger}*_{N}\mathcal{\widehat{S}}_{1},\
\mathcal{\widehat{B}}_{11}=R_{\mathcal{\widehat{N}}_{1}}*_{M}\mathcal{\widehat{D}}_{1}*_{M}\mathcal{\widehat{B}}_{1}^{\dagger},\\
\label{abc66}
&\begin{array}{l}
 \mathcal{E}_{11}=\mathcal{\widehat{A}}_{1}^{\dagger}*_{N}\mathcal{\widehat{E}}_{1}*_{M}\mathcal{\widehat{B}}_{1}^{\dagger}
-\mathcal{\widehat{A}}_{1}^{\dagger}*_{N}\mathcal{\widehat{C}}_{1}*_{N}\mathcal{\widehat{M}}_{1}^{\dagger}*_{N}\mathcal{\widehat{E}}_{1}
*_{M}\mathcal{\widehat{B}}_{1}^{\dagger}
-\mathcal{\widehat{A}}_{1}^{\dagger}*_{N}\mathcal{\widehat{S}}_{1}*_{N}\mathcal{\widehat{C}}_{1}^{\dagger}*_{N}\mathcal{\widehat{E}}_{1}\\
\ \ \ \ \ *_{M}\mathcal{\widehat{N}}_{1}^{\dagger}*_{M}\mathcal{\widehat{D}}_{1}*_{M}\mathcal{\widehat{B}}_{1}^{\dagger}
-\mathcal{F}_{4}^{\dagger}*_{N}\mathcal{E}_{4}*_{M}\mathcal{G}_{4}^{\dagger},
 \end{array}\\
\label{abc77}
&\mathcal{A}_{22}=\begin{bmatrix}\mathcal{L}_{\mathcal{H}_{4}} & -\mathcal{L}_{\mathcal{\widehat{M}}_{3}}*_{N}\mathcal{L}_{\mathcal{\widehat{S}}_{3}}\end{bmatrix},\
\mathcal{D}_{22}=\begin{bmatrix}\mathcal{R}_{\mathcal{J}_{4}} \\ -\mathcal{R}_{\mathcal{\widehat{D}}_{3}}\end{bmatrix},\
\mathcal{\widehat{A}}_{22}=\mathcal{L}_{\mathcal{\widehat{M}}_{3}},\
\mathcal{\widehat{B}}_{22}=R_{\mathcal{\widehat{N}}_{3}},\\
\label{abc88}
&\begin{array}{l}
 \mathcal{E}_{22}=\mathcal{\widehat{M}}_{3}^{\dagger}*_{N}\mathcal{\widehat{E}}_{3}*_{M}\mathcal{\widehat{D}}_{3}^{\dagger}
 +\mathcal{\widehat{S}}_{3}^{\dagger}*_{N}\mathcal{\widehat{S}}_{3}*_{N}\mathcal{\widehat{C}}_{3}^{\dagger}
*_{N}\mathcal{\widehat{E}}_{3}*_{M}\mathcal{\widehat{N}}_{3}^{\dagger}
-\mathcal{H}_{4}^{\dagger}*_{N}\mathcal{E}_{5}*_{M}\mathcal{J}_{4}^{\dagger},
 \end{array}\\
 \label{abc99}
&\mathcal{\widehat{\widehat{A}}}_{ii}=\mathcal{R}_{\mathcal{A}_{ii}}*_{N}\mathcal{\widehat{A}}_{ii},\
\mathcal{\widehat{\widehat{B}}}_{ii}=\mathcal{\widehat{B}}_{ii}*_{M}\mathcal{L}_{\mathcal{D}_{ii}},\
\mathcal{\widehat{\widehat{E}}}_{ii}=\mathcal{R}_{\mathcal{A}_{ii}}*_{N}\mathcal{E}_{ii}*_{M}\mathcal{L}_{\mathcal{D}_{ii}},\ (i=1,2),\\
 \label{abc9911}
&\mathcal{\overline{A}}_{1}=\begin{bmatrix}-\mathcal{L}_{\mathcal{\widehat{M}}_{1}}*_{N}\mathcal{L}_{\mathcal{\widehat{S}}_{1}}& \mathcal{L}_{\mathcal{\widehat{A}}_{2}}\end{bmatrix},\ \mathcal{\overline{A}}_{2}=\begin{bmatrix}-\mathcal{L}_{\mathcal{\widehat{M}}_{2}}*_{N}\mathcal{L}_{\mathcal{\widehat{S}}_{2}}& \mathcal{L}_{\mathcal{\widehat{A}}_{3}}\end{bmatrix},\ \mathcal{\overline{F}}_{1}=\mathcal{\widehat{A}}^{\dagger}_{2}*_{N}\mathcal{\widehat{S}}_{2},\\
\label{abc9912}
&\mathcal{\overline{B}}_{1}=\begin{bmatrix}-\mathcal{R}_{\mathcal{\widehat{D}}_{1}}\\ \mathcal{R}_{\mathcal{\widehat{B}}_{2}}\end{bmatrix},\ \mathcal{\overline{B}}_{2}=\begin{bmatrix}-\mathcal{R}_{\mathcal{\widehat{D}}_{2}}\\ \mathcal{R}_{\mathcal{\widehat{B}}_{3}}\end{bmatrix},
\mathcal{\overline{F}}_{2}=\mathcal{\widehat{A}}^{\dagger}_{3}*_{N}\mathcal{\widehat{S}}_{3},\
\mathcal{\overline{G}}_{1}=\mathcal{\widehat{D}}_{2}*_{N}\mathcal{\widehat{B}}^{\dagger}_{2},\
\mathcal{\overline{G}}_{2}=\mathcal{\widehat{D}}_{3}*_{N}\mathcal{\widehat{B}}^{\dagger}_{3},\\
\label{abc99131}
&\mathcal{\overline{H}}_{1}=\mathcal{L}_{\mathcal{\widehat{M}}_{1}},\ \mathcal{\overline{J}}_{1}=\mathcal{R}_{\mathcal{\widehat{N}}_{1}},\
\mathcal{\overline{H}}_{2}=\mathcal{L}_{\mathcal{\widehat{M}}_{2}},\ \mathcal{\overline{J}}_{2}=\mathcal{R}_{\mathcal{\widehat{N}}_{2}},\\
\label{abc99132}
&\begin{array}{l}
 \mathcal{\overline{E}}_{1}=-\mathcal{\widehat{M}}_{1}^{\dagger}*_{N}\mathcal{\widehat{E}}_{1}*_{M}\mathcal{\widehat{D}}_{1}^{\dagger}
 -\mathcal{\widehat{S}}_{1}^{\dagger}*_{N}\mathcal{\widehat{S}}_{1}*_{N}\mathcal{\widehat{C}}_{1}^{\dagger}
*_{N}\mathcal{\widehat{E}}*_{M}\mathcal{\widehat{N}}_{1}^{\dagger}
+\mathcal{\widehat{A}}_{2}^{\dagger}*_{N}\mathcal{\widehat{E}}_{2}*_{M}\mathcal{\widehat{B}}_{2}^{\dagger}
-\mathcal{\widehat{A}}_{2}^{\dagger}*_{N}\mathcal{\widehat{C}}_{2}\\
\ \ \ \ \ \ \ \ \ \ \ *_{N}\mathcal{\widehat{M}}_{2}^{\dagger}*_{N}\mathcal{\widehat{E}}_{2}
*_{M}\mathcal{\widehat{B}}_{2}^{\dagger}
-\mathcal{\widehat{A}}_{2}^{\dagger}*_{N}\mathcal{\widehat{S}}_{2}*_{N}\mathcal{\widehat{C}}_{2}^{\dagger}*_{N}\mathcal{\widehat{E}}_{2}
*_{M}\mathcal{\widehat{N}}_{2}^{\dagger}*_{M}\mathcal{\widehat{D}}_{2}*_{M}\mathcal{\widehat{B}}_{2}^{\dagger},
 \end{array}
 \end{align}
\end{subequations}
\begin{subequations}
\begin{align}
\label{abc991322}
&\begin{array}{l}
 \mathcal{\overline{E}}_{2}=-\mathcal{\widehat{M}}_{2}^{\dagger}*_{N}\mathcal{\widehat{E}}_{2}*_{M}\mathcal{\widehat{D}}_{2}^{\dagger}
 -\mathcal{\widehat{S}}_{2}^{\dagger}*_{N}\mathcal{\widehat{S}}_{2}*_{N}\mathcal{\widehat{C}}_{2}^{\dagger}
*_{N}\mathcal{\widehat{E}}*_{M}\mathcal{\widehat{N}}_{2}^{\dagger}
+\mathcal{\widehat{A}}_{3}^{\dagger}*_{N}\mathcal{\widehat{E}}_{3}*_{M}\mathcal{\widehat{B}}_{3}^{\dagger}
-\mathcal{\widehat{A}}_{3}^{\dagger}*_{N}\mathcal{\widehat{C}}_{3}\\
\ \ \ \ \ \ \ \ \ \ \ *_{N}\mathcal{\widehat{M}}_{3}^{\dagger}*_{N}\mathcal{\widehat{E}}_{3}
*_{M}\mathcal{\widehat{B}}_{3}^{\dagger}
-\mathcal{\widehat{A}}_{3}^{\dagger}*_{N}\mathcal{\widehat{S}}_{3}*_{N}\mathcal{\widehat{C}}_{3}^{\dagger}*_{N}\mathcal{\widehat{E}}_{3}
*_{M}\mathcal{\widehat{N}}_{3}^{\dagger}*_{M}\mathcal{\widehat{D}}_{3}*_{M}\mathcal{\widehat{B}}_{3}^{\dagger},
 \end{array}\\
\label{abxc22}
&\mathcal{\overline{F}}_{ii}=\mathcal{R}_{\mathcal{\overline{A}}_{i}}*_{N}\mathcal{\overline{F}}_{i},\ \mathcal{\overline{G}}_{ii}=\mathcal{\overline{G}}_{i}*_{M}\mathcal{L}_{\mathcal{\overline{B}}_{i}},\ \mathcal{\overline{H}}_{ii}=\mathcal{R}_{\mathcal{\overline{A}}_{i}}*_{N}\mathcal{\overline{H}}_{i},\ \mathcal{\overline{J}}_{ii}=\mathcal{\overline{J}}_{i}*_{M}\mathcal{L}_{\mathcal{\overline{B}}_{i}},\\
\label{abxc33}
&\mathcal{\overline{E}}_{ii}=\mathcal{R}_{\mathcal{\overline{A}}_{i}}*_{N}\mathcal{\overline{E}}_{i}*_{M}\mathcal{L}_{\mathcal{\overline{B}}_{i}},\ \mathcal{\overline{M}}_{ii}=\mathcal{R}_{\mathcal{\overline{F}}_{ii}}*_{N}\mathcal{\overline{H}}_{ii},\ \mathcal{\overline{N}}_{ii}=\mathcal{\overline{J}}_{ii}*_{M}\mathcal{L}_{\mathcal{\overline{G}}_{ii}},\ \mathcal{\overline{S}}_{ii}=\mathcal{\overline{H}}_{ii}*_{N}\mathcal{L}_{\mathcal{\overline{M}}_{ii}},\\
\label{abcD55}
&\mathcal{\overline{\overline{A}}}_{1}=\begin{bmatrix}\mathcal{L}_{\mathcal{\overline{F}}_{11}} &
 -\mathcal{L}_{\mathcal{\overline{M}}_{22}}*_{N}\mathcal{L}_{\mathcal{\overline{S}}_{22}}\end{bmatrix},\
\mathcal{\overline{\overline{B}}}_{1}=\begin{bmatrix}\mathcal{R}_{\mathcal{\overline{G}}_{11}} \\ -\mathcal{R}_{\mathcal{\overline{J}}_{11}}\end{bmatrix},\
\mathcal{\overline{\overline{F}}}_{1}=\mathcal{\overline{F}}_{11}^{\dagger}*_{N}\mathcal{\overline{S}}_{11},\\
\label{abcDE55}
&\mathcal{\overline{\overline{G}}}_{1}=R_{\mathcal{\overline{N}}_{11}}*_{M}\mathcal{\overline{J}}_{11}*_{M}\mathcal{\overline{G}}_{11}^{\dagger},\
\mathcal{\overline{\overline{H}}}_{1}=\mathcal{L}_{\mathcal{\overline{M}}_{22}},\
\mathcal{\overline{\overline{J}}}_{1}=\mathcal{R}_{\mathcal{\overline{N}}_{22}},\\
\label{abcD66}
&\begin{array}{l}
 \mathcal{\overline{\overline{E}}}_{1}=\mathcal{\overline{F}}_{11}^{\dagger}*_{N}\mathcal{\overline{E}}_{11}*_{M}\mathcal{\overline{G}}_{11}^{\dagger}
-\mathcal{\overline{F}}_{11}^{\dagger}*_{N}\mathcal{\overline{H}}_{11}*_{N}\mathcal{\overline{M}}_{11}^{\dagger}*_{N}\mathcal{\overline{E}}_{11}
*_{M}\mathcal{\overline{G}}_{11}^{\dagger}
-\mathcal{\overline{F}}_{11}^{\dagger}*_{N}\mathcal{\overline{S}}_{11}*_{N}\mathcal{\overline{H}}_{11}^{\dagger}\\ \ \ \ \ \ \ \ \ *_{N}\mathcal{\overline{E}}_{11}*_{M}\mathcal{\overline{N}}_{11}^{\dagger}*_{M}\mathcal{\overline{J}}_{11}*_{M}\mathcal{\overline{G}}_{11}^{\dagger}-\mathcal{\overline{M}}_{22}^{\dagger}*_{N}\mathcal{\overline{E}}_{22}*_{M}\mathcal{\overline{J}}_{22}^{\dagger}
 -\mathcal{\overline{S}}_{22}^{\dagger}*_{N}\mathcal{\overline{S}}_{22}*_{N}\mathcal{\overline{H}}_{22}^{\dagger}\\
\ \  \ \ \ \ \ \ \ *_{N}\mathcal{\overline{E}}_{22}*_{M}\mathcal{\overline{N}}_{22}^{\dagger}
 \end{array}\\
\label{abKxc22}
&\mathcal{\overline{\overline{F}}}_{11}=\mathcal{R}_{\mathcal{\overline{\overline{A}}}_{1}}*_{N}\mathcal{\overline{\overline{F}}}_{1},\ \mathcal{\overline{\overline{G}}}_{11}=\mathcal{\overline{\overline{G}}}_{1}*_{M}\mathcal{L}_{\mathcal{\overline{\overline{B}}}_{1}},\ \mathcal{\overline{\overline{H}}}_{11}=\mathcal{R}_{\mathcal{\overline{\overline{A}}}_{1}}*_{N}\mathcal{\overline{\overline{H}}}_{1},\ \mathcal{\overline{\overline{J}}}_{11}=\mathcal{\overline{\overline{J}}}_{1}*_{M}\mathcal{L}_{\mathcal{\overline{\overline{B}}}_{1}},\\
\label{abKKxAc33}
&\mathcal{\overline{\overline{E}}}_{11}=\mathcal{R}_{\mathcal{\overline{\overline{A}}}_{1}}*_{N}\mathcal{\overline{\overline{E}}}_{1}
*_{M}\mathcal{L}_{\mathcal{\overline{\overline{B}}}_{1}},\ \mathcal{\overline{\overline{M}}}_{11}=\mathcal{R}_{\mathcal{\overline{\overline{F}}}_{11}}*_{N}\mathcal{\overline{\overline{H}}}_{11},\ \mathcal{\overline{\overline{N}}}_{11}=\mathcal{\overline{\overline{J}}}_{11}*_{M}\mathcal{L}_{\mathcal{\overline{\overline{G}}}_{11}},\\
\label{abc99911}
&\mathcal{\overline{\overline{S}}}_{11}=\mathcal{\overline{\overline{H}}}_{11}*_{N}\mathcal{L}_{\mathcal{\overline{\overline{M}}}_{11}},\
\mathcal{\widetilde{A}}_{1}=\begin{bmatrix}\mathcal{L}_{\mathcal{\overline{M}}_{11}}*_{N}\mathcal{L}_{\mathcal{\overline{S}}_{11}}& -\mathcal{L}_{\mathcal{\widehat{\widehat{A}}}_{11}}\end{bmatrix},\ \mathcal{\widetilde{A}}_{2}=\begin{bmatrix}\mathcal{L}_{\mathcal{\widehat{\widehat{A}}}_{22}}& -\mathcal{L}_{\mathcal{\overline{F}}_{22}}\end{bmatrix},\\
\label{abc99912}
&\mathcal{\widetilde{B}}_{1}=\begin{bmatrix}\mathcal{R}_{\mathcal{\overline{J}}_{11}}\\-\mathcal{R}_{\mathcal{\widehat{\widehat{B}}}_{11}}\end{bmatrix},\ \mathcal{\widetilde{B}}_{2}=\begin{bmatrix}\mathcal{R}_{\mathcal{\widehat{\widehat{B}}}_{22}}\\ -\mathcal{R}_{\mathcal{\overline{G}}_{22}}\end{bmatrix},
\mathcal{\widetilde{C}}_{1}=\mathcal{L}_{\mathcal{\overline{M}}_{11}}
\mathcal{\widetilde{D}}_{1}=\mathcal{R}_{\mathcal{\overline{N}}_{11}},\\
\label{abc999131}
&\mathcal{\widetilde{C}}_{2}=\mathcal{\overline{F}}^{\dagger}_{22}*_{N}\mathcal{\overline{S}}_{22},\
\mathcal{\widetilde{D}}_{2}=\mathcal{R}_{\mathcal{\overline{N}}_{22}}*_{M}\mathcal{\overline{J}}_{22}*_{M}\mathcal{\overline{G}}^{\dagger}_{22},\\
\label{abcF131}
&\mathcal{\widetilde{E}}_{1}=\widehat{\widehat{\mathcal{A}}}_{11}^{\dagger}*_{N}\widehat{\widehat{\mathcal{E}}}_{11}
*_{M}\widehat{\widehat{\mathcal{B}}}_{11}^{\dagger}-\mathcal{\overline{M}}_{11}^{\dagger}*_{N}\mathcal{\overline{E}}_{11}*_{M}\mathcal{\overline{J}}_{11}^{\dagger}
 -\mathcal{\overline{S}}_{11}^{\dagger}*_{N}\mathcal{\overline{S}}_{11}*_{N}\mathcal{\overline{H}}_{11}^{\dagger}
*_{N}\mathcal{\overline{E}}_{11}*_{M}\mathcal{\overline{N}}_{11}^{\dagger},\\
\label{abcD66}
&\begin{array}{l}
\mathcal{\widetilde{E}}_{2}=\mathcal{\overline{F}}_{22}^{\dagger}*_{N}\mathcal{\overline{E}}_{22}*_{M}\mathcal{\overline{G}}_{22}^{\dagger}
-\mathcal{\overline{F}}_{22}^{\dagger}*_{N}\mathcal{\overline{H}}_{22}*_{N}\mathcal{\overline{M}}_{22}^{\dagger}*_{N}\mathcal{\overline{E}}_{22}
*_{M}\mathcal{\overline{G}}_{22}^{\dagger}
-\mathcal{\overline{F}}_{22}^{\dagger}*_{N}\mathcal{\overline{S}}_{22}*_{N}\mathcal{\overline{H}}_{22}^{\dagger}\\
\ \ \ \ \ *_{N}\mathcal{\overline{E}}_{22}*_{M}\mathcal{\overline{N}}_{22}^{\dagger}*_{M}\mathcal{\overline{J}}_{22}*_{M}\mathcal{\overline{G}}_{22}^{\dagger}
-\widehat{\widehat{\mathcal{A}}}_{22}^{\dagger}*_{N}\widehat{\widehat{\mathcal{E}}}_{22}
*_{M}\widehat{\widehat{\mathcal{B}}}_{22}^{\dagger},
 \end{array}\\
\label{abc999WW11}
& \mathcal{\widetilde{F}}_{1}=\begin{bmatrix}\mathcal{L}_{\mathcal{\widetilde{C}}_{11}}& -\mathcal{L}_{\mathcal{\overline{\overline{F}}}_{11}}\end{bmatrix},\
\mathcal{\widetilde{F}}_{2}=\begin{bmatrix}\mathcal{L}_{\mathcal{\overline{\overline{M}}}_{11}}*_{N}\mathcal{L}_{\mathcal{\overline{\overline{S}}}_{11}}& -\mathcal{L}_{\mathcal{\widetilde{C}}_{22}}\end{bmatrix},\
\mathcal{\widetilde{H}}_{1}=\mathcal{\overline{\overline{F}}}_{11}^{\dagger}*_{N}\mathcal{\overline{\overline{S}}}_{11},\\
\label{abc999WWW12}
&\mathcal{\widetilde{J}}_{1}=\mathcal{R}_{\mathcal{\overline{\overline{N}}}_{11}}
*_{M}\mathcal{\overline{\overline{J}}}_{11}*_{M}\mathcal{\overline{\overline{G}}}_{11}^{\dagger}\
\mathcal{\widetilde{G}}_{1}=\begin{bmatrix}\mathcal{R}_{\mathcal{\widetilde{D}}_{11}}\\-\mathcal{R}_{\mathcal{\overline{\overline{G}}}_{11}}\end{bmatrix},\ \mathcal{\widetilde{G}}_{2}=\begin{bmatrix}\mathcal{R}_{\mathcal{\overline{\overline{J}}}_{11}}\\ -\mathcal{R}_{\mathcal{\widetilde{D}}_{22}}\end{bmatrix},
\mathcal{\widetilde{H}}_{2}=\mathcal{L}_{\mathcal{\overline{\overline{M}}}_{11}},\
\mathcal{\widetilde{J}}_{2}=\mathcal{R}_{\mathcal{\overline{\overline{N}}}_{11}},\\
\label{abcDS66}
&\begin{array}{l}
 \mathcal{\widetilde{\widetilde{E}}}_{1}=\mathcal{\overline{\overline{F}}}_{11}^{\dagger}*_{N}\mathcal{\overline{\overline{E}}}_{11}*_{M}\mathcal{\overline{\overline{G}}}_{11}^{\dagger}
-\mathcal{\overline{\overline{F}}}_{11}^{\dagger}*_{N}\mathcal{\overline{\overline{H}}}_{11}*_{N}\mathcal{\overline{\overline{M}}}_{11}^{\dagger}
*_{N}\mathcal{\overline{\overline{E}}}_{11}
*_{M}\mathcal{\overline{\overline{G}}}_{11}^{\dagger}
-\mathcal{\overline{\overline{F}}}_{11}^{\dagger}*_{N}\mathcal{\overline{\overline{S}}}_{11}\\
\ \ \ \ \ *_{N}\mathcal{\overline{\overline{H}}}_{11}^{\dagger}*_{N}\mathcal{\overline{\overline{E}}}_{11}*_{M}\mathcal{\overline{\overline{N}}}_{11}^{\dagger}*_{M}\mathcal{\overline{\overline{J}}}_{11}
*_{M}\mathcal{\overline{\overline{G}}}_{11}^{\dagger}-\mathcal{\widetilde{C}}_{11}^{\dagger}*_{N}\mathcal{\widetilde{E}}_{11}
*_{M}\mathcal{\widetilde{D}}_{11}^{\dagger},
 \end{array}\\
\label{abcAAD66}
&\mathcal{\widetilde{\widetilde{E}}}_{2}= \mathcal{\widetilde{C}}_{22}^{\dagger}*_{N}\mathcal{\widetilde{E}}_{22}
*_{M}\mathcal{\widetilde{D}}_{22}^{\dagger}
-\mathcal{\overline{\overline{M}}}_{11}^{\dagger}*_{N}\mathcal{\overline{\overline{E}}}_{11}*_{M}\mathcal{\overline{\overline{J}}}_{11}^{\dagger}
 -\mathcal{\overline{\overline{S}}}_{11}^{\dagger}*_{N}\mathcal{\overline{\overline{S}}}_{11}*_{N}\mathcal{\overline{\overline{H}}}_{11}^{\dagger}
*_{N}\mathcal{\overline{\overline{E}}}_{11}*_{M}\mathcal{\overline{\overline{N}}}_{11}^{\dagger},\\
\label{abcDGG66}
&\mathcal{\widetilde{H}}_{11}=\mathcal{R}_{\mathcal{\widetilde{F}}_{1}}*_{N}\mathcal{\widetilde{H}}_{1},\
\mathcal{\widetilde{H}}_{22}=\mathcal{R}_{\mathcal{\widetilde{F}}_{2}}*_{N}\mathcal{\widetilde{H}}_{2},\
\mathcal{\widetilde{J}}_{11}=\mathcal{\widetilde{J}}_{1}*_{M}\mathcal{L}_{\mathcal{\widetilde{G}}_{1}},\
\mathcal{\widetilde{J}}_{22}=\mathcal{\widetilde{J}}_{2}*_{M}\mathcal{L}_{\mathcal{\widetilde{G}}_{2}},\\
\label{abcHHD6611}
&\begin{array}{l}
\mathcal{\widetilde{\widetilde{E}}}_{11}=\mathcal{R}_{\mathcal{\widetilde{F}}_{1}}*_{N}\mathcal{\widetilde{\widetilde{E}}}_{1}
*_{M}\mathcal{L}_{\mathcal{\widetilde{G}}_{1}},\
\mathcal{\widetilde{\widetilde{E}}}_{22}=\mathcal{R}_{\mathcal{\widetilde{F}}_{2}}*_{N}\mathcal{\widetilde{\widetilde{E}}}_{2}
*_{M}\mathcal{L}_{\mathcal{\widetilde{G}}_{2}},\
\mathcal{\widetilde{A}}=\begin{bmatrix}\mathcal{L}_{\mathcal{\widetilde{H}}_{11}}& -\mathcal{L}_{\mathcal{\widetilde{H}}_{22}}\end{bmatrix},\\
\mathcal{\widetilde{B}}=\begin{bmatrix}\mathcal{R}_{\mathcal{\widetilde{J}}_{11}}& -\mathcal{R}_{\mathcal{\widetilde{J}}_{22}}\end{bmatrix},\
\mathcal{\widetilde{E}}=\mathcal{\widetilde{H}}_{22}^{\dagger}*_{N}\mathcal{\widetilde{\widetilde{E}}}_{22}
*_{M}\mathcal{\widetilde{J}}_{22}^{\dagger}-
\mathcal{\widetilde{H}}_{11}^{\dagger}*_{N}\mathcal{\widetilde{\widetilde{E}}}_{11}
*_{M}\mathcal{\widetilde{J}}_{11}^{\dagger}.
 \end{array}
\end{align}
\end{subequations}
Then the system \eqref{1.4aa} is consistent if and only if
\begin{align}
\label{b11}
&\mathcal{R}_{\mathcal{M}_{i}}*_{N}\mathcal{R}_{\mathcal{A}_{i}}*_{N}\mathcal{E}_{i}=0,\ \mathcal{E}_{i}*_{M}\mathcal{L}_{\mathcal{B}_{i}}*_{M}\mathcal{L}_{\mathcal{N}_{i}}=0,\
\mathcal{R}_{\mathcal{C}_{i}}*_{N}\mathcal{E}_{i}*_{M}\mathcal{L}_{\mathcal{B}_{i}}=0,\\
\label{b22}
&\mathcal{R}_{\mathcal{\widehat{M}}_{i}}*_{N}\mathcal{R}_{\mathcal{\widehat{A}}_{i}}*_{N}\mathcal{\widehat{E}}_{i}=0,\ \mathcal{\widehat{E}}_{i}*_{M}\mathcal{L}_{\mathcal{\widehat{B}}_{i}}*_{M}\mathcal{L}_{\mathcal{\widehat{N}}_{i}}=0,
\end{align}
\begin{align}
\label{b331}
&\mathcal{R}_{\mathcal{\widehat{A}}_{i}}*_{N}\mathcal{\widehat{E}}_{i}*_{M}\mathcal{L}_{\mathcal{\widehat{D}}_{i}}=0,\
\mathcal{R}_{\mathcal{\widehat{C}}_{i}}*_{N}\mathcal{\widehat{E}}_{i}*_{M}\mathcal{L}_{\mathcal{\widehat{B}}_{i}}=0,\ (i=\overline{1,3}),\\
\label{b441}
&\mathcal{R}_{\mathcal{F}_{4}}*_{N}\mathcal{E}_{4}=0,\ \mathcal{E}_{4}*_{M}\mathcal{L}_{\mathcal{G}_{4}}=0,\
\mathcal{R}_{\mathcal{H}_{4}}*_{N}\mathcal{E}_{5}=0,\ \mathcal{E}_{5}*_{M}\mathcal{L}_{\mathcal{J}_{4}}=0,\\
\label{b5511}
&\mathcal{R}_{\mathcal{\widehat{\widehat{A}}}_{kk}}*_{N}\mathcal{\widehat{\widehat{E}}}_{kk}=0,\
\mathcal{\widehat{\widehat{E}}}_{kk}*_{M}\mathcal{L}_{\mathcal{\widehat{\widehat{B}}}_{kk}}=0,\ \\
\label{b6611}
&\mathcal{R}_{\mathcal{\overline{M}}_{kk}}*_{N}\mathcal{R}_{\mathcal{\overline{F}}_{kk}}*_{N}\mathcal{\overline{E}}_{kk}=0,\ \mathcal{\overline{E}}_{kk}*_{M}\mathcal{L}_{\mathcal{\overline{G}}_{kk}}*_{M}\mathcal{L}_{\mathcal{\overline{N}}_{kk}}=0,\\
\label{b7711}
&\mathcal{R}_{\mathcal{\overline{F}}_{kk}}*_{N}\mathcal{\overline{E}}_{kk}*_{M}\mathcal{L}_{\mathcal{\overline{J}}_{kk}}=0,\
\mathcal{R}_{\mathcal{\overline{H}}_{kk}}*_{N}\mathcal{\overline{E}}_{kk}*_{M}\mathcal{L}_{\mathcal{\overline{G}}_{kk}}=0,\ (k=1,2),\\
\label{b8811}
&\mathcal{R}_{\mathcal{\overline{\overline{M}}}_{11}}*_{N}\mathcal{R}_{\mathcal{\overline{\overline{F}}}_{11}}*_{N}\mathcal{\overline{\overline{E}}}_{11}=0,\ \mathcal{\overline{\overline{E}}}_{11}*_{M}\mathcal{L}_{\mathcal{\overline{\overline{G}}}_{11}}*_{M}\mathcal{L}_{\mathcal{\overline{\overline{N}}}_{11}}=0,\\
\label{b9911}
&\mathcal{R}_{\mathcal{\overline{\overline{F}}}_{11}}*_{N}\mathcal{\overline{\overline{E}}}_{11}*_{M}\mathcal{L}_{\mathcal{\overline{\overline{J}}}_{11}}=0,\
\mathcal{R}_{\mathcal{\overline{\overline{H}}}_{11}}*_{N}\mathcal{\overline{\overline{E}}}_{11}*_{M}\mathcal{L}_{\mathcal{\overline{\overline{G}}}_{11}}=0,\\
\label{bb111d}
&\mathcal{R}_{\mathcal{\widetilde{C}}_{jj}}*_{N}\mathcal{\widetilde{E}}_{jj}=0,\
\mathcal{\widetilde{E}}_{jj}*_{M}\mathcal{L}_{\mathcal{\widetilde{D}}_{jj}}=0\ (j=1,2),\\
\label{bb1121d}
&\mathcal{R}_{\mathcal{\widetilde{H}}_{ll}}*_{N}\mathcal{\widetilde{\widetilde{E}}}_{ll}=0,\
\mathcal{\widetilde{\widetilde{E}}}_{ll}*_{M}\mathcal{L}_{\mathcal{\widetilde{G}}_{ll}}=0,\
\mathcal{R}_{\mathcal{\widetilde{A}}}*_{N}\mathcal{\widetilde{E}}*_{M}\mathcal{L}_{\mathcal{\widetilde{B}}},\ (l=1,2),\\
\label{bb122121d}
&\mathcal{R}_{\mathcal{\widetilde{A}}}*_{N}\mathcal{\widetilde{E}}*_{M}\mathcal{L}_{\mathcal{\widetilde{B}}}=0.
\end{align}
Under these conditions, the general solution to system \eqref{1.4aa} can be expressed as follows:
\begin{align}
\label{sssa222p1}
&\begin{array}{l}
 \mathcal{X}_{i}=\mathcal{A}_{i}^{\dagger}*_{N}\mathcal{\grave{E}}_{i}*_{M}\mathcal{B}_{i}^{\dagger}
-\mathcal{A}_{i}^{\dagger}*_{N}\mathcal{C}_{i}*_{N}\mathcal{M}_{i}^{\dagger}*_{N}\mathcal{\grave{E}}_{i}*_{M}\mathcal{B}_{i}^{\dagger}
-\mathcal{A}_{i}^{\dagger}*_{N}\mathcal{S}_{i}*_{N}\mathcal{C}_{i}^{\dagger}*_{N}\mathcal{\grave{E}}_{i}\\
\ \ \ \ \ *_{M}\mathcal{N}_{i}^{\dagger}*_{M}\mathcal{D}_{i}*_{M}\mathcal{B}_{i}^{\dagger}-\mathcal{A}_{i}^{\dagger}*_{N}\mathcal{S}_{i}
*_{N}\mathcal{U}_{2i}*_{M}\mathcal{R}_{\mathcal{N}_{i}}*_{M}\mathcal{D}_{i}*_{M}\mathcal{B}_{i}^{\dagger}
+\mathcal{L}_{\mathcal{A}_{i}}*_{N}\mathcal{U}_{4i}\\
\ \ \ \ \ +\mathcal{U}_{5i}*_{M}\mathcal{R}_{\mathcal{B}_{i}},
 \end{array}\\
\label{sssb}
&\begin{array}{l}
 \mathcal{Y}_{i}=\mathcal{M}_{i}^{\dagger}*_{N}\mathcal{\grave{E}}_{i}*_{M}\mathcal{D}_{i}^{\dagger}
 +\mathcal{S}_{i}^{\dagger}*_{N}\mathcal{S}_{i}*_{N}\mathcal{C}_{i}^{\dagger}
*_{N}\mathcal{\grave{E}}_{i}*_{M}\mathcal{N}_{i}^{\dagger}+\mathcal{L}_{\mathcal{M}_{i}}*_{N}\mathcal{L}_{\mathcal{S}_{i}}*_{N}\mathcal{U}_{1i}\\
\ \ \ \ \ \ +\mathcal{L}_{\mathcal{M}_{i}}*_{N}\mathcal{U}_{2i}*_{M}\mathcal{R}_{\mathcal{N}_{i}}+\mathcal{U}_{3i}*_{M}\mathcal{R}_{\mathcal{D}_{1}},
 \end{array}\\
&\mathcal{Z}_{1}=\mathcal{F}_{4}^{\dagger}*_{N}\mathcal{E}_{4}*_{M}\mathcal{G}_{4}^{\dagger}
+\mathcal{L}_{\mathcal{F}_{4}}*_{N}\mathcal{W}_{1}
+\mathcal{W}_{2}*_{M}\mathcal{R}_{\mathcal{G}_{4}},\\
&\mathcal{Z}_{4}=\mathcal{H}_{4}^{\dagger}*_{N}\mathcal{E}_{5}*_{M}\mathcal{J}_{4}^{\dagger}
+\mathcal{L}_{\mathcal{H}_{4}}*_{N}\mathcal{\acute{W}}_{1}
+\mathcal{W}_{3}*_{M}\mathcal{R}_{\mathcal{J}_{4}}, \\
&\begin{array}{l}
 \mathcal{Z}_{2}=\mathcal{\widehat{M}}_{1}^{\dagger}*_{N}\mathcal{\widehat{E}}_{1}*_{M}\mathcal{\widehat{D}}_{1}^{\dagger}
 +\mathcal{\widehat{S}}_{1}^{\dagger}*_{N}\mathcal{\widehat{S}}_{1}*_{N}\mathcal{\widehat{C}}_{1}^{\dagger}
*_{N}\mathcal{\widehat{E}}_{1}*_{M}\mathcal{\widehat{N}}_{1}^{\dagger}+\mathcal{L}_{\mathcal{\widehat{M}}_{1}}
*_{N}\mathcal{L}_{\mathcal{\widehat{S}}_{1}}*_{N}\mathcal{\widehat{U}}_{1}\\
\ \ \ \ \ \ +\mathcal{L}_{\mathcal{\widehat{M}}_{1}}*_{N}\mathcal{\widehat{U}}_{2}*_{M}\mathcal{R}_{\mathcal{\widehat{N}}_{1}}
+\mathcal{\widehat{U}}_{3}*_{M}\mathcal{R}_{\mathcal{\widehat{D}}_{1}},
 \end{array}\\
&\begin{array}{l}
or\ \mathcal{Z}_{2}=\mathcal{\widehat{A}}_{2}^{\dagger}*_{N}\mathcal{\widehat{E}}_{2}*_{M}\mathcal{\widehat{B}}_{2}^{\dagger}
-\mathcal{\widehat{A}}_{2}^{\dagger}*_{N}\mathcal{\widehat{C}}_{2}*_{N}\mathcal{\widehat{M}}_{2}^{\dagger}*_{N}\mathcal{\widehat{E}}_{2}
*_{M}\mathcal{\widehat{B}}_{2}^{\dagger}
-\mathcal{\widehat{A}}_{2}^{\dagger}*_{N}\mathcal{\widehat{S}}_{2}*_{N}\mathcal{\widehat{C}}_{2}^{\dagger}*_{N}\mathcal{\widehat{E}}_{2}\\
\ \ \ \ \ *_{M}\mathcal{\widehat{N}}_{2}^{\dagger}*_{M}\mathcal{\widehat{D}}_{2}*_{M}\mathcal{\widehat{B}}_{2}^{\dagger}-\mathcal{\widehat{A}}_{2}^{\dagger}
*_{N}\mathcal{\widehat{S}}_{2}
*_{N}\mathcal{\widehat{V}}_{2}*_{M}\mathcal{R}_{\mathcal{\widehat{N}}_{2}}*_{M}\mathcal{\widehat{D}}_{2}*_{M}\mathcal{\widehat{B}}_{2}^{\dagger}
+\mathcal{L}_{\mathcal{\widehat{A}}_{2}}*_{N}\mathcal{\widehat{V}}_{4}\\
\ \ \ \ \ +\mathcal{\widehat{V}}_{5}*_{M}\mathcal{R}_{\mathcal{\widehat{B}}_{2}},
 \end{array}\\
&\begin{array}{l}
 \mathcal{Z}_{3}=\mathcal{\widehat{M}}_{2}^{\dagger}*_{N}\mathcal{\widehat{E}}_{2}*_{M}\mathcal{\widehat{D}}_{2}^{\dagger}
 +\mathcal{\widehat{S}}_{2}^{\dagger}*_{N}\mathcal{\widehat{S}}_{2}*_{N}\mathcal{\widehat{C}}_{2}^{\dagger}
*_{N}\mathcal{\widehat{E}}_{2}*_{M}\mathcal{\widehat{N}}_{2}^{\dagger}+\mathcal{L}_{\mathcal{\widehat{M}}_{2}}
*_{N}\mathcal{L}_{\mathcal{\widehat{S}}_{2}}*_{N}\mathcal{\widehat{V}}_{1}\\
\ \ \ \ \ \ +\mathcal{L}_{\mathcal{\widehat{M}}_{2}}*_{N}\mathcal{\widehat{V}}_{2}*_{M}\mathcal{R}_{\mathcal{\widehat{N}}_{2}}
+\mathcal{\widehat{V}}_{3}*_{M}\mathcal{R}_{\mathcal{\widehat{D}}_{2}},
 \end{array}\\
\label{sssa333}
&\begin{array}{l}
or\  \mathcal{Z}_{3}=\mathcal{\widehat{A}}_{3}^{\dagger}*_{N}\mathcal{\widehat{E}}_{3}*_{M}\mathcal{\widehat{B}}_{3}^{\dagger}
-\mathcal{\widehat{A}}_{3}^{\dagger}*_{N}\mathcal{\widehat{C}}_{3}*_{N}\mathcal{\widehat{M}}_{3}^{\dagger}*_{N}\mathcal{\widehat{E}}_{3}
*_{M}\mathcal{\widehat{B}}_{3}^{\dagger}
-\mathcal{\widehat{A}}_{3}^{\dagger}*_{N}\mathcal{\widehat{S}}_{3}*_{N}\mathcal{\widehat{C}}_{3}^{\dagger}*_{N}\mathcal{\widehat{E}}_{3}\\
\ \ \ \ \ *_{M}\mathcal{\widehat{N}}_{3}^{\dagger}*_{M}\mathcal{\widehat{D}}_{3}*_{M}\mathcal{\widehat{B}}_{3}^{\dagger}-\mathcal{\widehat{A}}_{3}^{\dagger}
*_{N}\mathcal{\widehat{S}}_{3}
*_{N}\mathcal{\widehat{K}}_{2}*_{M}\mathcal{R}_{\mathcal{\widehat{N}}_{3}}*_{M}\mathcal{\widehat{D}}_{3}*_{M}\mathcal{\widehat{B}}_{3}^{\dagger}
+\mathcal{L}_{\mathcal{\widehat{A}}_{3}}*_{N}\mathcal{\widehat{K}}_{4}\\
\ \ \ \ \ +\mathcal{\widehat{K}}_{5}*_{M}\mathcal{R}_{\mathcal{\widehat{B}}_{3}},\ (i=\overline{1,3}),\ \ \  where
 \end{array}\\
 \label{Xyqs1}
&\mathcal{W}_{1}=\begin{bmatrix}\mathcal{I}& 0\end{bmatrix}*_{N}[\mathcal{A}_{11}^{\dagger}*_{N}
(\mathcal{E}_{11}-\mathcal{\widehat{A}}_{11}*_{N}\mathcal{\widehat{U}}_{2}*_{M}\mathcal{\widehat{B}}_{11})
-\mathcal{V}_{11}*_{M}\mathcal{D}_{11}+\mathcal{L}_{\mathcal{A}_{11}}*_{N}\mathcal{V}_{22}],\\
\label{Xyqs2}
&\mathcal{\widehat{U}}_{4}=\begin{bmatrix}0&\mathcal{I}\end{bmatrix}*_{N}[\mathcal{A}_{11}^{\dagger}*_{N}
(\mathcal{E}_{11}-\mathcal{\widehat{A}}_{11}*_{N}\mathcal{\widehat{U}}_{2}*_{M}\mathcal{\widehat{B}}_{11})
-\mathcal{V}_{11}*_{M}\mathcal{D}_{11}+\mathcal{L}_{\mathcal{A}_{11}}*_{N}\mathcal{V}_{22}],\\
\label{Xyqs3}
&\begin{array}{l}
\mathcal{W}_{2}=[\mathcal{R}_{\mathcal{A}_{11}}*_{N}
(\mathcal{E}_{11}-\mathcal{\widehat{A}}_{11}*_{N}\mathcal{\widehat{U}}_{2}*_{M}\mathcal{\widehat{B}}_{11})*_{M}\mathcal{D}_{11}^{\dagger}
+\mathcal{A}_{11}*_{N}\mathcal{V}_{11}\\
\ \ \ \ \ \ \ \ \ +\mathcal{V}_{33}*_{M}\mathcal{R}_{\mathcal{D}_{11}}]*_{M}\begin{bmatrix}\mathcal{I}\\ 0\end{bmatrix},
\end{array}
\end{align}
\begin{subequations}
\begin{align}
\label{Xyqs4}
&\begin{array}{l}
\mathcal{\widehat{U}}_{5}=[\mathcal{R}_{\mathcal{A}_{11}}*_{N}
(\mathcal{E}_{11}-\mathcal{\widehat{A}}_{11}*_{N}\mathcal{\widehat{U}}_{2}*_{M}\mathcal{\widehat{B}}_{11})*_{M}\mathcal{D}_{11}^{\dagger}
+\mathcal{A}_{11}*_{N}\mathcal{V}_{11}\\
\ \ \ \ \ \ \ +\mathcal{V}_{33}*_{M}\mathcal{R}_{\mathcal{D}_{11}}]*_{M}\begin{bmatrix}0\\ \mathcal{I}\end{bmatrix},
\end{array}\\
\label{Xyqs5}
&\mathcal{\grave{W}}_{1}=\begin{bmatrix}\mathcal{I}& 0\end{bmatrix}*_{N}[\mathcal{A}_{22}^{\dagger}*_{N}
(\mathcal{E}_{22}-\mathcal{\widehat{A}}_{22}*_{N}\mathcal{\widehat{K}}_{2}*_{M}\mathcal{\widehat{B}}_{22})
-\mathcal{V}_{44}*_{M}\mathcal{D}_{22}+\mathcal{L}_{\mathcal{A}_{22}}*_{N}\mathcal{V}_{55}],\\
\label{Xyqs6}
&\mathcal{\widehat{K}}_{1}=\begin{bmatrix}0&\mathcal{I}\end{bmatrix}*_{N}[\mathcal{A}_{22}^{\dagger}*_{N}
(\mathcal{E}_{22}-\mathcal{\widehat{A}}_{22}*_{N}\mathcal{\widehat{K}}_{2}*_{M}\mathcal{\widehat{B}}_{22})
-\mathcal{V}_{44}*_{M}\mathcal{D}_{22}+\mathcal{L}_{\mathcal{A}_{22}}*_{N}\mathcal{V}_{55}],\\
\label{Xyqs7}
&\begin{array}{l}
\mathcal{W}_{3}=[\mathcal{R}_{\mathcal{A}_{22}}*_{N}
(\mathcal{E}_{22}-\mathcal{\widehat{A}}_{22}*_{N}\mathcal{\widehat{K}}_{2}*_{M}\mathcal{\widehat{B}}_{22})*_{M}\mathcal{D}_{22}^{\dagger}
+\mathcal{A}_{22}*_{N}\mathcal{V}_{44}\\
\ \ \ \ \ \ \ +\mathcal{V}_{66}*_{M}\mathcal{R}_{\mathcal{D}_{22}}]*_{M}\begin{bmatrix}\mathcal{I}\\ 0\end{bmatrix},
\end{array}\\
\label{Xyqs8}
&\begin{array}{l}
\mathcal{\widehat{K}}_{3}=[\mathcal{R}_{\mathcal{A}_{22}}*_{N}
(\mathcal{E}_{22}-\mathcal{\widehat{A}}_{22}*_{N}\mathcal{\widehat{K}}_{2}*_{M}\mathcal{\widehat{B}}_{22})*_{M}\mathcal{D}_{22}^{\dagger}
+\mathcal{A}_{22}*_{N}\mathcal{V}_{44}\\
\ \ \ \ \ \ \ \ \ +\mathcal{V}_{66}*_{M}\mathcal{R}_{\mathcal{D}_{22}}]*_{M}\begin{bmatrix}0\\ \mathcal{I}\end{bmatrix},
\end{array}\\
&\widehat{\mathcal{U}}_{2}=\widehat{\widehat{\mathcal{A}}}_{11}^{\dagger}*_{N}\widehat{\widehat{\mathcal{E}}}_{11}
*_{M}\widehat{\widehat{\mathcal{B}}}_{11}^{\dagger}+\mathcal{L}_{\widehat{\widehat{\mathcal{A}}}_{11}}*_{N}\mathcal{V}_{77}
+\mathcal{V}_{88}*_{M}\mathcal{R}_{\widehat{\widehat{\mathcal{B}}}_{11}},\\
&\widehat{\mathcal{K}}_{2}=\widehat{\widehat{\mathcal{A}}}_{22}^{\dagger}*_{N}\widehat{\widehat{\mathcal{E}}}_{22}
*_{M}\widehat{\widehat{\mathcal{B}}}_{22}^{\dagger}+\mathcal{L}_{\widehat{\widehat{\mathcal{A}}}_{22}}*_{N}\mathcal{V}_{99}
+\mathcal{W}_{11}*_{M}\mathcal{R}_{\widehat{\widehat{\mathcal{B}}}_{22}},\\
\label{Xya991}
&\begin{array}{l}
\mathcal{\widehat{U}}_{1}=\begin{bmatrix}\mathcal{I}& 0\end{bmatrix}*_{N}[\mathcal{\overline{A}}_{1}*_{N}(
\mathcal{\overline{F}}_{1}*_{N}\mathcal{\widehat{V}}_{2}*_{M}\mathcal{\overline{G}}_{1}
+\mathcal{\overline{H}}_{1}*_{N}\mathcal{\widehat{U}}_{2}*_{M}\mathcal{\overline{J}}_{1}-\mathcal{\overline{E}}_{1})\\
\ \ \ \ \ \ \ \ \ \ \ \ \ \ +\mathcal{P}_{11}*_{M}\mathcal{\overline{B}}_{1}+\mathcal{L}_{\mathcal{\overline{A}}_{1}}*_{N}\mathcal{Q}_{11}],
\end{array}\\
\label{Xya992}
&\begin{array}{l}
\mathcal{\widehat{V}}_{4}=\begin{bmatrix}0& \mathcal{I}\end{bmatrix}*_{N}[\mathcal{\overline{A}}_{1}*_{N}(
\mathcal{\overline{F}}_{1}*_{N}\mathcal{\widehat{V}}_{2}*_{M}\mathcal{\overline{G}}_{1}
+\mathcal{\overline{H}}_{1}*_{N}\mathcal{\widehat{U}}_{2}*_{M}\mathcal{\overline{J}}_{1}-\mathcal{\overline{E}}_{1})\\
\ \ \ \ \ \ \ \ \ \ \ \ \ \  \ \ \ +\mathcal{P}_{11}*_{M}\mathcal{\overline{B}}_{1}+\mathcal{L}_{\mathcal{\overline{A}}_{1}}*_{N}\mathcal{Q}_{11}],
\end{array}\\
\label{Xya993}
&\begin{array}{l}
\mathcal{\widehat{U}}_{3}=[\mathcal{R}_{\mathcal{\overline{A}}_{1}}*_{N}
(
\mathcal{\overline{F}}_{1}*_{N}\mathcal{\widehat{V}}_{2}*_{M}\mathcal{\overline{G}}_{1}
+\mathcal{\overline{H}}_{1}*_{N}\mathcal{\widehat{U}}_{2}*_{M}\mathcal{\overline{J}}_{1}-\mathcal{\overline{E}}_{1})\\
\ \ \ \ \ \ \ \ \ \ \ *_{M}
\mathcal{\overline{B}}_{1}^{\dagger}+\mathcal{\overline{A}}_{1}*_{N}\mathcal{P}_{11}+
\mathcal{P}_{33}*_{M}\mathcal{R}_{\mathcal{\overline{B}}_{1}}]*_{M}\begin{bmatrix}\mathcal{I}\\ 0\end{bmatrix},
\end{array}\\
\label{Xya994}
&\begin{array}{l}
\mathcal{\widehat{V}}_{5}=[\mathcal{R}_{\mathcal{\overline{A}}_{1}}*_{N}
(
\mathcal{\overline{F}}_{1}*_{N}\mathcal{\widehat{V}}_{2}*_{M}\mathcal{\overline{G}}_{1}
+\mathcal{\overline{H}}_{1}*_{N}\mathcal{\widehat{U}}_{2}*_{M}\mathcal{\overline{J}}_{1}-\mathcal{\overline{E}}_{1})\\
\ \ \ \ \ \ \ \ \ \ \ *_{M}
\mathcal{\overline{B}}_{1}^{\dagger}+\mathcal{\overline{A}}_{1}*_{N}\mathcal{P}_{11}+
\mathcal{P}_{33}*_{M}\mathcal{R}_{\mathcal{\overline{B}}_{1}}]*_{M}\begin{bmatrix}0\\ \mathcal{I}\end{bmatrix},
\end{array}\\
\label{Xya995}
&\begin{array}{l}
\mathcal{\widehat{V}}_{1}=\begin{bmatrix}\mathcal{I}& 0\end{bmatrix}*_{N}[\mathcal{\overline{A}}_{2}*_{N}(
\mathcal{\overline{F}}_{2}*_{N}\mathcal{\widehat{K}}_{2}*_{M}\mathcal{\overline{G}}_{2}
+\mathcal{\overline{H}}_{2}*_{N}\mathcal{\widehat{V}}_{2}*_{M}\mathcal{\overline{J}}_{2}-\mathcal{\overline{E}}_{2})\\
\ \ \ \ \ \ \ \ \ \ \ \ \ \ +\mathcal{P}_{22}*_{M}\mathcal{\overline{B}}_{2}+\mathcal{L}_{\mathcal{\overline{A}}_{2}}*_{N}\mathcal{Q}_{22}],
\end{array}\\
\label{Xya996}
&\begin{array}{l}
\mathcal{\widehat{K}}_{4}=\begin{bmatrix}0& \mathcal{I}\end{bmatrix}*_{N}[\mathcal{\overline{A}}_{2}*_{N}(
                                 \mathcal{\overline{F}}_{2}*_{N}\mathcal{\widehat{K}}_{2}*_{M}\mathcal{\overline{G}}_{2}
                                 +\mathcal{\overline{H}}_{2}*_{N}\mathcal{\widehat{V}}_{2}*_{M}\mathcal{\overline{J}}_{2}-\mathcal{\overline{E}}_{2})\\
\ \ \ \ \ \ \ \ \ \ \ \ \ \ \ +\mathcal{P}_{22}*_{M}\mathcal{\overline{B}}_{2}+\mathcal{L}_{\mathcal{\overline{A}}_{2}}*_{N}\mathcal{Q}_{22}],
\end{array}\\
\label{Xya997}
&\begin{array}{l}
\mathcal{\widehat{V}}_{3}=[\mathcal{R}_{\mathcal{\overline{A}}_{2}}*_{N}(
                                 \mathcal{\overline{F}}_{2}*_{N}\mathcal{\widehat{K}}_{2}*_{M}\mathcal{\overline{G}}_{2}
                                 +\mathcal{\overline{H}}_{2}*_{N}\mathcal{\widehat{V}}_{2}*_{M}\mathcal{\overline{J}}_{2}-\mathcal{\overline{E}}_{2})\\
\ \ \ \ \ \ \ \ \ \ \ \ *_{M}
                                 \mathcal{\overline{B}}_{2}^{\dagger}+\mathcal{\overline{A}}_{2}*_{N}\mathcal{P}_{22}+
                                 \mathcal{Q}_{33}*_{M}\mathcal{R}_{\mathcal{\overline{B}}_{2}}]*_{M}\begin{bmatrix}\mathcal{I}\\ 0\end{bmatrix},
\end{array}\\
\label{Xya998}
&\begin{array}{l}
\mathcal{\widehat{K}}_{5}=[\mathcal{R}_{\mathcal{\overline{A}}_{2}}*_{N}(\mathcal{\overline{E}}_{2}
                                 +\mathcal{\overline{F}}_{2}*_{N}\mathcal{\widehat{K}}_{2}*_{M}\mathcal{\overline{G}}_{2}
                                 +\mathcal{\overline{H}}_{2}*_{N}\mathcal{\widehat{V}}_{2}*_{M}\mathcal{\overline{J}}_{2})\\
\ \ \ \ \ \ \ \ \ \ \ \ *_{M}
                                 \mathcal{\overline{B}}_{2}^{\dagger}+\mathcal{\overline{A}}_{2}*_{N}\mathcal{P}_{22}+
                                 \mathcal{Q}_{33}*_{M}\mathcal{R}_{\mathcal{\overline{B}}_{2}}]*_{M}\begin{bmatrix}0\\ \mathcal{I}\end{bmatrix},
\end{array}\\
\label{XZyw993}
&\begin{array}{l}
 \mathcal{\widehat{V}}_{2}=\mathcal{\overline{F}}_{11}^{\dagger}*_{N}\mathcal{\overline{E}}_{11}*_{M}\mathcal{\overline{G}}_{11}^{\dagger}
-\mathcal{\overline{F}}_{11}^{\dagger}*_{N}\mathcal{\overline{H}}_{11}*_{N}\mathcal{\overline{M}}_{11}^{\dagger}*_{N}\mathcal{\overline{E}}_{11}
*_{M}\mathcal{\overline{G}}_{11}^{\dagger}
-\mathcal{\overline{F}}_{11}^{\dagger}*_{N}\mathcal{\overline{S}}_{11}*_{N}\mathcal{\overline{H}}_{11}^{\dagger}\\
\ \ \ \ \ *_{N}\mathcal{\overline{E}}_{11}*_{M}\mathcal{\overline{N}}_{11}^{\dagger}*_{M}\mathcal{\overline{J}}_{11}*_{M}\mathcal{\overline{G}}_{11}^{\dagger}
-\mathcal{\overline{F}}_{11}^{\dagger}*_{N}\mathcal{\overline{S}}_{11}
*_{N}\mathcal{P}_{44}*_{M}\mathcal{R}_{\mathcal{\overline{N}}_{11}}*_{M}\mathcal{\overline{J}}_{11}*_{M}\mathcal{\overline{G}}_{11}^{\dagger}
\\
\ \ \ \ \ +\mathcal{L}_{\mathcal{\overline{F}}_{11}}*_{N}\mathcal{P}_{55}+\mathcal{P}_{66}*_{M}\mathcal{R}_{\mathcal{\overline{G}}_{11}},
 \end{array}
\end{align}
\begin{align}
\label{XZyw996}
&\begin{array}{l}
or\  \mathcal{\widehat{V}}_{2}=\mathcal{\overline{M}}_{22}^{\dagger}*_{N}\mathcal{\overline{E}}_{22}*_{M}\mathcal{\overline{J}}_{22}^{\dagger}
 +\mathcal{\overline{S}}_{22}^{\dagger}*_{N}\mathcal{\overline{S}}_{22}*_{N}\mathcal{\overline{H}}_{22}^{\dagger}
*_{N}\mathcal{\overline{E}}_{22}*_{M}\mathcal{\overline{N}}_{22}^{\dagger}+\mathcal{L}_{\mathcal{\overline{M}}_{22}}
\\
\ \ \ \ \ \ *_{N}\mathcal{L}_{\mathcal{\overline{S}}_{22}}*_{N}\mathcal{Q}_{77}+\mathcal{L}_{\mathcal{\overline{M}}_{22}}*_{N}\mathcal{Q}_{55}
*_{M}\mathcal{R}_{\mathcal{\overline{N}}_{22}}+\mathcal{Q}_{88}
*_{M}\mathcal{R}_{\mathcal{\overline{J}}_{22}},
 \end{array}\\
\label{yWW992}
&\begin{array}{l}
\mathcal{P}_{55}=\begin{bmatrix}0& \mathcal{I}\end{bmatrix}*_{N}[\mathcal{\overline{A}}_{1}*_{N}(-\mathcal{\overline{\overline{E}}}_{1}
+\mathcal{\overline{\overline{F}}}_{1}*_{N}\mathcal{P}_{44}*_{M}\mathcal{\overline{\overline{G}}}_{1}
+\mathcal{\overline{\overline{H}}}_{1}*_{N}\mathcal{Q}_{55}*_{M}\mathcal{\overline{\overline{J}}}_{1})\\
\ \ \ \ \ \  \ \ \ \  \ \ \ \  \ \ +\mathcal{K}_{11}*_{M}\mathcal{\overline{\overline{B}}}_{1}
+\mathcal{L}_{\mathcal{\overline{\overline{A}}}_{1}}*_{N}\mathcal{K}_{22}],
\end{array}\\
\label{yWW9S92}
&\begin{array}{l}
\mathcal{Q}_{77}=\begin{bmatrix}\mathcal{I}&0 \end{bmatrix}*_{N}[\mathcal{\overline{A}}_{1}*_{N}(-\mathcal{\overline{\overline{E}}}_{1}
+\mathcal{\overline{\overline{F}}}_{1}*_{N}\mathcal{P}_{44}*_{M}\mathcal{\overline{\overline{G}}}_{1}
+\mathcal{\overline{\overline{H}}}_{1}*_{N}\mathcal{Q}_{55}*_{M}\mathcal{\overline{\overline{J}}}_{1})\\
\ \ \ \ \ \  \ \ \ \  \ \ \ \  \ \ +\mathcal{K}_{11}*_{M}\mathcal{\overline{\overline{B}}}_{1}
+\mathcal{L}_{\mathcal{\overline{\overline{A}}}_{1}}*_{N}\mathcal{K}_{22}],
\end{array}\\
\label{y9WWW93}
&\begin{array}{l}
\mathcal{P}_{66}=[\mathcal{R}_{\mathcal{\overline{\overline{A}}}_{1}}*_{N}
(-\mathcal{\overline{\overline{E}}}_{1}
                                 +\mathcal{\overline{\overline{F}}}_{1}*_{N}\mathcal{P}_{44}*_{M}\mathcal{\overline{\overline{G}}}_{1}
                                 +\mathcal{\overline{\overline{H}}}_{1}*_{N}\mathcal{Q}_{55}*_{M}\mathcal{\overline{\overline{J}}}_{1})\\
\ \ \ \ \ \ \ \ \ \ \ \ \  \ \ \ \ \ *_{M}
\mathcal{\overline{\overline{B}}}_{1}^{\dagger}+\mathcal{\overline{\overline{A}}}_{1}*_{N}\mathcal{K}_{11}+
\mathcal{K}_{33}*_{M}\mathcal{R}_{\mathcal{\overline{\overline{B}}}_{1}}]*_{M}\begin{bmatrix}\mathcal{I}\\ 0\end{bmatrix},
\end{array}\\
\label{y9WWW93}
&\begin{array}{l}
\mathcal{Q}_{88}=[\mathcal{R}_{\mathcal{\overline{\overline{A}}}_{1}}*_{N}
(-\mathcal{\overline{\overline{E}}}_{1}
                                 +\mathcal{\overline{\overline{F}}}_{1}*_{N}\mathcal{P}_{44}*_{M}\mathcal{\overline{\overline{G}}}_{1}
                                 +\mathcal{\overline{\overline{H}}}_{1}*_{N}\mathcal{Q}_{55}*_{M}\mathcal{\overline{\overline{J}}}_{1})\\
\ \ \ \ \ \ \ \ \ \ \ \ \  \ \ \ \ \ *_{M}
\mathcal{\overline{\overline{B}}}_{1}^{\dagger}+\mathcal{\overline{\overline{A}}}_{1}*_{N}\mathcal{K}_{11}+
\mathcal{K}_{33}*_{M}\mathcal{R}_{\mathcal{\overline{\overline{B}}}_{1}}]*_{M}\begin{bmatrix}0\\ \mathcal{I}\end{bmatrix},
\end{array}\\
\label{x992QQQ}
&\mathcal{Q}_{44}= \begin{bmatrix}\mathcal{I}&
0\end{bmatrix}*_{N}[\mathcal{\widetilde{A}}_{1}*_{N}(\mathcal{\widetilde{E}}_{1}
-\mathcal{\widetilde{C}}_{1}*_{N}\mathcal{P}_{44}*_{M}\mathcal{\widetilde{D}}_{1})+
\mathcal{W}_{22}*_{M}\mathcal{\widetilde{B}}_{1}+\mathcal{L}_{\mathcal{\widetilde{A}}_{1}}*_{N}\mathcal{W}_{33}],\\
\label{x992QQQ}
&\mathcal{V}_{77}=\begin{bmatrix}0 &\mathcal{I}\end{bmatrix}*_{N}[\mathcal{\widetilde{A}}_{1}*_{N}(\mathcal{\widetilde{E}}_{1}
-\mathcal{\widetilde{C}}_{1}*_{N}\mathcal{P}_{44}*_{M}\mathcal{\widetilde{D}}_{1})+
\mathcal{W}_{22}*_{M}\mathcal{\widetilde{B}}_{1}+\mathcal{L}_{\mathcal{\widetilde{A}}_{1}}*_{N}\mathcal{W}_{33}],\\
\label{x993qq}
&\mathcal{Q}_{66}=[\mathcal{R}_{\mathcal{\widetilde{A}}_{1}}*_{N}(\mathcal{\widetilde{E}}_{1}
-\mathcal{\widetilde{C}}_{1}*_{N}\mathcal{P}_{44}*_{M}\mathcal{\widetilde{D}}_{1})*_{M}\mathcal{\widetilde{B}}_{1}^{\dagger}+
\mathcal{\widetilde{A}}_{1}*_{N}\mathcal{W}_{22}+\mathcal{W}_{44}*_{M}\mathcal{R}_{\mathcal{\widetilde{B}}_{1}}]*_{M}\begin{bmatrix}\mathcal{I} \\ 0\end{bmatrix},\\
\label{x993qrq}
&\mathcal{V}_{88}=[\mathcal{R}_{\mathcal{\widetilde{A}}_{1}}*_{N}(\mathcal{\widetilde{E}}_{1}
-\mathcal{\widetilde{C}}_{1}*_{N}\mathcal{P}_{44}*_{M}\mathcal{\widetilde{D}}_{1})*_{M}\mathcal{\widetilde{B}}_{1}^{\dagger}+
\mathcal{\widetilde{A}}_{1}*_{N}\mathcal{W}_{22}+\mathcal{W}_{44}*_{M}\mathcal{R}_{\mathcal{\widetilde{B}}_{1}}]*_{M}\begin{bmatrix}0 \\ \mathcal{I}\end{bmatrix},\\
\label{x99qq2LQQ}
&\mathcal{V}_{99} =\begin{bmatrix}\mathcal{I}&
0\end{bmatrix}*_{N}[\mathcal{\widetilde{A}}_{2}*_{N}(\mathcal{\widetilde{E}}_{2}
-\mathcal{\widetilde{C}}_{2}*_{N}\mathcal{P}_{44}*_{M}\mathcal{\widetilde{D}}_{2})+
\mathcal{W}_{55}*_{M}\mathcal{\widetilde{B}}_{2}+\mathcal{L}_{\mathcal{\widetilde{A}}_{2}}*_{N}\mathcal{W}_{66}],
\end{align}
\end{subequations}
\vspace*{-\baselineskip}
\begin{subequations}
\begin{align}
\label{x99qq2LQQ}
&\mathcal{P}_{77}=\begin{bmatrix}0&
\mathcal{I}\end{bmatrix}*_{N}[\mathcal{\widetilde{A}}_{2}*_{N}(\mathcal{\widetilde{E}}_{2}
-\mathcal{\widetilde{C}}_{2}*_{N}\mathcal{P}_{44}*_{M}\mathcal{\widetilde{D}}_{2})+
\mathcal{W}_{55}*_{M}\mathcal{\widetilde{B}}_{2}+\mathcal{L}_{\mathcal{\widetilde{A}}_{2}}*_{N}\mathcal{W}_{66}],\\
\label{xqq9W93}
&\mathcal{W}_{11}=[\mathcal{R}_{\mathcal{\widetilde{A}}_{2}}*_{N}(\mathcal{\widetilde{E}}_{2}
-\mathcal{\widetilde{C}}_{2}*_{N}\mathcal{P}_{44}*_{M}\mathcal{\widetilde{D}}_{2})*_{M}\mathcal{\widetilde{B}}_{2}^{\dagger}+
\mathcal{\widetilde{A}}_{2}*_{N}\mathcal{W}_{55}+\mathcal{W}_{77}*_{M}\mathcal{R}_{\mathcal{\widetilde{B}}_{2}}]*_{M}\begin{bmatrix}\mathcal{I} \\ 0\end{bmatrix},\\
\label{xqq9W93}
&\mathcal{P}_{88}=[\mathcal{R}_{\mathcal{\widetilde{A}}_{2}}*_{N}(\mathcal{\widetilde{E}}_{2}
-\mathcal{\widetilde{C}}_{2}*_{N}\mathcal{P}_{44}*_{M}\mathcal{\widetilde{D}}_{2})*_{M}\mathcal{\widetilde{B}}_{2}^{\dagger}+
\mathcal{\widetilde{A}}_{2}*_{N}\mathcal{W}_{55}+\mathcal{W}_{77}*_{M}\mathcal{R}_{\mathcal{\widetilde{B}}_{2}}]*_{M}\begin{bmatrix}0 \\ \mathcal{I}\end{bmatrix},\\
\label{XZEEyT993}
&\mathcal{P}_{44}=\mathcal{\widetilde{C}}_{11}^{\dagger}*_{N}\mathcal{\widetilde{E}}_{11}
*_{M}\mathcal{\widetilde{D}}_{11}^{\dagger}+\mathcal{L}_{\mathcal{\widetilde{C}}_{11}}*_{N}\mathcal{W}_{88}
+\mathcal{W}_{99}*_{M}\mathcal{R}_{\mathcal{\widetilde{D}}_{11}},\\
\label{XZEEyTTT993}
&\mathcal{Q}_{55}=\mathcal{\widetilde{C}}_{22}^{\dagger}*_{N}\mathcal{\widetilde{E}}_{22}
*_{M}\mathcal{\widetilde{D}}_{22}^{\dagger}+\mathcal{L}_{\mathcal{\widetilde{C}}_{22}}*_{N}\mathcal{T}_{11}
+\mathcal{T}_{22}*_{M}\mathcal{R}_{\mathcal{\widetilde{D}}_{22}},\\
\label{XZEEyTTTr993}
&\mathcal{W}_{88}=\begin{bmatrix}\mathcal{I}&
0\end{bmatrix}*_{N}[\mathcal{\widetilde{F}}_{1}*_{N}(\mathcal{\widetilde{\widetilde{E}}}_{1}
                                 -\mathcal{\widetilde{H}}_{1}*_{N}\mathcal{K}_{44}*_{M}\mathcal{\widetilde{J}}_{1})+
\mathcal{T}_{33}*_{M}\mathcal{\widetilde{G}}_{1}+\mathcal{L}_{\mathcal{\widetilde{F}}_{1}}*_{N}\mathcal{T}_{44}],\\
\label{XZEEyTTTe993}
&\mathcal{K}_{55}=\begin{bmatrix}0&
\mathcal{I}\end{bmatrix}*_{N}[\mathcal{\widetilde{F}}_{1}*_{N}(\mathcal{\widetilde{\widetilde{E}}}_{1}
                                 -\mathcal{\widetilde{H}}_{1}*_{N}\mathcal{K}_{44}*_{M}\mathcal{\widetilde{J}}_{1})+
\mathcal{T}_{33}*_{M}\mathcal{\widetilde{G}}_{1}+\mathcal{L}_{\mathcal{\widetilde{F}}_{1}}*_{N}\mathcal{T}_{44}],\\
\label{XZEEyTTTw993}
&\mathcal{W}_{99}=[\mathcal{R}_{\mathcal{\widetilde{F}}_{1}}*_{N}(\mathcal{\widetilde{\widetilde{E}}}_{1}
                                 -\mathcal{\widetilde{H}}_{1}*_{N}\mathcal{K}_{44}*_{M}\mathcal{\widetilde{J}}_{1})*_{M}\mathcal{\widetilde{G}}_{1}^{\dagger}+
\mathcal{\widetilde{F}}_{1}*_{N}\mathcal{T}_{33}+\mathcal{T}_{55}*_{M}\mathcal{R}_{\mathcal{\widetilde{G}}_{1}}]*_{M}\begin{bmatrix}\mathcal{I} \\ 0\end{bmatrix},\\
\label{XZEEyTTT9q93}
&\mathcal{K}_{66}=[\mathcal{R}_{\mathcal{\widetilde{F}}_{1}}*_{N}(\mathcal{\widetilde{\widetilde{E}}}_{1}
                                 -\mathcal{\widetilde{H}}_{1}*_{N}\mathcal{K}_{44}*_{M}\mathcal{\widetilde{J}}_{1})*_{M}\mathcal{\widetilde{G}}_{1}^{\dagger}+
\mathcal{\widetilde{F}}_{1}*_{N}\mathcal{T}_{33}+\mathcal{T}_{55}*_{M}\mathcal{R}_{\mathcal{\widetilde{G}}_{1}}]*_{M}\begin{bmatrix}0 \\ \mathcal{I}\end{bmatrix},
\end{align}
\begin{align}
\label{XZEEyTTTw1993}
&\mathcal{K}_{77}=\begin{bmatrix}\mathcal{I}&
0\end{bmatrix}*_{N}[\mathcal{\widetilde{F}}_{2}*_{N}(\mathcal{\widetilde{\widetilde{E}}}_{2}
                                 -\mathcal{\widetilde{H}}_{2}*_{N}\mathcal{K}_{44}*_{M}\mathcal{\widetilde{J}}_{2})+
\mathcal{T}_{66}*_{M}\mathcal{\widetilde{B}}_{2}+\mathcal{L}_{\mathcal{\widetilde{F}}_{2}}*_{N}\mathcal{T}_{77}],\\
\label{ppp}
&\mathcal{T}_{11}=\begin{bmatrix}0&
\mathcal{I}\end{bmatrix}*_{N}[\mathcal{\widetilde{F}}_{2}*_{N}(\mathcal{\widetilde{\widetilde{E}}}_{2}
                                 -\mathcal{\widetilde{H}}_{2}*_{N}\mathcal{K}_{44}*_{M}\mathcal{\widetilde{J}}_{2})+
\mathcal{T}_{66}*_{M}\mathcal{\widetilde{B}}_{2}+\mathcal{L}_{\mathcal{\widetilde{F}}_{2}}*_{N}\mathcal{T}_{77}],\\
\label{kkk}
&\mathcal{K}_{88}=[\mathcal{R}_{\mathcal{\widetilde{F}}_{2}}*_{N}(\mathcal{\widetilde{\widetilde{E}}}_{2}
                                 -\mathcal{\widetilde{H}}_{2}*_{N}\mathcal{K}_{44}*_{M}\mathcal{\widetilde{J}}_{2})*_{M}\mathcal{\widetilde{G}}_{2}^{\dagger}+
\mathcal{\widetilde{F}}_{2}*_{N}\mathcal{T}_{66}+\mathcal{T}_{88}*_{M}\mathcal{R}_{\mathcal{\widetilde{G}}_{2}}]*_{M}\begin{bmatrix}\mathcal{I} \\ 0\end{bmatrix},\\
\label{kklk}
&\mathcal{T}_{22}=[\mathcal{R}_{\mathcal{\widetilde{F}}_{2}}*_{N}(\mathcal{\widetilde{\widetilde{E}}}_{2}
                                 -\mathcal{\widetilde{H}}_{2}*_{N}\mathcal{K}_{44}*_{M}\mathcal{\widetilde{J}}_{2})*_{M}\mathcal{\widetilde{G}}_{2}^{\dagger}+
\mathcal{\widetilde{F}}_{2}*_{N}\mathcal{T}_{66}+\mathcal{T}_{88}*_{M}\mathcal{R}_{\mathcal{\widetilde{G}}_{2}}]*_{M}\begin{bmatrix}0 \\ \mathcal{I}\end{bmatrix},\\
\label{sss1}
&\mathcal{K}_{44}=\mathcal{\widetilde{H}}_{11}^{\dagger}*_{N}\mathcal{\widetilde{\widetilde{E}}}_{11}
*_{M}\mathcal{\widetilde{J}}_{11}^{\dagger}+\mathcal{L}_{\mathcal{\widetilde{H}}_{11}}*_{N}\mathcal{\acute{W}}_{2}
+\mathcal{\acute{W}}_{3}*_{M}\mathcal{R}_{\mathcal{\widetilde{J}}_{11}},\\
\label{sss2}
&or\ \mathcal{K}_{44}=\mathcal{\widetilde{H}}_{22}^{\dagger}*_{N}\mathcal{\widetilde{\widetilde{E}}}_{22}
*_{M}\mathcal{\widetilde{J}}_{22}^{\dagger}+\mathcal{L}_{\mathcal{\widetilde{H}}_{22}}*_{N}\mathcal{\acute{W}}_{4}
+\mathcal{\acute{W}}_{5}*_{M}\mathcal{R}_{\mathcal{\widetilde{J}}_{22}},\\
\label{sss22}
&\mathcal{T}_{33}=\begin{bmatrix}\mathcal{I}&
0\end{bmatrix}*_{N}[\mathcal{\widetilde{A}}*_{N}\mathcal{\widetilde{E}}+
\mathcal{\acute{W}}_{6}*_{M}\mathcal{\widetilde{B}}+\mathcal{L}_{\mathcal{\widetilde{A}}}*_{N}\mathcal{\acute{W}}_{7}],\\
\label{sss322}
&\mathcal{T}_{55}=\begin{bmatrix}0&
\mathcal{I}\end{bmatrix}*_{N}[\mathcal{\widetilde{A}}*_{N}\mathcal{\widetilde{E}}+
\mathcal{\acute{W}}_{6}*_{M}\mathcal{\widetilde{B}}+\mathcal{L}_{\mathcal{\widetilde{A}}}*_{N}\mathcal{\acute{W}}_{7}],\\
\label{sss3H22}
&\mathcal{T}_{44}=[\mathcal{R}_{\mathcal{\widetilde{A}}}*_{N}\mathcal{\widetilde{E}}*_{M}\mathcal{\widetilde{B}}^{\dagger}+
\mathcal{\widetilde{A}}*_{N}\mathcal{\acute{W}}_{6}+\mathcal{\acute{W}}_{8}*_{M}\mathcal{R}_{\mathcal{\widetilde{B}}}]*_{M}\begin{bmatrix}\mathcal{I} \\ 0\end{bmatrix},\\
\label{sss3HY22}
&\mathcal{T}_{66}=[\mathcal{R}_{\mathcal{\widetilde{A}}}*_{N}\mathcal{\widetilde{E}}*_{M}\mathcal{\widetilde{B}}^{\dagger}+
\mathcal{\widetilde{A}}*_{N}\mathcal{\acute{W}}_{6}+\mathcal{\acute{W}}_{8}*_{M}\mathcal{R}_{\mathcal{\widetilde{B}}}]*_{M}\begin{bmatrix}0 \\ \mathcal{I}\end{bmatrix}.
\end{align}
\end{subequations}
Where $\mathcal{U}_{ji}$, $\mathcal{W}_{kk}$, $\mathcal{V}_{kk}$, $\acute{\mathcal{W}}_{k}$, $\mathcal{P}_{ii}$, $\mathcal{Q}_{ii}$, $\mathcal{K}_{ii}$ and $\mathcal{T}_{ss}$,  $(i=\overline{1,3},\ j=\overline{1,5},\ k=\overline{1,6},\ s=\overline{3,8})$ are arbitrary quaternion tensors with appropriate sizes.
\end{theorem}
\begin{proof}
The system \eqref{1.4aa} can divide to the following two-sided Sylvester-like  tensor equations:
\begin{equation}
\label{k11a}
\mathcal{A}_{1}*_{N}\mathcal{X}_{1}*_{M}\mathcal{B}_{1}+\mathcal{C}_{1}*_{N}\mathcal{Y}_{1}*_{M}\mathcal{D}_{1}\\
+\mathcal{C}_{1}*_{N}(\mathcal{F}_{1}*_{N}\mathcal{Z}_{1}*_{M}\mathcal{G}_{1}
+\mathcal{H}_{1}*_{N}\mathcal{Z}_{2}*_{M}\mathcal{J}_{1})*_{M}\mathcal{B}_{1}=\mathcal{E}_{1},
\end{equation}
\begin{equation}
\label{k12a}
\mathcal{A}_{2}*_{N}\mathcal{X}_{2}*_{M}\mathcal{B}_{2}+\mathcal{C}_{2}*_{N}\mathcal{Y}_{2}*_{M}\mathcal{D}_{2}\\
+\mathcal{C}_{2}*_{N}(\mathcal{F}_{2}*_{N}\mathcal{Z}_{2}*_{M}\mathcal{G}_{2}
+\mathcal{H}_{2}*_{N}\mathcal{Z}_{3}*_{M}\mathcal{J}_{2})*_{M}\mathcal{B}_{2}=\mathcal{E}_{2},
\end{equation}
\begin{equation}
\label{k13a}
\mathcal{A}_{3}*_{N}\mathcal{X}_{3}*_{M}\mathcal{B}_{3}+\mathcal{C}_{3}*_{N}\mathcal{Y}_{3}*_{M}\mathcal{D}_{3}\\
+\mathcal{C}_{3}*_{N}(\mathcal{F}_{3}*_{N}\mathcal{Z}_{3}*_{M}\mathcal{G}_{3}
+\mathcal{H}_{3}*_{N}\mathcal{Z}_{4}*_{M}\mathcal{J}_{3})*_{M}\mathcal{B}_{3}=\mathcal{E}_{3},
\end{equation}
\begin{equation}
\label{k14}
\mathcal{F}_{4}*_{N}\mathcal{Z}_{1}*_{M}\mathcal{G}_{4}=\mathcal{E}_{4},\
\end{equation}
\begin{equation}
\label{k15}
\mathcal{H}_{4}*_{N}\mathcal{Z}_{4}*_{M}\mathcal{J}_{4}=\mathcal{E}_{5}.
\end{equation}
The main idea is to implement the conditions of consistency to enable this group to have a solution, and hence, we  establish an expression of this  solution. Applying $Lemma$ $\ref{lma 2.3}$, we have that the Sylvester-like tensor equation \eqref{k11a} is consistent if and only if 
\begin{align}
\begin{gathered}
\mathcal{R}_{\mathcal{M}_{1}}*_{N}\mathcal{R}_{\mathcal{A}_{1}}*_{N}\mathcal{E}_{1}=0,\ \mathcal{E}_{1}*_{M}\mathcal{L}_{\mathcal{B}_{1}}*_{M}\mathcal{L}_{\mathcal{N}_{1}}=0,\
\mathcal{R}_{\mathcal{C}_{1}}*_{N}\mathcal{E}_{1}*_{M}\mathcal{L}_{\mathcal{B}_{1}}=0\\
\mathcal{R}_{\mathcal{\widehat{M}}_{1}}*_{N}\mathcal{R}_{\mathcal{\widehat{A}}_{1}}*_{N}\mathcal{\widehat{E}}_{1}=0,\ \mathcal{\widehat{E}}_{1}*_{M}\mathcal{L}_{\mathcal{\widehat{B}}_{1}}*_{M}\mathcal{L}_{\mathcal{\widehat{N}}_{1}}=0,\\ \mathcal{R}_{\mathcal{\widehat{A}}_{1}}*_{N}\mathcal{\widehat{E}}_{1}*_{M}\mathcal{L}_{\mathcal{\widehat{D}}_{1}}=0,\
\mathcal{R}_{\mathcal{\widehat{C}}_{1}}*_{N}\mathcal{\widehat{E}}_{1}*_{M}\mathcal{L}_{\mathcal{\widehat{B}}_{1}}=0.
\end{gathered}
\end{align}
In that case, the general solution can be expressed as
\begin{subequations}
\begin{align}
&\begin{array}{l}
 \mathcal{X}_{1}=\mathcal{A}_{1}^{\dagger}*_{N}\mathcal{\grave{E}}_{1}*_{M}\mathcal{B}_{1}^{\dagger}
-\mathcal{A}_{1}^{\dagger}*_{N}\mathcal{C}_{1}*_{N}\mathcal{M}_{1}^{\dagger}*_{N}\mathcal{\grave{E}}_{1}*_{M}\mathcal{B}_{1}^{\dagger}
-\mathcal{A}_{1}^{\dagger}*_{N}\mathcal{S}_{1}*_{N}\mathcal{C}_{1}^{\dagger}*_{N}\mathcal{\grave{E}}_{1}\\
\ \ \ \ \ *_{M}\mathcal{N}_{1}^{\dagger}*_{M}\mathcal{D}_{1}*_{M}\mathcal{B}_{1}^{\dagger}-\mathcal{A}_{1}^{\dagger}*_{N}\mathcal{S}_{1}
*_{N}\mathcal{U}_{21}*_{M}\mathcal{R}_{\mathcal{N}_{1}}*_{M}\mathcal{D}_{1}*_{M}\mathcal{B}_{1}^{\dagger}
+\mathcal{L}_{\mathcal{A}_{1}}*_{N}\mathcal{U}_{41}\\
\ \ \ \ \ +\mathcal{U}_{51}*_{M}\mathcal{R}_{\mathcal{B}_{1}},
 \end{array}\\
&\begin{array}{l}
 \mathcal{Y}_{1}=\mathcal{M}_{1}^{\dagger}*_{N}\mathcal{\grave{E}}_{1}*_{M}\mathcal{D}_{1}^{\dagger}
 +\mathcal{S}_{1}^{\dagger}*_{N}\mathcal{S}_{1}*_{N}\mathcal{C}_{1}^{\dagger}
*_{N}\mathcal{\grave{E}}_{1}*_{M}\mathcal{N}_{1}^{\dagger}+\mathcal{L}_{\mathcal{M}_{1}}*_{N}\mathcal{L}_{\mathcal{S}_{1}}*_{N}\mathcal{U}_{11}\\
\ \ \ \ \ \ +\mathcal{L}_{\mathcal{M}_{1}}*_{N}\mathcal{U}_{21}*_{M}\mathcal{R}_{\mathcal{N}_{1}}+\mathcal{U}_{31}*_{M}\mathcal{R}_{\mathcal{D}_{1}}.
 \end{array}\\
 \label{x11}
&\begin{array}{l}
 \mathcal{Z}_{1}=\mathcal{\widehat{A}}_{1}^{\dagger}*_{N}\mathcal{\widehat{E}}_{1}*_{M}\mathcal{\widehat{B}}_{1}^{\dagger}
-\mathcal{\widehat{A}}_{1}^{\dagger}*_{N}\mathcal{\widehat{C}}_{1}*_{N}\mathcal{\widehat{M}}_{1}^{\dagger}*_{N}\mathcal{\widehat{E}}_{1}
*_{M}\mathcal{\widehat{B}}_{1}^{\dagger}
-\mathcal{\widehat{A}}_{1}^{\dagger}*_{N}\mathcal{\widehat{S}}_{1}*_{N}\mathcal{\widehat{C}}_{1}^{\dagger}*_{N}\mathcal{\widehat{E}}_{1}\\
\ \ \ \ \ *_{M}\mathcal{\widehat{N}}_{1}^{\dagger}*_{M}\mathcal{\widehat{D}}_{1}*_{M}\mathcal{\widehat{B}}_{1}^{\dagger}-\mathcal{\widehat{A}}_{1}^{\dagger}
*_{N}\mathcal{\widehat{S}}_{1}
*_{N}\mathcal{\widehat{U}}_{2}*_{M}\mathcal{R}_{\mathcal{\widehat{N}}_{1}}*_{M}\mathcal{\widehat{D}}_{1}*_{M}\mathcal{\widehat{B}}_{1}^{\dagger}
+\mathcal{L}_{\mathcal{\widehat{A}}_{1}}*_{N}\mathcal{\widehat{U}}_{4}\\
\ \ \ \ \ +\mathcal{\widehat{U}}_{5}*_{M}\mathcal{R}_{\mathcal{\widehat{B}}_{1}},
 \end{array}  \\
&\begin{array}{l}
\label{x221}
 \mathcal{Z}_{2}=\mathcal{\widehat{M}}_{1}^{\dagger}*_{N}\mathcal{\widehat{E}}_{1}*_{M}\mathcal{\widehat{D}}_{1}^{\dagger}
 +\mathcal{\widehat{S}}_{1}^{\dagger}*_{N}\mathcal{\widehat{S}}_{1}*_{N}\mathcal{\widehat{C}}_{1}^{\dagger}
*_{N}\mathcal{\widehat{E}}*_{M}\mathcal{\widehat{N}}_{1}^{\dagger}+\mathcal{L}_{\mathcal{\widehat{M}}_{1}}
*_{N}\mathcal{L}_{\mathcal{\widehat{S}}_{1}}*_{N}\mathcal{\widehat{U}}_{1}\\
\ \ \ \ \ \ +\mathcal{L}_{\mathcal{\widehat{M}}_{1}}*_{N}\mathcal{\widehat{U}}_{2}*_{M}\mathcal{R}_{\mathcal{\widehat{N}}_{1}}
+\mathcal{\widehat{U}}_{3}*_{M}\mathcal{R}_{\mathcal{\widehat{D}}_{1}},
 \end{array}
\end{align}
\end{subequations}
where $\mathcal{\widehat{A}}_{1}$, $\mathcal{\widehat{B}}_{1}$, $\mathcal{\widehat{C}}_{1}$, $\mathcal{\widehat{D}}_{1}$, $\mathcal{\widehat{E}}_{1}$, $\mathcal{\widehat{M}}_{1}$, $\mathcal{\widehat{N}}_{1}$ and $\mathcal{\widehat{S}}_{1}$ given by \eqref{abc22}-\eqref{abc44} whenever $i=1$.
It can follow the same technique to determine the consistency conditions and the general solution to the Sylvester-like quaternion tensor equation \eqref{k12a}. So, we have that Eq.\eqref{k12a} is solvable if and only if the conditions \eqref{b11}-\eqref{b331} satisfy whenever $i=2$. In this case, the quaternion tensors $\mathcal{X}_{2}$ and $\mathcal{Y}_{2}$ can be given by \eqref{sssa222p1}-\eqref{sssb} whenever $i=2$ and
\begin{subequations}
\begin{align}
&\begin{array}{l}
 \label{x331}
 \mathcal{Z}_{2}=\mathcal{\widehat{A}}_{2}^{\dagger}*_{N}\mathcal{\widehat{E}}_{2}*_{M}\mathcal{\widehat{B}}_{2}^{\dagger}
-\mathcal{\widehat{A}}_{2}^{\dagger}*_{N}\mathcal{\widehat{C}}_{2}*_{N}\mathcal{\widehat{M}}_{2}^{\dagger}*_{N}\mathcal{\widehat{E}}_{2}
*_{M}\mathcal{\widehat{B}}_{2}^{\dagger}
-\mathcal{\widehat{A}}_{2}^{\dagger}*_{N}\mathcal{\widehat{S}}_{2}*_{N}\mathcal{\widehat{C}}_{2}^{\dagger}*_{N}\mathcal{\widehat{E}}_{2}\\
\ \ \ \ \ *_{M}\mathcal{\widehat{N}}_{2}^{\dagger}*_{M}\mathcal{\widehat{D}}_{2}*_{M}\mathcal{\widehat{B}}_{2}^{\dagger}-\mathcal{\widehat{A}}_{2}^{\dagger}
*_{N}\mathcal{\widehat{S}}_{2}
*_{N}\mathcal{\widehat{V}}_{2}*_{M}\mathcal{R}_{\mathcal{\widehat{N}}_{2}}*_{M}\mathcal{\widehat{D}}_{2}*_{M}\mathcal{\widehat{B}}_{2}^{\dagger}
+\mathcal{L}_{\mathcal{\widehat{A}}_{2}}*_{N}\mathcal{\widehat{V}}_{4}\\
\ \ \ \ \ +\mathcal{\widehat{V}}_{5}*_{M}\mathcal{R}_{\mathcal{\widehat{B}}_{2}},
 \end{array}  \\
&\begin{array}{l}
\label{x441}
 \mathcal{Z}_{3}=\mathcal{\widehat{M}}_{2}^{\dagger}*_{N}\mathcal{\widehat{E}}_{2}*_{M}\mathcal{\widehat{D}}_{2}^{\dagger}
 +\mathcal{\widehat{S}}_{2}^{\dagger}*_{N}\mathcal{\widehat{S}}_{2}*_{N}\mathcal{\widehat{C}}_{2}^{\dagger}
*_{N}\mathcal{\widehat{E}}_{2}*_{M}\mathcal{\widehat{N}}_{2}^{\dagger}+\mathcal{L}_{\mathcal{\widehat{M}}_{2}}
*_{N}\mathcal{L}_{\mathcal{\widehat{S}}_{2}}*_{N}\mathcal{\widehat{V}}_{1}\\
\ \ \ \ \ \ +\mathcal{L}_{\mathcal{\widehat{M}}_{2}}*_{N}\mathcal{\widehat{V}}_{2}*_{M}\mathcal{R}_{\mathcal{\widehat{N}}_{2}}
+\mathcal{\widehat{V}}_{3}*_{M}\mathcal{R}_{\mathcal{\widehat{D}}_{2}},
 \end{array}
\end{align}
\end{subequations}
where $\mathcal{\widehat{A}}_{2}$, $\mathcal{\widehat{B}}_{2}$, $\mathcal{\widehat{C}}_{2}$, $\mathcal{\widehat{D}}_{2}$, $\mathcal{\widehat{E}}_{2}$, $\mathcal{\widehat{M}}_{2}$, $\mathcal{\widehat{N}}_{2}$ and $\mathcal{\widehat{S}}_{2}$ given by \eqref{abc22}-\eqref{abc44} whenever $i=2$.\par
Similarly, we can provide that E.q. \eqref{k13a} is solvable if and only if
the conditions \eqref{b11}-\eqref{b331} are satisfying whenever $i=3$.
In this case, the quaternion tensors $\mathcal{X}_{3}$ and $\mathcal{Y}_{3}$ can be given by \eqref{sssa222p1}-\eqref{sssb}, whenever $i=3$ and
\begin{subequations}
\begin{align}
 &\begin{array}{l}
 \label{x551}
 \mathcal{Z}_{3}=\mathcal{\widehat{A}}_{3}^{\dagger}*_{N}\mathcal{\widehat{E}}_{3}*_{M}\mathcal{\widehat{B}}_{3}^{\dagger}
-\mathcal{\widehat{A}}_{3}^{\dagger}*_{N}\mathcal{\widehat{C}}_{3}*_{N}\mathcal{\widehat{M}}_{3}^{\dagger}*_{N}\mathcal{\widehat{E}}_{3}
*_{M}\mathcal{\widehat{B}}_{3}^{\dagger}
-\mathcal{\widehat{A}}_{3}^{\dagger}*_{N}\mathcal{\widehat{S}}_{3}*_{N}\mathcal{\widehat{C}}_{3}^{\dagger}*_{N}\mathcal{\widehat{E}}_{3}\\
\ \ \ \ \ *_{M}\mathcal{\widehat{N}}_{3}^{\dagger}*_{M}\mathcal{\widehat{D}}_{3}*_{M}\mathcal{\widehat{B}}_{3}^{\dagger}-\mathcal{\widehat{A}}_{3}^{\dagger}
*_{N}\mathcal{\widehat{S}}_{3}
*_{N}\mathcal{\widehat{K}}_{2}*_{M}\mathcal{R}_{\mathcal{\widehat{N}}_{3}}*_{M}\mathcal{\widehat{D}}_{3}*_{M}\mathcal{\widehat{B}}_{3}^{\dagger}
+\mathcal{L}_{\mathcal{\widehat{A}}_{3}}*_{N}\mathcal{\widehat{K}}_{4}\\
\ \ \ \ \ +\mathcal{\widehat{K}}_{5}*_{M}\mathcal{R}_{\mathcal{\widehat{B}}_{3}},
 \end{array}  \\
&\begin{array}{l}
\label{x661}
 \mathcal{Z}_{4}=\mathcal{\widehat{M}}_{3}^{\dagger}*_{N}\mathcal{\widehat{E}}_{3}*_{M}\mathcal{\widehat{D}}_{3}^{\dagger}
 +\mathcal{\widehat{S}}_{3}^{\dagger}*_{N}\mathcal{\widehat{S}}_{3}*_{N}\mathcal{\widehat{C}}_{3}^{\dagger}
*_{N}\mathcal{\widehat{E}}_{3}*_{M}\mathcal{\widehat{N}}_{3}^{\dagger}+\mathcal{L}_{\mathcal{\widehat{M}}_{3}}
*_{N}\mathcal{L}_{\mathcal{\widehat{S}}_{3}}*_{N}\mathcal{\widehat{K}}_{1}\\
\ \ \ \ \ \ +\mathcal{L}_{\mathcal{\widehat{M}}_{3}}*_{N}\mathcal{\widehat{K}}_{2}*_{M}\mathcal{R}_{\mathcal{\widehat{N}}_{3}}
+\mathcal{\widehat{K}}_{3}*_{M}\mathcal{R}_{\mathcal{\widehat{D}}_{3}},
 \end{array}
\end{align}
\end{subequations}
where $\mathcal{\widehat{A}}_{3}$, $\mathcal{\widehat{B}}_{3}$, $\mathcal{\widehat{C}}_{3}$, $\mathcal{\widehat{D}}_{3}$, $\mathcal{\widehat{E}}_{3}$, $\mathcal{\widehat{M}}_{3}$, $\mathcal{\widehat{N}}_{3}$ and $\mathcal{\widehat{S}}_{3}$ given by \eqref{abc22}-\eqref{abc44}, whenever $i=3$.\par
It follows from $Lemma$ $\ref{lma 2.3}$ that the necessary and sufficient conditions
for the Sylvester-like quaternion tensor equation \eqref{k14} and \eqref{k15} to be consistent are given by \eqref{b441}, respectively. Consequently, the solutions to these two equations are expressed as
\begin{subequations}
\begin{align}
\label{x77}
 &\mathcal{Z}_{1}=\mathcal{F}_{4}^{\dagger}*_{N}\mathcal{E}_{4}*_{M}\mathcal{G}_{4}^{\dagger}
+\mathcal{L}_{\mathcal{F}_{4}}*_{N}\mathcal{W}_{1}
+\mathcal{W}_{2}*_{M}\mathcal{R}_{\mathcal{G}_{4}},\\
\label{x881}
&\mathcal{Z}_{4}=\mathcal{H}_{4}^{\dagger}*_{N}\mathcal{E}_{5}*_{M}\mathcal{J}_{4}^{\dagger}
+\mathcal{L}_{\mathcal{H}_{4}}*_{N}\mathcal{\acute{W}}_{1}
+\mathcal{\acute{W}}_{2}*_{M}\mathcal{R}_{\mathcal{J}_{4}},
\end{align}
\end{subequations}\par
Let $Z_{1}$ in \eqref{x11} be equal to $Z_{1}$ in \eqref{x77}, and $Z_{4}$ in \eqref{x661} be equal to $Z_{4}$ in \eqref{x881}. Then we have the following equations:
\begin{align}
\label{x88}
&\mathcal{A}_{11}*_{N}\begin{bmatrix}
                                \mathcal{W}_{1} \\
                                 \mathcal{\widehat{U}}_{4}\end{bmatrix}+\begin{bmatrix}
                                 \mathcal{W}_{2}&\mathcal{\widehat{U}}_{5}\end{bmatrix}*_{M} \mathcal{D}_{11}=\mathcal{E}_{11}-\mathcal{\widehat{A}}_{11}*_{N}\mathcal{\widehat{U}}_{2}*_{M}\mathcal{\widehat{B}}_{11},\\
\label{x99}
&\mathcal{A}_{22}*_{N}\begin{bmatrix}
                                \mathcal{\grave{W}}_{1} \\
                                 \mathcal{\widehat{K}}_{1}\end{bmatrix}+\begin{bmatrix}
                                 \mathcal{W}_{3}&\mathcal{\widehat{K}}_{3}\end{bmatrix}*_{M} \mathcal{D}_{22}=\mathcal{E}_{22}-\mathcal{\widehat{A}}_{22}*_{N}\mathcal{\widehat{K}}_{2}*_{M}\mathcal{\widehat{B}}_{22},
\end{align}
Apply $Lemma$ $\ref{lma 2.3},$ to Eq.\eqref{x88}, we have that it is solvable if and only if there exists a quaternion tensor $\mathcal{\widehat{U}}_{2}$ satisfies
\begin{align}
\label{x9911}
&\mathcal{\widehat{\widehat{A}}}_{11}*_{N}\mathcal{\widehat{U}}_{2}*_{M}\mathcal{\widehat{\widehat{B}}}_{11}
=\mathcal{\widehat{\widehat{E}}}_{11}.
\end{align}
In that case, the general solution can be express as
\begin{align}
\label{x992aa}
&\begin{bmatrix}
\mathcal{W}_{1} \\
\mathcal{\widehat{U}}_{4}\end{bmatrix}=\mathcal{A}_{11}^{\dagger}*_{N}
(\mathcal{E}_{11}-\mathcal{\widehat{A}}_{11}*_{N}\mathcal{\widehat{U}}_{2}*_{M}\mathcal{\widehat{B}}_{11})
-\mathcal{V}_{11}*_{M}\mathcal{D}_{11}+\mathcal{L}_{\mathcal{A}_{11}}*_{N}\mathcal{V}_{22},\\
\label{x993aa}
&\begin{bmatrix}\mathcal{W}_{2}&\mathcal{\widehat{U}}_{5}\end{bmatrix}=
\mathcal{R}_{\mathcal{A}_{11}}*_{N}
(\mathcal{E}_{11}-\mathcal{\widehat{A}}_{11}*_{N}\mathcal{\widehat{U}}_{2}*_{M}\mathcal{\widehat{B}}_{11})*_{M}\mathcal{D}_{11}^{\dagger}
+\mathcal{A}_{11}*_{N}\mathcal{V}_{11}+\mathcal{V}_{33}*_{M}\mathcal{R}_{\mathcal{D}_{11}}.
\end{align}By applying $Proposition$ \ref{lma 2.22} to equations \eqref{x992aa}-\eqref{x993aa}, we can find expressions for quaternion tensors $\mathcal{W}_{1}$, $\mathcal{\widehat{U}}_{4}$, $\mathcal{W}_{2}$, and $\mathcal{\widehat{U}}_{5}$ in \eqref{Xyqs1}-\eqref{Xyqs4}.
In the same way, we have that Eq.\eqref{x99} is solvable if and only if there exists a quaternion tensor $\mathcal{\widehat{K}}_{2}$ satisfies
\begin{align}
\label{x991}
&\mathcal{\widehat{\widehat{A}}}_{22}*_{N}\mathcal{\widehat{K}}_{2}*_{M}\mathcal{\widehat{\widehat{B}}}_{22}
=\mathcal{\widehat{\widehat{E}}}_{22}.
\end{align}
In that case, the general solution can be express as
\begin{align}
\label{x9921}
&\begin{bmatrix}
\mathcal{\grave{W}}_{1} \\
\mathcal{\widehat{K}}_{1}\end{bmatrix}=\mathcal{A}_{22}^{\dagger}*_{N}
(\mathcal{E}_{22}-\mathcal{\widehat{A}}_{22}*_{N}\mathcal{\widehat{K}}_{2}*_{M}\mathcal{\widehat{B}}_{22})
-\mathcal{V}_{44}*_{M}\mathcal{D}_{22}+\mathcal{L}_{\mathcal{A}_{22}}*_{N}\mathcal{V}_{55},\\
\label{x9931}
&\begin{bmatrix}\mathcal{W}_{3}&\mathcal{\widehat{K}}_{3}\end{bmatrix}=
\mathcal{R}_{\mathcal{A}_{22}}*_{N}
(\mathcal{E}_{22}-\mathcal{\widehat{A}}_{22}*_{N}\mathcal{\widehat{K}}_{2}*_{M}\mathcal{\widehat{B}}_{22})*_{M}\mathcal{D}_{22}^{\dagger}
+\mathcal{A}_{22}*_{N}\mathcal{V}_{44}+\mathcal{V}_{66}*_{M}\mathcal{R}_{\mathcal{D}_{22}}.
\end{align}
 It can be utilized $Proposition$ \ref{lma 2.22} to equations \eqref{x9921}-\eqref{x9931}, we can get quaternion tensors $\mathcal{\acute{W}}_{1}$, $\mathcal{\widehat{K}}_{1}$, $\mathcal{\acute{W}}_{2}$ and $\mathcal{\widehat{K}}_{3}$ in \eqref{Xyqs5}-\eqref{Xyqs8}.
Meanwhile, the quaternion tensor equations \eqref{x9911} and \eqref{x991} are solvable if and only if the conditions \eqref{b5511} are satisfying, respectively, for $k=1,2$. In that case, the general solution can be written as
\begin{subequations}
\begin{align}
 \label{x9944}
&\widehat{\mathcal{U}}_{2}=\widehat{\widehat{\mathcal{A}}}_{11}^{\dagger}*_{N}\widehat{\widehat{\mathcal{E}}}_{11}
*_{M}\widehat{\widehat{\mathcal{B}}}_{11}^{\dagger}+\mathcal{L}_{\widehat{\widehat{\mathcal{A}}}_{11}}*_{N}\mathcal{V}_{77}
+\mathcal{V}_{88}*_{M}\mathcal{R}_{\widehat{\widehat{\mathcal{B}}}_{11}},\\
\label{x9955}
&\widehat{\mathcal{K}}_{2}=\widehat{\widehat{\mathcal{A}}}_{22}^{\dagger}*_{N}\widehat{\widehat{\mathcal{E}}}_{22}
*_{M}\widehat{\widehat{\mathcal{B}}}_{22}^{\dagger}+\mathcal{L}_{\widehat{\widehat{\mathcal{A}}}_{22}}*_{N}\mathcal{V}_{99}
+\mathcal{W}_{11}*_{M}\mathcal{R}_{\widehat{\widehat{\mathcal{B}}}_{22}}.
\end{align}
\end{subequations}\par
Now, $Z_{2}$ in \eqref{x221} should be equal to $Z_{2}$ in \eqref{x331} and $Z_{3}$ in \eqref{x441} should be equal to $Z_{3}$ in \eqref{x551}, yields:
\begin{align}
\label{y11}
&\mathcal{\overline{A}}_{1}*_{N}\begin{bmatrix}
                                \mathcal{\widehat{U}}_{1} \\
                                 \mathcal{\widehat{V}}_{4}\end{bmatrix}+\begin{bmatrix}
                                 \mathcal{\widehat{U}}_{3}&\mathcal{\widehat{V}}_{5}\end{bmatrix}*_{M} \mathcal{\overline{B}}_{1}=-\mathcal{\overline{E}}_{1}
                                 +\mathcal{\overline{F}}_{1}*_{N}\mathcal{\widehat{V}}_{2}*_{M}\mathcal{\overline{G}}_{1}
                                 +\mathcal{\overline{H}}_{1}*_{N}\mathcal{\widehat{U}}_{2}*_{M}\mathcal{\overline{J}}_{1},\\
\label{y22}
&\mathcal{\overline{A}}_{2}*_{N}\begin{bmatrix}
                                \mathcal{\widehat{V}}_{1} \\
                                 \mathcal{\widehat{K}}_{4}\end{bmatrix}+\begin{bmatrix}
                                 \mathcal{\widehat{V}}_{3}&\mathcal{\widehat{K}}_{5}\end{bmatrix}*_{M} \mathcal{\overline{B}}_{2}=-\mathcal{\overline{E}}_{2}
                                 +\mathcal{\overline{F}}_{2}*_{N}\mathcal{\widehat{K}}_{2}*_{M}\mathcal{\overline{G}}_{2}
                                 +\mathcal{\overline{H}}_{2}*_{N}\mathcal{\widehat{V}}_{2}*_{M}\mathcal{\overline{J}}_{2}.
\end{align}
The  system of tensor equations which consists of the two equations \eqref{y11} and \eqref {y22} is consistent if and only if there exists quaternion tensors $\mathcal{\widehat{V}}_{2}$, $\mathcal{\widehat{U}}_{2}$, and $\mathcal{\widehat{K}}_{2}$ satisfy the following system:
\begin{align}
\label{y111}
&\mathcal{\overline{F}}_{11}*_{N}\mathcal{\widehat{V}}_{2}*_{M}\mathcal{\overline{G}}_{11}
+\mathcal{\overline{H}}_{11}*_{N}\mathcal{\widehat{U}}_{2}*_{M}\mathcal{\overline{J}}_{11}
=\mathcal{\overline{E}}_{11},\\
\label{y222}
&\mathcal{\overline{F}}_{22}*_{N}\mathcal{\widehat{K}}_{2}*_{M}\mathcal{\overline{G}}_{22}
+\mathcal{\overline{H}}_{22}*_{N}\mathcal{\widehat{V}}_{2}*_{M}\mathcal{\overline{J}}_{22}
=\mathcal{\overline{E}}_{22},
\end{align}
In utilizing $Lemma$ \ref{lma 2.3}, we have that the general solution to quaternion tensor equations \eqref{y11}-\eqref{y22} can be given by \eqref{Xya991}-\eqref{Xya998}. The quaternion system of tensor equations \eqref{y111}-\eqref{y222} is solvable if and only if Eq. \eqref{y111} is solvable, Eq. \eqref{y222} is solvable, and the quaternion tensor $\mathcal{\widehat{V}}_{2}$ in \eqref{y111} coincide with $\mathcal{\widehat{V}}_{2}$ in \eqref{y222}. So, Eq. \eqref{y111} and Eq. \eqref{y222} are solvable if and only if the conditions \eqref{b6611}-\eqref{b7711} are satisfying, respectively, for $k=1,2$. In this case, the general solution can be express as
\begin{subequations}
\begin{align}
\label{XZy993}
&\begin{array}{l}
 \mathcal{\widehat{V}}_{2}=\mathcal{\overline{F}}_{11}^{\dagger}*_{N}\mathcal{\overline{E}}_{11}*_{M}\mathcal{\overline{G}}_{11}^{\dagger}
-\mathcal{\overline{F}}_{11}^{\dagger}*_{N}\mathcal{\overline{H}}_{11}*_{N}\mathcal{\overline{M}}_{11}^{\dagger}*_{N}\mathcal{\overline{E}}_{11}
*_{M}\mathcal{\overline{G}}_{11}^{\dagger}
-\mathcal{\overline{F}}_{11}^{\dagger}*_{N}\mathcal{\overline{S}}_{11}*_{N}\mathcal{\overline{H}}_{11}^{\dagger}\\
\ \ \ \ \ *_{N}\mathcal{\overline{E}}_{11}*_{M}\mathcal{\overline{N}}_{11}^{\dagger}*_{M}\mathcal{\overline{J}}_{11}*_{M}\mathcal{\overline{G}}_{11}^{\dagger}
-\mathcal{\overline{F}}_{11}^{\dagger}*_{N}\mathcal{\overline{S}}_{11}
*_{N}\mathcal{P}_{44}*_{M}\mathcal{R}_{\mathcal{\overline{N}}_{11}}*_{M}\mathcal{\overline{J}}_{11}*_{M}\mathcal{\overline{G}}_{11}^{\dagger}
\\
\ \ \ \ \ +\mathcal{L}_{\mathcal{\overline{F}}_{11}}*_{N}\mathcal{P}_{55}+\mathcal{P}_{66}*_{M}\mathcal{R}_{\mathcal{\overline{G}}_{11}},
 \end{array}\\
\label{XZy994}
&\begin{array}{l}
 \mathcal{\widehat{U}}_{2}=\mathcal{\overline{M}}_{11}^{\dagger}*_{N}\mathcal{\overline{E}}_{11}*_{M}\mathcal{\overline{J}}_{11}^{\dagger}
 +\mathcal{\overline{S}}_{11}^{\dagger}*_{N}\mathcal{\overline{S}}_{11}*_{N}\mathcal{\overline{H}}_{11}^{\dagger}
*_{N}\mathcal{\overline{E}}_{11}*_{M}\mathcal{\overline{N}}_{11}^{\dagger}+\mathcal{L}_{\mathcal{\overline{M}}_{11}}
\\
\ \ \ \ \ \ *_{N}\mathcal{L}_{\mathcal{\overline{S}}_{11}}*_{N}\mathcal{Q}_{44}+\mathcal{L}_{\mathcal{\overline{M}}_{11}}*_{N}\mathcal{P}_{44}*_{M}\mathcal{R}_{\mathcal{\overline{N}}_{11}}+\mathcal{Q}_{66}
*_{M}\mathcal{R}_{\mathcal{\overline{J}}_{11}},
 \end{array}
\end{align}
\begin{align}
\label{XZy995}
&\begin{array}{l}
 \mathcal{\widehat{K}}_{2}=\mathcal{\overline{F}}_{22}^{\dagger}*_{N}\mathcal{\overline{E}}_{22}*_{M}\mathcal{\overline{G}}_{22}^{\dagger}
-\mathcal{\overline{F}}_{22}^{\dagger}*_{N}\mathcal{\overline{H}}_{22}*_{N}\mathcal{\overline{M}}_{22}^{\dagger}*_{N}\mathcal{\overline{E}}_{22}
*_{M}\mathcal{\overline{G}}_{22}^{\dagger}
-\mathcal{\overline{F}}_{22}^{\dagger}*_{N}\mathcal{\overline{S}}_{22}*_{N}\mathcal{\overline{H}}_{22}^{\dagger}\\
\ \ \ \ \ *_{N}\mathcal{\overline{E}}_{22}*_{M}\mathcal{\overline{N}}_{22}^{\dagger}*_{M}\mathcal{\overline{J}}_{22}*_{M}\mathcal{\overline{G}}_{22}^{\dagger}
-\mathcal{\overline{F}}_{22}^{\dagger}*_{N}\mathcal{\overline{S}}_{22}
*_{N}\mathcal{Q}_{55}*_{M}\mathcal{R}_{\mathcal{\overline{N}}_{22}}*_{M}\mathcal{\overline{J}}_{22}*_{M}\mathcal{\overline{G}}_{22}^{\dagger}
\\
\ \ \ \ \ +\mathcal{L}_{\mathcal{\overline{F}}_{22}}*_{N}\mathcal{P}_{77}+\mathcal{P}_{88}*_{M}\mathcal{R}_{\mathcal{\overline{G}}_{22}},
 \end{array}\\
\label{XZy996}
&\begin{array}{l}
 \mathcal{\widehat{V}}_{2}=\mathcal{\overline{M}}_{22}^{\dagger}*_{N}\mathcal{\overline{E}}_{22}*_{M}\mathcal{\overline{J}}_{22}^{\dagger}
 +\mathcal{\overline{S}}_{22}^{\dagger}*_{N}\mathcal{\overline{S}}_{22}*_{N}\mathcal{\overline{H}}_{22}^{\dagger}
*_{N}\mathcal{\overline{E}}_{22}*_{M}\mathcal{\overline{N}}_{22}^{\dagger}+\mathcal{L}_{\mathcal{\overline{M}}_{22}}
\\
\ \ \ \ \ \ *_{N}\mathcal{L}_{\mathcal{\overline{S}}_{22}}*_{N}\mathcal{Q}_{77}+\mathcal{L}_{\mathcal{\overline{M}}_{22}}*_{N}\mathcal{Q}_{55}
*_{M}\mathcal{R}_{\mathcal{\overline{N}}_{22}}+\mathcal{Q}_{88}
*_{M}\mathcal{R}_{\mathcal{\overline{J}}_{22}},
 \end{array}
\end{align}
\end{subequations}
By equating $\mathcal{\widehat{V}}_{2}$ in \eqref{XZy993} with $\mathcal{\widehat{V}}_{2}$ in \eqref{XZy996}, we have the following equation:
\begin{align}
\label{yyt11}
&\mathcal{\overline{\overline{A}}}_{1}*_{N}\begin{bmatrix}
                                \mathcal{P}_{55} \\
                                 \mathcal{Q}_{77}\end{bmatrix}+\begin{bmatrix}
                                 \mathcal{P}_{66}&\mathcal{Q}_{88}\end{bmatrix}*_{M} \mathcal{\overline{\overline{B}}}_{1}=-\mathcal{\overline{\overline{E}}}_{1}
                                 +\mathcal{\overline{\overline{F}}}_{1}*_{N}\mathcal{P}_{44}*_{M}\mathcal{\overline{\overline{G}}}_{1}
                                 +\mathcal{\overline{\overline{H}}}_{1}*_{N}\mathcal{Q}_{55}*_{M}\mathcal{\overline{\overline{J}}}_{1},\
\end{align}
It follows from $Lemma$ $\ref{lma 2.3}$  that Eq.\eqref{yyt11} is solvable if and only if there exist quaternion tensors $\mathcal{P}_{44}$ and $\mathcal{Q}_{55}$ satisfy
\begin{align}
\label{yytt11}
&\mathcal{\overline{\overline{F}}}_{11}*_{N}\mathcal{P}_{44}*_{M}\mathcal{\overline{\overline{G}}}_{11}
+\mathcal{\overline{\overline{H}}}_{11}*_{N}\mathcal{Q}_{55}*_{M}\mathcal{\overline{\overline{J}}}_{11}=\mathcal{\overline{\overline{E}}}_{11},\
\end{align}
In utilizing $Lemma$ \ref{lma 2.3}, we have that the general solution to quaternion tensor equation \eqref{yyt11} can be given by \eqref{yWW992}-\eqref{y9WWW93}. On applying $Lemma$ $\ref{lma 2.3}$, we have that Eq.\eqref{yytt11} is solvable if and only if conditions \eqref{b8811}-\eqref{b9911} satisfy. In that case, the general solution can be given by
\begin{align}
\label{XZyEEEE994}
&\begin{array}{l}
 \mathcal{P}_{44}=\mathcal{\overline{\overline{F}}}_{11}^{\dagger}*_{N}\mathcal{\overline{\overline{E}}}_{11}*_{M}\mathcal{\overline{\overline{G}}}_{11}^{\dagger}
-\mathcal{\overline{\overline{F}}}_{11}^{\dagger}*_{N}\mathcal{\overline{\overline{H}}}_{11}*_{N}\mathcal{\overline{\overline{M}}}_{11}^{\dagger}
*_{N}\mathcal{\overline{\overline{E}}}_{11}
*_{M}\mathcal{\overline{\overline{G}}}_{11}^{\dagger}
-\mathcal{\overline{\overline{F}}}_{11}^{\dagger}*_{N}\mathcal{\overline{\overline{S}}}_{11}\\
\ \ \ \ \ *_{N}\mathcal{\overline{\overline{H}}}_{11}^{\dagger}*_{N}\mathcal{\overline{\overline{E}}}_{11}*_{M}\mathcal{\overline{\overline{N}}}_{11}^{\dagger}*_{M}\mathcal{\overline{\overline{J}}}_{11}
*_{M}\mathcal{\overline{\overline{G}}}_{11}^{\dagger}
-\mathcal{\overline{\overline{F}}}_{11}^{\dagger}*_{N}\mathcal{\overline{\overline{S}}}_{11}
*_{N}\mathcal{K}_{44}*_{M}\mathcal{R}_{\mathcal{\overline{\overline{N}}}_{11}}
\\
\ \ \ \ \ *_{M}\mathcal{\overline{\overline{J}}}_{11}
*_{M}\mathcal{\overline{\overline{G}}}_{11}^{\dagger}+\mathcal{L}_{\mathcal{\overline{\overline{F}}}_{11}}*_{N}\mathcal{K}_{55}+\mathcal{K}_{66}*_{M}\mathcal{R}_{\mathcal{\overline{\overline{G}}}_{11}},
 \end{array}\\
\label{XZyEEE994}
&\begin{array}{l}
 \mathcal{Q}_{55}=\mathcal{\overline{\overline{M}}}_{11}^{\dagger}*_{N}\mathcal{\overline{\overline{E}}}_{11}*_{M}\mathcal{\overline{\overline{J}}}_{11}^{\dagger}
 +\mathcal{\overline{\overline{S}}}_{11}^{\dagger}*_{N}\mathcal{\overline{\overline{S}}}_{11}*_{N}\mathcal{\overline{\overline{H}}}_{11}^{\dagger}
*_{N}\mathcal{\overline{\overline{E}}}_{11}*_{M}\mathcal{\overline{\overline{N}}}_{11}^{\dagger}+\mathcal{L}_{\mathcal{\overline{\overline{M}}}_{11}}
\\
\ \ \ \ \ \ *_{N}\mathcal{L}_{\mathcal{\overline{\overline{S}}}_{11}}*_{N}\mathcal{K}_{77}+\mathcal{L}_{\mathcal{\overline{\overline{M}}}_{11}}
*_{N}\mathcal{K}_{44}*_{M}\mathcal{R}_{\mathcal{\overline{\overline{N}}}_{11}}+\mathcal{K}_{88}
*_{M}\mathcal{R}_{\mathcal{\overline{\overline{J}}}_{11}}.
 \end{array}
\end{align}\par
Quaternion tensors $\mathcal{\widehat{U}}_{2}$ in \eqref{x9944} and  $\mathcal{\widehat{K}}_{2}$ in \eqref{x9955} should coincide  with  $\mathcal{\widehat{U}}_{2}$ in \eqref{XZy994} and  $\mathcal{\widehat{K}}_{2}$ in \eqref{XZy995}, respectively. In that case, we have the following equations:
\begin{align}
\label{x88QQ}
&\mathcal{\widetilde{A}}_{1}*_{N}\begin{bmatrix}
                                \mathcal{Q}_{44} \\
                                 \mathcal{V}_{77}\end{bmatrix}+\begin{bmatrix}
                                 \mathcal{Q}_{66}&\mathcal{V}_{88}\end{bmatrix}*_{M} \mathcal{\widetilde{B}}_{1}=\mathcal{\widetilde{E}}_{1}
                                 -\mathcal{\widetilde{C}}_{1}*_{N}\mathcal{P}_{44}*_{M}\mathcal{\widetilde{D}}_{1},\\
\label{x99QQ}
&\mathcal{\widetilde{A}}_{2}*_{N}\begin{bmatrix}
                                \mathcal{V}_{99} \\
                                 \mathcal{P}_{77}\end{bmatrix}+\begin{bmatrix}
                                 \mathcal{W}_{11}&\mathcal{P}_{88}\end{bmatrix}*_{M} \mathcal{\widetilde{B}}_{2}=\mathcal{\widetilde{E}}_{2}
                                 -\mathcal{\widetilde{C}}_{2}*_{N}\mathcal{Q}_{55}*_{M}\mathcal{\widetilde{D}}_{2},
\end{align}
Apply $Lemma$ $\ref{lma 2.3}$ to Eq. \eqref{x88QQ} and Eq.\eqref{x99QQ}. we have that \eqref{x88QQ} and Eq.\eqref{x99QQ} are solvable if and only if there are quaternion tensors $\mathcal{P}_{44}$ and $\mathcal{Q}_{55}$ satisfy
\begin{align}
\label{x9911QQQ}
&\mathcal{\widetilde{C}}_{11}*_{N}\mathcal{P}_{44}*_{M}\mathcal{\widetilde{D}}_{11}
=\mathcal{\widetilde{E}}_{11},\\
\label{x9911QQQQ}
&\mathcal{\widetilde{C}}_{22}*_{N}\mathcal{Q}_{55}*_{M}\mathcal{\widetilde{D}}_{22}
=\mathcal{\widetilde{E}}_{22},
\end{align}
In that case, the general solution to \eqref{x88QQ}-\eqref{x99QQ} can be given by \eqref{x992QQQ}-\eqref{xqq9W93}. Meanwhile, the quaternion tensor equations \eqref{x9911QQQ} and \eqref{x9911QQQQ} are solvable if and only if the conditions \eqref{bb111d} are satisfying, respectively, for $j=1,2$.
In that case, the general solution can be given by
\begin{align}
 \label{x99WWW441}
&\mathcal{P}_{44}=\mathcal{\widetilde{C}}_{11}^{\dagger}*_{N}\mathcal{\widetilde{E}}_{11}
*_{M}\mathcal{\widetilde{D}}_{11}^{\dagger}+\mathcal{L}_{\mathcal{\widetilde{C}}_{11}}*_{N}\mathcal{W}_{88}
+\mathcal{W}_{99}*_{M}\mathcal{R}_{\mathcal{\widetilde{D}}_{11}},\\
\label{x99WWW551}
&\mathcal{Q}_{55}=\mathcal{\widetilde{C}}_{22}^{\dagger}*_{N}\mathcal{\widetilde{E}}_{22}
*_{M}\mathcal{\widetilde{D}}_{22}^{\dagger}+\mathcal{L}_{\mathcal{\widetilde{C}}_{22}}*_{N}\mathcal{T}_{11}
+\mathcal{T}_{22}*_{M}\mathcal{R}_{\mathcal{\widetilde{D}}_{22}}.
\end{align}
Quaternion tensors $\mathcal{P}_{44}$ in \eqref{XZyEEEE994} and $\mathcal{Q}_{55}$ in \eqref{XZyEEE994} should be equal to quaternion tensors $\mathcal{P}_{44}$ in \eqref{x99WWW441} and $\mathcal{Q}_{55}$ in \eqref{x99WWW551}, respectively. Then we have the following system of tensor equations:
\begin{align}
\label{x88QwwQ}
&\mathcal{\widetilde{F}}_{1}*_{N}\begin{bmatrix}
                                \mathcal{W}_{88} \\
                                 \mathcal{K}_{55}\end{bmatrix}+\begin{bmatrix}
                                 \mathcal{W}_{99}&\mathcal{K}_{66}\end{bmatrix}*_{M} \mathcal{\widetilde{G}}_{1}=\mathcal{\widetilde{\widetilde{E}}}_{1}
                                 -\mathcal{\widetilde{H}}_{1}*_{N}\mathcal{K}_{44}*_{M}\mathcal{\widetilde{J}}_{1},\\
\label{x99QwwQ}
&\mathcal{\widetilde{F}}_{2}*_{N}\begin{bmatrix}
                                \mathcal{K}_{77} \\
                                 \mathcal{T}_{11}\end{bmatrix}+\begin{bmatrix}
                                 \mathcal{K}_{88}&\mathcal{T}_{22}\end{bmatrix}*_{M} \mathcal{\widetilde{G}}_{2}=\mathcal{\widetilde{\widetilde{E}}}_{2}
                                 -\mathcal{\widetilde{H}}_{2}*_{N}\mathcal{K}_{44}*_{M}\mathcal{\widetilde{J}}_{2},
\end{align}
Apply $Lemma$ $\ref{lma 2.3},$ to Eq.\eqref{x88QwwQ} and Eq.\eqref{x99QwwQ}. Consequently, we have that
\begin{align}
\label{x9911QeQQ}
&\mathcal{\widetilde{H}}_{11}*_{N}\mathcal{K}_{44}*_{M}\mathcal{\widetilde{J}}_{11}
=\mathcal{\widetilde{\widetilde{E}}}_{11},\\
\label{x9911QQQeQe}
&\mathcal{\widetilde{H}}_{22}*_{N}\mathcal{K}_{44}*_{M}\mathcal{\widetilde{J}}_{22}
=\mathcal{\widetilde{\widetilde{E}}}_{22},
\end{align}
In that case, the general solution to Equations \eqref{x88QwwQ} and \eqref{x99QwwQ} can be given by \eqref{XZEEyTTTr993}-\eqref{kklk}. Meanwhile, the quaternion system of tensor equations \eqref{x9911QeQQ}-\eqref{x9911QQQeQe} is solvable if and only if Eq. \eqref{x9911QeQQ}, Eq. \eqref{x9911QQQeQe} are solvable, and $\mathcal{K}_{44}$ in \eqref{x9911QeQQ} coincide with $\mathcal{K}_{44}$ in \eqref{x9911QQQeQe}. Eq. \eqref{x9911QeQQ} and  Eq. \eqref{x9911QQQeQe} are solvable if and only if the conditions \eqref{bb1121d} are satisfying, respectively, for $l=1,2$. In that case, the general solution to \eqref{x9911QeQQ} and \eqref{x9911QQQeQe} can be given by
\begin{align}
 \label{x99WWW44}
&\mathcal{K}_{44}=\mathcal{\widetilde{H}}_{11}^{\dagger}*_{N}\mathcal{\widetilde{\widetilde{E}}}_{11}
*_{M}\mathcal{\widetilde{J}}_{11}^{\dagger}+\mathcal{L}_{\mathcal{\widetilde{H}}_{11}}*_{N}\mathcal{\acute{W}}_{2}
+\mathcal{\acute{W}}_{3}*_{M}\mathcal{R}_{\mathcal{\widetilde{J}}_{11}},\\
\label{x99WWW55}
&\mathcal{K}_{44}=\mathcal{\widetilde{H}}_{22}^{\dagger}*_{N}\mathcal{\widetilde{\widetilde{E}}}_{22}
*_{M}\mathcal{\widetilde{J}}_{22}^{\dagger}+\mathcal{L}_{\mathcal{\widetilde{H}}_{22}}*_{N}\mathcal{\acute{W}}_{4}
+\mathcal{\acute{W}}_{5}*_{M}\mathcal{R}_{\mathcal{\widetilde{J}}_{22}}.
\end{align}
Ultimately, equating $\mathcal{K}_{44}$ in \eqref{x99WWW44} by $\mathcal{K}_{44}$ in \eqref{x99WWW55}, yield:
\begin{align}
\label{x88QwwwwQ}
&\mathcal{\widetilde{A}}*_{N}\begin{bmatrix}
                                \mathcal{T}_{33} \\
                                 \mathcal{T}_{55}\end{bmatrix}+\begin{bmatrix}
                                 \mathcal{T}_{44}&\mathcal{T}_{66}\end{bmatrix}*_{M} \mathcal{\widetilde{B}}=\mathcal{\widetilde{E}}.
\end{align}
Apply $Lemma$ $\ref{lma 2.3}$ to Eq.\eqref{x88QwwwwQ}, we have that Eq.\eqref{x88QwwwwQ} is solvable if and only if condition \eqref{bb1121d} satisfies. In that case the general solution can be given by \eqref{sss22}-\eqref{sss3HY22}.
\end{proof}
\begin{algorithm}The general solution to the system of two-sided four coupled Sylvester-like quaternion tensor equations \eqref{1.4aa} gives by the following:
\begin{enumerate}
\item \textbf{Input} the system of two-sided four coupled Sylvester-like quaternion tensor equations \eqref{1.4aa} with viable orders over $\mathbb{H}$.
\item Compute all quaternion tensors, which appeared in \eqref{abc11}-\eqref{abcHHD6611}.
 \item Check whether the Moore-Penrose inverses conditions in $Theorem$ $\ref{system 3.33AA}$ are satisfied or not. If not, return ``The system \eqref{1.4aa} is inconsistent".
 \item Else compute the quaternion unknowns $\mathcal{X}_{i}, \mathcal{Y}_{i}, \mathcal{Z}_{j}$, where $(i=\overline{1,3})$ and $(j=\overline{1,4})$ by \eqref{sssa222p1}-\eqref{sssa333}.
 \item \textbf{Output} the general solution of the system \eqref{1.4aa} is $\mathcal{X}_{i}, \mathcal{Y}_{i}, \mathcal{Z}_{j}$.
\end{enumerate}
\end{algorithm}
We give an example to illustrate $Theorem$ $\ref{system 3.33AA}$.
\begin{example}Consider the two-sided four coupled Sylvester-like quaternion system of tensor equations \eqref{1.4aa}, where
\begin{small}
\begin{align*}
\begin{gathered}
\mathcal{F}_{4}(:,:,1,1)=\begin{bmatrix}
             \mathbf{0} & \mathbf{i} \\
            \mathbf{0} & \mathbf{0}
           \end{bmatrix},\
\mathcal{F}_{4}(:,:,1,2)=\begin{bmatrix}
             \mathbf{1} & \mathbf{0} \\
            \mathbf{0} & \mathbf{j}
           \end{bmatrix},\
\mathcal{F}_{4}(:,:,2,1)=\begin{bmatrix}
             \mathbf{0} & \mathbf{0} \\
            \mathbf{0} & \mathbf{-k}
           \end{bmatrix},\
 \mathcal{F}_{4}(:,:,2,2)=\begin{bmatrix}
             \mathbf{1} & \mathbf{i} \\
            \mathbf{0} & \mathbf{0}
           \end{bmatrix},\\
\mathcal{G}_{4}(:,:,1,1)=\begin{bmatrix}
             \mathbf{0} & \mathbf{0} \\
            \mathbf{2} & \mathbf{i}
           \end{bmatrix},\
\mathcal{G}_{4}(:,:,1,2)=\begin{bmatrix}
             \mathbf{2} & \mathbf{0} \\
            \mathbf{0} & \mathbf{2i}
           \end{bmatrix},\
\mathcal{G}_{4}(:,:,2,1)=\begin{bmatrix}
             \mathbf{0} & \mathbf{0} \\
            \mathbf{0} & \mathbf{k}
           \end{bmatrix},\
 \mathcal{G}_{4}(:,:,2,2)=\begin{bmatrix}
             \mathbf{0} & \mathbf{3-j} \\
            \mathbf{0} & \mathbf{0}
           \end{bmatrix},\\
\mathcal{E}_{4}(:,:,1,1)=\begin{bmatrix}
             \mathbf{-2j} & \mathbf{2i} \\
            \mathbf{0} & \mathbf{-2j}
           \end{bmatrix},\
\mathcal{E}_{4}(:,:,1,2)=\begin{bmatrix}
             \mathbf{0} & \mathbf{-4i} \\
            \mathbf{0} & \mathbf{2i-6j}
           \end{bmatrix},\
\mathcal{E}_{4}(:,:,2,1)=\begin{bmatrix}
             \mathbf{0} & \mathbf{0} \\
            \mathbf{0} & \mathbf{2+k}
           \end{bmatrix},\\
\mathcal{E}_{4}(:,:,2,2)=\begin{bmatrix}
             \mathbf{6-2j} & \mathbf{-1-3j} \\
            \mathbf{0} & \mathbf{2-3i+6j+k}
           \end{bmatrix},\
\mathcal{H}_{4}(:,:,1,1)=\begin{bmatrix}
             \mathbf{i} & \mathbf{j} \\
            \mathbf{0} & \mathbf{0}
           \end{bmatrix},\
\mathcal{H}_{4}(:,:,1,2)=\begin{bmatrix}
             \mathbf{j} & \mathbf{k} \\
            \mathbf{0} & \mathbf{0}
           \end{bmatrix},\\
\mathcal{H}_{4}(:,:,2,1)=\begin{bmatrix}
             \mathbf{i} & \mathbf{k} \\
            \mathbf{0} & \mathbf{0}
           \end{bmatrix},\
 \mathcal{H}_{4}(:,:,2,2)=\begin{bmatrix}
             \mathbf{1} & \mathbf{0} \\
            \mathbf{2-i} & \mathbf{0}
           \end{bmatrix},\
\mathcal{J}_{4}(:,:,1,1)=\begin{bmatrix}
             \mathbf{0} & \mathbf{0} \\
            \mathbf{k} & \mathbf{i}
           \end{bmatrix},\
\mathcal{J}_{4}(:,:,1,2)=\begin{bmatrix}
             \mathbf{0} & \mathbf{2} \\
            \mathbf{0} & \mathbf{-j}
           \end{bmatrix},\\
\mathcal{J}_{4}(:,:,2,1)=\begin{bmatrix}
             \mathbf{0} & \mathbf{0} \\
            \mathbf{2-i} & \mathbf{2+j}
           \end{bmatrix},\
 \mathcal{J}_{4}(:,:,2,2)=\begin{bmatrix}
             \mathbf{0} & \mathbf{3-k} \\
            \mathbf{0} & \mathbf{3+k}
           \end{bmatrix},\
\mathcal{E}_{5}(:,:,1,1)=\begin{bmatrix}
             \mathbf{i+2j-k} & \mathbf{1+i} \\
            \mathbf{1+2i+2j-k} & \mathbf{0}
           \end{bmatrix},\\
\mathcal{E}_{5}(:,:,1,2)=\begin{bmatrix}
             \mathbf{10i+2j} & \mathbf{-2i+6j+k} \\
            \mathbf{3+6i} & \mathbf{0}
           \end{bmatrix},\
\mathcal{E}_{5}(:,:,2,1)=\begin{bmatrix}
             \mathbf{-2+i+4j+5k} & \mathbf{1-2i+2j+k} \\
            \mathbf{9j+3k} & \mathbf{0}
           \end{bmatrix},\\
\mathcal{E}_{5}(:,:,2,2)=\begin{bmatrix}
             \mathbf{-2+11i+7j+6k} & \mathbf{-9i+9j} \\
            \mathbf{1+7i+5j+5k} & \mathbf{0}
           \end{bmatrix},\
\mathcal{A}_{1}(:,:,1,1)=\begin{bmatrix}
             \mathbf{i} & \mathbf{-i} \\
            \mathbf{0} & \mathbf{0}
           \end{bmatrix},\
\mathcal{A}_{1}(:,:,1,2)=\begin{bmatrix}
             \mathbf{j} & \mathbf{-j} \\
            \mathbf{0} & \mathbf{0}
           \end{bmatrix},\\
\mathcal{A}_{1}(:,:,2,1)=\begin{bmatrix}
             \mathbf{k} & \mathbf{-k} \\
            \mathbf{0} & \mathbf{0}
           \end{bmatrix},\
 \mathcal{A}_{1}(:,:,2,2)=\begin{bmatrix}
             \mathbf{0} & \mathbf{0} \\
            \mathbf{i} & \mathbf{-i}
           \end{bmatrix},\
\mathcal{B}_{1}(:,:,1,1)=\begin{bmatrix}
             \mathbf{0} & \mathbf{0} \\
            \mathbf{j} & \mathbf{-j}
           \end{bmatrix},\
\mathcal{B}_{1}(:,:,1,2)=\begin{bmatrix}
             \mathbf{0} & \mathbf{0} \\
            \mathbf{k} & \mathbf{-k}
           \end{bmatrix},\\
\mathcal{B}_{1}(:,:,2,1)=\begin{bmatrix}
             \mathbf{i} & \mathbf{i+j} \\
            \mathbf{0} & \mathbf{0}
           \end{bmatrix},\
 \mathcal{B}_{1}(:,:,2,2)=\begin{bmatrix}
             \mathbf{j} & \mathbf{j+k} \\
            \mathbf{0} & \mathbf{0}
           \end{bmatrix},\
\mathcal{C}_{1}(:,:,1,1)=\begin{bmatrix}
             \mathbf{k} & \mathbf{i+k} \\
            \mathbf{0} & \mathbf{0}
           \end{bmatrix},
\mathcal{C}_{1}(:,:,2,1)=\begin{bmatrix}
             \mathbf{5i} & \mathbf{1} \\
            \mathbf{0} & \mathbf{0}
           \end{bmatrix},\\
\mathcal{C}_{1}(:,:,1,2)=\begin{bmatrix}
             \mathbf{2-i} & \mathbf{0} \\
            \mathbf{0} & \mathbf{2k}
           \end{bmatrix},\
 \mathcal{C}_{1}(:,:,2,2)=\begin{bmatrix}
             \mathbf{0} & \mathbf{i} \\
            \mathbf{k} & \mathbf{0}
           \end{bmatrix},\
\mathcal{D}_{1}(:,:,1,1)=\begin{bmatrix}
             \mathbf{k} & \mathbf{2-k} \\
            \mathbf{j} & \mathbf{0}
           \end{bmatrix},
\mathcal{D}_{1}(:,:,1,2)=\begin{bmatrix}
             \mathbf{0} & \mathbf{0} \\
            \mathbf{0} & \mathbf{3i}
           \end{bmatrix},           
\end{gathered}  
\end{align*}
\begin{align*} 
\begin{gathered}          
\mathcal{D}_{1}(:,:,2,1)=\begin{bmatrix}
             \mathbf{i} & \mathbf{j} \\
            \mathbf{k} & \mathbf{2}
           \end{bmatrix},\
 \mathcal{D}_{1}(:,:,2,2)=\begin{bmatrix}
             \mathbf{2} & \mathbf{i} \\
            \mathbf{0} & \mathbf{j}
           \end{bmatrix},\
\mathcal{F}_{1}(:,:,1,1)=\begin{bmatrix}
             \mathbf{0} & \mathbf{k} \\
            \mathbf{0} & \mathbf{2i}
           \end{bmatrix},
\mathcal{F}_{1}(:,:,1,2)=\begin{bmatrix}
             \mathbf{-1} & \mathbf{-1+i} \\
            \mathbf{0} & \mathbf{0}
           \end{bmatrix},\\
\mathcal{F}_{1}(:,:,2,1)=\begin{bmatrix}
             \mathbf{0} & \mathbf{1} \\
            \mathbf{0} & \mathbf{1-i}
           \end{bmatrix},\
 \mathcal{F}_{1}(:,:,2,2)=\begin{bmatrix}
             \mathbf{0} & \mathbf{0} \\
            \mathbf{2} & \mathbf{2-i}
           \end{bmatrix},\
\mathcal{G}_{1}(:,:,1,2)=\begin{bmatrix}
             \mathbf{i} & \mathbf{2-i} \\
            \mathbf{j} & \mathbf{0}
           \end{bmatrix},\
\mathcal{H}_{1}(:,:,1,2)=\begin{bmatrix}
             \mathbf{0} & \mathbf{i} \\
            \mathbf{j} & \mathbf{0}
           \end{bmatrix},\\
\mathcal{G}_{1}(:,:,1,1)=\begin{bmatrix}
             \mathbf{0} & \mathbf{0} \\
            \mathbf{0} & \mathbf{i+j+k}
           \end{bmatrix},\
\mathcal{G}_{1}(:,:,2,1)=\begin{bmatrix}
             \mathbf{j} & \mathbf{2-j} \\
            \mathbf{k} & \mathbf{0}
           \end{bmatrix},\
\mathcal{H}_{1}(:,:,2,2)=\begin{bmatrix}
             \mathbf{0} & \mathbf{k} \\
            \mathbf{i} & \mathbf{0}
           \end{bmatrix},\
\mathcal{J}_{1}(:,:,1,1)=\begin{bmatrix}
             \mathbf{0} & \mathbf{j} \\
            \mathbf{k} & \mathbf{0}
           \end{bmatrix},\\
\mathcal{G}_{1}(:,:,2,2)=\begin{bmatrix}
             \mathbf{0} & \mathbf{0} \\
            \mathbf{0} & \mathbf{i+k}
           \end{bmatrix},\
\mathcal{H}_{1}(:,:,1,1)=\begin{bmatrix}
             \mathbf{i-k} & \mathbf{0} \\
            \mathbf{0} & \mathbf{-2k}
           \end{bmatrix},\
\mathcal{J}_{1}(:,:,1,2)=\begin{bmatrix}
             \mathbf{j} & \mathbf{0} \\
            \mathbf{0} & \mathbf{k}
           \end{bmatrix},\
\mathcal{J}_{1}(:,:,2,1)=\begin{bmatrix}
             \mathbf{i} & \mathbf{0} \\
            \mathbf{0} & \mathbf{j}
           \end{bmatrix},\\
\mathcal{H}_{1}(:,:,2,1)=\begin{bmatrix}
             \mathbf{i} & \mathbf{j-k} \\
            \mathbf{0} & \mathbf{0}
           \end{bmatrix},\
\mathcal{E}_{1}(:,:,1,1)=\begin{bmatrix}
             \mathbf{49+19i+5j+23k} & \mathbf{3-10i-3j-14k} \\
            \mathbf{1+3i+2j} & \mathbf{-12+6i+3j-7k}
           \end{bmatrix},\\
\mathcal{J}_{1}(:,:,2,2)=\begin{bmatrix}
             \mathbf{1-i} & \mathbf{0} \\
            \mathbf{0} & \mathbf{j}
           \end{bmatrix},\
\mathcal{E}_{1}(:,:,1,2)=\begin{bmatrix}
             \mathbf{14-44i-19j-2k} & \mathbf{-6-i+16j} \\
            \mathbf{5-5i-3j+3k} & \mathbf{12+10i+7j+3k}
           \end{bmatrix},\\
\mathcal{E}_{1}(:,:,2,1)=\begin{bmatrix}
             \mathbf{-2-4i-39j+82k} & \mathbf{-3-i+7j+27k} \\
            \mathbf{3-15i5j-k} & \mathbf{-14-2i-14j-4k}
           \end{bmatrix},\
\mathcal{A}_{2}(:,:,1,1)=\begin{bmatrix}
             \mathbf{2} & \mathbf{i} \\
            \mathbf{0} & \mathbf{0}
           \end{bmatrix},\\
\mathcal{E}_{1}(:,:,2,2)=\begin{bmatrix}
             \mathbf{58-18i-45j+34k} & \mathbf{14-21i+11k} \\
            \mathbf{-4-10i-8j-12k} & \mathbf{-2+22i-14j-4k}
           \end{bmatrix},\
\mathcal{A}_{2}(:,:,1,2)=\begin{bmatrix}
             \mathbf{0} & \mathbf{-5k} \\
            \mathbf{0} & \mathbf{k}
           \end{bmatrix},\\
\mathcal{A}_{2}(:,:,2,1)=\begin{bmatrix}
             \mathbf{3} & \mathbf{0} \\
            \mathbf{-i} & \mathbf{0}
           \end{bmatrix},\
\mathcal{A}_{2}(:,:,2,2)=\begin{bmatrix}
             \mathbf{0} & \mathbf{-6} \\
            \mathbf{0} & \mathbf{-i}
           \end{bmatrix},\
\mathcal{B}_{2}(:,:,1,1)=\begin{bmatrix}
             \mathbf{4} & \mathbf{j} \\
            \mathbf{0} & \mathbf{0}
           \end{bmatrix},\
 \mathcal{B}_{2}(:,:,1,2)=\begin{bmatrix}
             \mathbf{-i} & \mathbf{-7} \\
            \mathbf{0} & \mathbf{0}
           \end{bmatrix},\\
\mathcal{B}_{2}(:,:,2,1)=\begin{bmatrix}
             \mathbf{5} & \mathbf{0} \\
            \mathbf{-j} & \mathbf{0}
           \end{bmatrix},\
\mathcal{B}_{2}(:,:,2,2)=\begin{bmatrix}
             \mathbf{j} & \mathbf{-8} \\
            \mathbf{0} & \mathbf{0}
           \end{bmatrix},\
\mathcal{C}_{2}(:,:,1,1)=\begin{bmatrix}
             \mathbf{6} & \mathbf{k} \\
            \mathbf{0} & \mathbf{-k}
           \end{bmatrix},\
 \mathcal{C}_{2}(:,:,1,2)=\begin{bmatrix}
             \mathbf{k} & \mathbf{0} \\
            \mathbf{9} & \mathbf{0}
           \end{bmatrix},\\
\mathcal{C}_{2}(:,:,2,1)=\begin{bmatrix}
             \mathbf{7} & \mathbf{0} \\
            \mathbf{0} & \mathbf{-3j}
           \end{bmatrix},\
\mathcal{C}_{2}(:,:,2,2)=\begin{bmatrix}
             \mathbf{2k} & \mathbf{0} \\
            \mathbf{8} & \mathbf{0}
           \end{bmatrix},\
\mathcal{D}_{2}(:,:,1,1)=\begin{bmatrix}
             \mathbf{8} & \mathbf{0} \\
            \mathbf{0} & \mathbf{3j}
           \end{bmatrix},\
 \mathcal{D}_{2}(:,:,1,2)=\begin{bmatrix}
             \mathbf{0} & \mathbf{0} \\
            \mathbf{7} & \mathbf{0}
           \end{bmatrix},\\
\mathcal{D}_{2}(:,:,2,1)=\begin{bmatrix}
             \mathbf{9} & \mathbf{i-k} \\
            \mathbf{0} & \mathbf{0}
           \end{bmatrix},\
\mathcal{D}_{2}(:,:,2,2)=\begin{bmatrix}
             \mathbf{i} & \mathbf{0} \\
            \mathbf{6} & \mathbf{-k}
           \end{bmatrix},\
\mathcal{F}_{2}(:,:,1,1)=\begin{bmatrix}
             \mathbf{0} & \mathbf{0} \\
            \mathbf{i} & \mathbf{i-k}
           \end{bmatrix},\\
\mathcal{F}_{2}(:,:,2,1)=\begin{bmatrix}
             \mathbf{0} & \mathbf{2i} \\
            \mathbf{k} & \mathbf{0}
           \end{bmatrix},\
\mathcal{F}_{2}(:,:,2,2)=\begin{bmatrix}
             \mathbf{0} & \mathbf{k-j} \\
            \mathbf{0} & \mathbf{j}
           \end{bmatrix},\
\mathcal{G}_{2}(:,:,1,1)=\begin{bmatrix}
             \mathbf{0} & \mathbf{-k} \\
            \mathbf{2k} & \mathbf{0}
           \end{bmatrix},\\
\mathcal{G}_{2}(:,:,2,1)=\begin{bmatrix}
             \mathbf{1} & \mathbf{0} \\
            \mathbf{2i} & \mathbf{0}
           \end{bmatrix},\
\mathcal{G}_{2}(:,:,2,2)=\begin{bmatrix}
             \mathbf{0} & \mathbf{5} \\
            \mathbf{0} & \mathbf{3j}
           \end{bmatrix},\
\mathcal{H}_{2}(:,:,1,1)=\begin{bmatrix}
             \mathbf{2i} & \mathbf{0} \\
            \mathbf{2j} & \mathbf{0}
           \end{bmatrix},\
 \mathcal{H}_{2}(:,:,1,2)=\begin{bmatrix}
             \mathbf{0} & \mathbf{i-k} \\
            \mathbf{0} & \mathbf{i+k}
           \end{bmatrix},\\
\mathcal{H}_{2}(:,:,2,1)=\begin{bmatrix}
             \mathbf{-i+j} & \mathbf{0} \\
            \mathbf{0} & \mathbf{i+k}
           \end{bmatrix},\
\mathcal{H}_{2}(:,:,2,2)=\begin{bmatrix}
             \mathbf{0} & \mathbf{3} \\
            \mathbf{3-i} & \mathbf{j}
           \end{bmatrix},\
\mathcal{J}_{2}(:,:,1,1)=\begin{bmatrix}
             \mathbf{i} & \mathbf{0} \\
            \mathbf{0} & \mathbf{2j}
           \end{bmatrix},\\
\mathcal{J}_{2}(:,:,2,1)=\begin{bmatrix}
             \mathbf{0} & \mathbf{0} \\
            \mathbf{0} & \mathbf{2i}
           \end{bmatrix},\
\mathcal{E}_{2}(:,:,1,1)=\begin{bmatrix}
             \mathbf{-89+649i-668j+5k} & \mathbf{-2+88i+25j-20k} \\
            \mathbf{-271-523i-241j+6k} & \mathbf{-33-94i+46j+252k}
           \end{bmatrix},\\
\mathcal{J}_{2}(:,:,2,2)=\begin{bmatrix}
             \mathbf{3i} & \mathbf{0} \\
            \mathbf{-j} & \mathbf{0}
           \end{bmatrix},\
 \mathcal{E}_{2}(:,:,1,2)=\begin{bmatrix}
             \mathbf{146+329i-635j-147k} & \mathbf{42+33j+35k} \\
            \mathbf{108-95i-860j+77k} & \mathbf{-318-74i-179j+15k}
           \end{bmatrix},\\
 \mathcal{E}_{2}(:,:,2,1)=\begin{bmatrix}
             \mathbf{-131+821i-879j+151k} & \mathbf{-20+105i+27j-25k} \\
            \mathbf{19-682i -160j-3k} & \mathbf{-24-113i+62j+309k}
           \end{bmatrix},\
\mathcal{A}_{3}(:,:,1,1)=\begin{bmatrix}
             \mathbf{0} & \mathbf{i-k} \\
            \mathbf{0} & \mathbf{j}
           \end{bmatrix},\           \\
\mathcal{E}_{2}(:,:,2,2)=\begin{bmatrix}
             \mathbf{166+368i-705j+192k} & \mathbf{17+3i+32j+53k} \\
            \mathbf{402-114i-975j-15k} & \mathbf{-337-148i-146j-8k}
           \end{bmatrix},\
\mathcal{A}_{3}(:,:,1,2)=\begin{bmatrix}
             \mathbf{i+j} & \mathbf{0} \\
            \mathbf{i-2j} & \mathbf{0}
           \end{bmatrix},
\end{gathered}             
 \end{align*}
\begin{align*} 
\begin{gathered}
\mathcal{A}_{3}(:,:,2,1)=\begin{bmatrix}
             \mathbf{0} & \mathbf{0} \\
            \mathbf{j+k} & \mathbf{i}
           \end{bmatrix},\
\mathcal{A}_{3}(:,:,2,2)=\begin{bmatrix}
             \mathbf{j} & \mathbf{-k} \\
            \mathbf{0} & \mathbf{0}
           \end{bmatrix},\
\mathcal{B}_{3}(:,:,1,1)=\begin{bmatrix}
             \mathbf{i} & \mathbf{0} \\
            \mathbf{j} & \mathbf{0}
           \end{bmatrix},\
 \mathcal{B}_{3}(:,:,1,2)=\begin{bmatrix}
             \mathbf{0} & \mathbf{j} \\
            \mathbf{k} & \mathbf{0}
           \end{bmatrix},\\
\mathcal{B}_{3}(:,:,2,1)=\begin{bmatrix}
             \mathbf{k} & \mathbf{i} \\
            \mathbf{0} & \mathbf{0}
           \end{bmatrix},\
\mathcal{B}_{3}(:,:,2,2)=\begin{bmatrix}
             \mathbf{i} & \mathbf{0} \\
            \mathbf{0} & \mathbf{i}
           \end{bmatrix},\
\mathcal{C}_{3}(:,:,1,1)=\begin{bmatrix}
             \mathbf{0} & \mathbf{j} \\
            \mathbf{j} & \mathbf{0}
           \end{bmatrix},\
 \mathcal{C}_{3}(:,:,1,2)=\begin{bmatrix}
             \mathbf{0} & \mathbf{k} \\
            \mathbf{0} & \mathbf{k}
           \end{bmatrix},\\         
\mathcal{C}_{3}(:,:,2,1)=\begin{bmatrix}
             \mathbf{i+j} & \mathbf{0} \\
            \mathbf{0} & \mathbf{k}
           \end{bmatrix},\
\mathcal{C}_{3}(:,:,2,2)=\begin{bmatrix}
             \mathbf{0} & \mathbf{0} \\
            \mathbf{j+k} & \mathbf{i}
           \end{bmatrix},\
\mathcal{D}_{3}(:,:,1,1)=\begin{bmatrix}
             \mathbf{0} & \mathbf{i+k} \\
            \mathbf{0} & \mathbf{j}
           \end{bmatrix},\\
\mathcal{D}_{3}(:,:,2,1)=\begin{bmatrix}
             \mathbf{0} & \mathbf{0} \\
            \mathbf{j-k} & \mathbf{i}
           \end{bmatrix},\
\mathcal{D}_{3}(:,:,2,2)=\begin{bmatrix}
             \mathbf{k-i} & \mathbf{j} \\
            \mathbf{0} & \mathbf{0}
           \end{bmatrix},\
\mathcal{F}_{3}(:,:,1,1)=\begin{bmatrix}
             \mathbf{2j} & \mathbf{0} \\
            \mathbf{3k} & \mathbf{0}
           \end{bmatrix},\\
\mathcal{F}_{3}(:,:,2,1)=\begin{bmatrix}
             \mathbf{0} & \mathbf{0} \\
            \mathbf{i+k} & \mathbf{-k}
           \end{bmatrix},\
\mathcal{F}_{3}(:,:,2,2)=\begin{bmatrix}
             \mathbf{k} & \mathbf{0} \\
            \mathbf{2k} & \mathbf{0}
           \end{bmatrix},\
\mathcal{G}_{3}(:,:,1,1)=\begin{bmatrix}
             \mathbf{0} & \mathbf{j} \\
            \mathbf{2j} & \mathbf{0}
           \end{bmatrix},\
 \mathcal{G}_{3}(:,:,1,2)=\begin{bmatrix}
             \mathbf{0} & \mathbf{i} \\
            \mathbf{3i} & \mathbf{0}
           \end{bmatrix},\\
\mathcal{G}_{3}(:,:,2,1)=\begin{bmatrix}
             \mathbf{i-j} & \mathbf{0} \\
            \mathbf{0} & \mathbf{-j}
           \end{bmatrix},\
\mathcal{G}_{3}(:,:,2,2)=\begin{bmatrix}
             \mathbf{j+k} & \mathbf{0} \\
            \mathbf{0} & \mathbf{-k}
           \end{bmatrix},\
\mathcal{H}_{3}(:,:,1,1)=\begin{bmatrix}
             \mathbf{0} & \mathbf{i+j} \\
            \mathbf{0} & \mathbf{k}
           \end{bmatrix},\\
\mathcal{H}_{3}(:,:,2,1)=\begin{bmatrix}
             \mathbf{i+j} & \mathbf{i} \\
            \mathbf{0} & \mathbf{0}
           \end{bmatrix},\
\mathcal{H}_{3}(:,:,2,2)=\begin{bmatrix}
             \mathbf{j+k} & \mathbf{0} \\
            \mathbf{j} & \mathbf{k}
           \end{bmatrix},\
\mathcal{J}_{3}(:,:,1,1)=\begin{bmatrix}
             \mathbf{i+j} & \mathbf{0} \\
            \mathbf{k} & \mathbf{0}
           \end{bmatrix},\\
\mathcal{J}_{3}(:,:,2,1)=\begin{bmatrix}
             \mathbf{0} & \mathbf{0} \\
            \mathbf{j-k} & \mathbf{-j}
           \end{bmatrix},\
\mathcal{E}_{3}(:,:,1,1)=\begin{bmatrix}
             \mathbf{4-7i-2j+10k} & \mathbf{9+2i+2j-9k} \\
            \mathbf{21-4i-5j-13k} & \mathbf{-1-5i+4j-4k}
           \end{bmatrix},\\
\mathcal{J}_{3}(:,:,2,2)=\begin{bmatrix}
             \mathbf{i+j} & \mathbf{i} \\
            \mathbf{0} & \mathbf{0}
           \end{bmatrix},\
\mathcal{E}_{3}(:,:,1,2)=\begin{bmatrix}
             \mathbf{-31+20i-18j+12k} & \mathbf{-17-9i+4j+10k} \\
            \mathbf{-14-13i+11j} & \mathbf{-22-22i+4j+7k}
           \end{bmatrix},\\
\mathcal{E}_{3}(:,:,2,1)=\begin{bmatrix}
             \mathbf{-10+4i+18j-28k} & \mathbf{-5+5i+14j-10k} \\
            \mathbf{1+4i+23j-11k} & \mathbf{3i+20j+10k}
           \end{bmatrix},\
 \mathcal{F}_{2}(:,:,1,2)=\begin{bmatrix}
             \mathbf{0} & \mathbf{j} \\
            \mathbf{0} & \mathbf{j-k}
           \end{bmatrix},\\
\mathcal{E}_{3}(:,:,2,1)=\begin{bmatrix}
             \mathbf{24+15i+8j-6k} & \mathbf{11+11i-7j-k} \\
            \mathbf{23+2j+6k} & \mathbf{27i-7j+6k}
           \end{bmatrix},\
\mathcal{J}_{2}(:,:,1,2)=\begin{bmatrix}
             \mathbf{k} & \mathbf{-2k} \\
            \mathbf{0} & \mathbf{0}
           \end{bmatrix},\\
\mathcal{D}_{3}(:,:,1,2)=\begin{bmatrix}
             \mathbf{i-j} & \mathbf{0} \\
            \mathbf{k} & \mathbf{0}
           \end{bmatrix},\
\mathcal{F}_{3}(:,:,1,2)=\begin{bmatrix}
             \mathbf{i-k} & \mathbf{-k} \\
            \mathbf{0} & \mathbf{0}
           \end{bmatrix},\
\mathcal{H}_{3}(:,:,1,2)=\begin{bmatrix}
             \mathbf{0} & \mathbf{j+k} \\
            \mathbf{j} & \mathbf{k}
           \end{bmatrix},\\
\mathcal{J}_{3}(:,:,1,2)=\begin{bmatrix}
             \mathbf{0} & \mathbf{i+k} \\
            \mathbf{i} & \mathbf{k}
           \end{bmatrix},\
 \mathcal{G}_{2}(:,:,1,2)=\begin{bmatrix}
             \mathbf{-2} & \mathbf{0} \\
            \mathbf{i-k} & \mathbf{0}
           \end{bmatrix}.
 \end{gathered}          
\end{align*}
\end{small}
We now look at the system \eqref{1.4aa}. Rendering of direct calculations
\begin{align*}
\begin{gathered}   
\mathcal{R}_{\mathcal{M}_{i}}*_{2}\mathcal{R}_{\mathcal{A}_{i}}*_{2}\mathcal{E}_{i}=0,\ \mathcal{E}_{i}*_{2}\mathcal{L}_{\mathcal{B}_{i}}*_{2}\mathcal{L}_{\mathcal{N}_{i}}=0,\
\mathcal{R}_{\mathcal{C}_{i}}*_{2}\mathcal{E}_{i}*_{2}\mathcal{L}_{\mathcal{B}_{i}}=0,\\
\mathcal{R}_{\mathcal{\widehat{M}}_{i}}*_{2}\mathcal{R}_{\mathcal{\widehat{A}}_{i}}*_{2}\mathcal{\widehat{E}}_{i}=0,\ \mathcal{\widehat{E}}_{i}*_{2}\mathcal{L}_{\mathcal{\widehat{B}}_{i}}*_{2}\mathcal{L}_{\mathcal{\widehat{N}}_{i}}=0,\\
\mathcal{R}_{\mathcal{\widehat{A}}_{i}}*_{2}\mathcal{\widehat{E}}_{i}*_{2}\mathcal{L}_{\mathcal{\widehat{D}}_{i}}=0,\
\mathcal{R}_{\mathcal{\widehat{C}}_{i}}*_{2}\mathcal{\widehat{E}}_{i}*_{2}\mathcal{L}_{\mathcal{\widehat{B}}_{i}}=0,\ (i=\overline{1,3}),\\
\mathcal{R}_{\mathcal{F}_{4}}*_{2}\mathcal{E}_{4}=0,\ \mathcal{E}_{4}*_{2}\mathcal{L}_{\mathcal{G}_{4}}=0,\
\mathcal{R}_{\mathcal{H}_{4}}*_{2}\mathcal{E}_{5}=0,\ \mathcal{E}_{5}*_{2}\mathcal{L}_{\mathcal{J}_{4}}=0,\\
\mathcal{R}_{\mathcal{\widehat{\widehat{A}}}_{kk}}*_{2}\mathcal{\widehat{\widehat{E}}}_{kk}=0,\
\mathcal{\widehat{\widehat{E}}}_{kk}*_{2}\mathcal{L}_{\mathcal{\widehat{\widehat{B}}}_{kk}}=0,\ \\
\mathcal{R}_{\mathcal{\overline{M}}_{kk}}*_{2}\mathcal{R}_{\mathcal{\overline{F}}_{kk}}*_{2}\mathcal{\overline{E}}_{kk}=0,\ \mathcal{\overline{E}}_{kk}*_{2}\mathcal{L}_{\mathcal{\overline{G}}_{kk}}*_{2}\mathcal{L}_{\mathcal{\overline{N}}_{kk}}=0,\\
\mathcal{R}_{\mathcal{\overline{F}}_{kk}}*_{2}\mathcal{\overline{E}}_{kk}*_{2}\mathcal{L}_{\mathcal{\overline{J}}_{kk}}=0,\
\mathcal{R}_{\mathcal{\overline{H}}_{kk}}*_{2}\mathcal{\overline{E}}_{kk}*_{2}\mathcal{L}_{\mathcal{\overline{G}}_{kk}}=0,\ (k=1,2),\\
\mathcal{R}_{\mathcal{\overline{\overline{M}}}_{11}}*_{2}\mathcal{R}_{\mathcal{\overline{\overline{F}}}_{11}}*_{2}\mathcal{\overline{\overline{E}}}_{11}=0,\ \mathcal{\overline{\overline{E}}}_{11}*_{2}\mathcal{L}_{\mathcal{\overline{\overline{G}}}_{11}}*_{2}\mathcal{L}_{\mathcal{\overline{\overline{N}}}_{11}}=0,
\end{gathered}
\end{align*}
\begin{align*}
\begin{gathered} 
\mathcal{R}_{\mathcal{\overline{\overline{F}}}_{11}}*_{2}\mathcal{\overline{\overline{E}}}_{11}*_{2}\mathcal{L}_{\mathcal{\overline{\overline{J}}}_{11}}=0,\
\mathcal{R}_{\mathcal{\overline{\overline{H}}}_{11}}*_{2}\mathcal{\overline{\overline{E}}}_{11}*_{2}\mathcal{L}_{\mathcal{\overline{\overline{G}}}_{11}}=0,\\
\mathcal{R}_{\mathcal{\widetilde{C}}_{jj}}*_{2}\mathcal{\widetilde{E}}_{jj}=0,\
\mathcal{\widetilde{E}}_{jj}*_{2}\mathcal{L}_{\mathcal{\widetilde{D}}_{jj}}=0\ (j=1,2),\\  
\mathcal{R}_{\mathcal{\widetilde{H}}_{ll}}*_{2}\mathcal{\widetilde{\widetilde{E}}}_{ll}=0,\
\mathcal{\widetilde{\widetilde{E}}}_{ll}*_{2}\mathcal{L}_{\mathcal{\widetilde{G}}_{ll}}=0,\
\mathcal{R}_{\mathcal{\widetilde{A}}}*_{2}\mathcal{\widetilde{E}}*_{2}\mathcal{L}_{\mathcal{\widetilde{B}}},\ (l=1,2),\\
\mathcal{R}_{\mathcal{\widetilde{A}}}*_{2}\mathcal{\widetilde{E}}*_{2}\mathcal{L}_{\mathcal{\widetilde{B}}}=0.
\end{gathered}   
\end{align*}
Consequently, in $Theorem$ $\ref{system 3.33AA},$ all Moore-Penrose inverse conditions hold, and the system \eqref{1.4aa} is thus consistent. Moreover it is simple to show that the following structures satisfy the system:
\begin{small}
\begin{align*}
\begin{gathered}   
\mathcal{Z}_{1}(:,:,1,1)=\begin{bmatrix}
             \mathbf{-2} & \mathbf{0} \\
            \mathbf{i} & \mathbf{0}
           \end{bmatrix},\
\mathcal{Z}_{1}(:,:,1,2)=\begin{bmatrix}
             \mathbf{0} & \mathbf{2-k} \\
            \mathbf{0} & \mathbf{k}
           \end{bmatrix},\
\mathcal{Z}_{1}(:,:,2,1)=\begin{bmatrix}
             \mathbf{1+j} & \mathbf{0} \\
            \mathbf{0} & \mathbf{-j}
           \end{bmatrix},\\
\mathcal{Z}_{4}(:,:,1,1)=\begin{bmatrix}
             \mathbf{2} & \mathbf{i} \\
            \mathbf{0} & \mathbf{0}
           \end{bmatrix},\
\mathcal{Z}_{4}(:,:,1,2)=\begin{bmatrix}
             \mathbf{3-k} & \mathbf{0} \\
            \mathbf{0} & \mathbf{i}
           \end{bmatrix},\
\mathcal{Z}_{4}(:,:,2,1)=\begin{bmatrix}
             \mathbf{0} & \mathbf{0} \\
            \mathbf{i} & \mathbf{j}
           \end{bmatrix},\
 \mathcal{Z}_{4}(:,:,2,2)=\begin{bmatrix}
             \mathbf{0} & \mathbf{0} \\
            \mathbf{j} & \mathbf{k}
           \end{bmatrix},\\
\mathcal{X}_{1}(:,:,1,1)=\begin{bmatrix}
             \mathbf{i-k} & \mathbf{0} \\
            \mathbf{0} & \mathbf{0}
           \end{bmatrix},\
\mathcal{X}_{1}(:,:,1,2)=\begin{bmatrix}
             \mathbf{j-k} & \mathbf{0} \\
            \mathbf{0} & \mathbf{0}
           \end{bmatrix},\
\mathcal{X}_{1}(:,:,2,1)=\begin{bmatrix}
             \mathbf{1+i} & \mathbf{0} \\
            \mathbf{0} & \mathbf{1-i}
           \end{bmatrix},\\
\mathcal{Y}_{1}(:,:,1,1)=\begin{bmatrix}
             \mathbf{0} & \mathbf{0} \\
            \mathbf{i-K} & \mathbf{0}
           \end{bmatrix},\
\mathcal{Y}_{1}(:,:,1,2)=\begin{bmatrix}
             \mathbf{0} & \mathbf{0} \\
            \mathbf{j} & \mathbf{j-k}
           \end{bmatrix},\
\mathcal{Y}_{1}(:,:,2,1)=\begin{bmatrix}
             \mathbf{0} & \mathbf{i} \\
            \mathbf{0} & \mathbf{2j}
           \end{bmatrix},\
 \mathcal{Y}_{1}(:,:,2,2)=\begin{bmatrix}
             \mathbf{0} & \mathbf{j} \\
            \mathbf{0} & \mathbf{2k}
           \end{bmatrix},\\
\mathcal{Z}_{2}(:,:,1,1)=\begin{bmatrix}
             \mathbf{k} & \mathbf{0} \\
            \mathbf{0} & \mathbf{i}
           \end{bmatrix},\
\mathcal{Z}_{2}(:,:,1,2)=\begin{bmatrix}
             \mathbf{j-1} & \mathbf{j} \\
            \mathbf{0} & \mathbf{1}
           \end{bmatrix},\
\mathcal{Z}_{2}(:,:,2,1)=\begin{bmatrix}
             \mathbf{1} & \mathbf{2} \\
            \mathbf{0} & \mathbf{i}
           \end{bmatrix},\
\mathcal{Z}_{2}(:,:,2,2)=\begin{bmatrix}
             \mathbf{1} & \mathbf{3} \\
            \mathbf{0} & \mathbf{j}
           \end{bmatrix},\\
\mathcal{X}_{2}(:,:,1,1)=\begin{bmatrix}
             \mathbf{-1} & \mathbf{-j+k} \\
            \mathbf{0} & \mathbf{0}
           \end{bmatrix},\
\mathcal{X}_{2}(:,:,1,2)=\begin{bmatrix}
             \mathbf{-3i} & \mathbf{0} \\
            \mathbf{5} & \mathbf{0}
           \end{bmatrix},\
\mathcal{X}_{2}(:,:,2,1)=\begin{bmatrix}
             \mathbf{2i} & \mathbf{-2} \\
            \mathbf{0} & \mathbf{0}
           \end{bmatrix},\\
\mathcal{Y}_{2}(:,:,1,1)=\begin{bmatrix}
             \mathbf{0} & \mathbf{-3} \\
            \mathbf{0} & \mathbf{i}
           \end{bmatrix},\
\mathcal{Y}_{2}(:,:,1,2)=\begin{bmatrix}
             \mathbf{0} & \mathbf{0} \\
            \mathbf{3} & \mathbf{0}
           \end{bmatrix},\
\mathcal{Y}_{2}(:,:,2,1)=\begin{bmatrix}
             \mathbf{0} & \mathbf{-4} \\
            \mathbf{0} & \mathbf{j}
           \end{bmatrix},\
 \mathcal{Y}_{2}(:,:,2,2)=\begin{bmatrix}
             \mathbf{2i} & \mathbf{0} \\
            \mathbf{2} & \mathbf{j}
           \end{bmatrix},\\
\mathcal{Z}_{3}(:,:,1,1)=\begin{bmatrix}
             \mathbf{0} & \mathbf{0} \\
            \mathbf{2i} & \mathbf{0}
           \end{bmatrix},\
\mathcal{Z}_{3}(:,:,1,2)=\begin{bmatrix}
             \mathbf{i-k} & \mathbf{0} \\
            \mathbf{0} & \mathbf{2i}
           \end{bmatrix},\
\mathcal{Z}_{3}(:,:,2,1)=\begin{bmatrix}
             \mathbf{0} & \mathbf{3i} \\
            \mathbf{0} & \mathbf{3j}
           \end{bmatrix},\
 \mathcal{Z}_{3}(:,:,2,2)=\begin{bmatrix}
             \mathbf{5i} & \mathbf{4j} \\
            \mathbf{0} & \mathbf{0}
           \end{bmatrix},\\          
\mathcal{X}_{3}(:,:,1,1)=\begin{bmatrix}
             \mathbf{0} & \mathbf{i} \\
            \mathbf{j} & \mathbf{0}
           \end{bmatrix},\
\mathcal{X}_{3}(:,:,1,2)=\begin{bmatrix}
             \mathbf{0} & \mathbf{2j} \\
            \mathbf{k} & \mathbf{0}
           \end{bmatrix},\
\mathcal{X}_{3}(:,:,2,1)=\begin{bmatrix}
             \mathbf{-i+k} & \mathbf{0} \\
            \mathbf{0} & \mathbf{j}
           \end{bmatrix},\
 \mathcal{X}_{3}(:,:,2,2)=\begin{bmatrix}
             \mathbf{j-k} & \mathbf{0} \\
            \mathbf{0} & \mathbf{i}
           \end{bmatrix},\\
\mathcal{Y}_{3}(:,:,1,1)=\begin{bmatrix}
             \mathbf{0} & \mathbf{i+j+k} \\
            \mathbf{0} & \mathbf{0}
           \end{bmatrix},\
\mathcal{Y}_{3}(:,:,1,2)=\begin{bmatrix}
             \mathbf{0} & \mathbf{0} \\
            \mathbf{i+j} & \mathbf{j+k}
           \end{bmatrix},\
\mathcal{Y}_{3}(:,:,2,1)=\begin{bmatrix}
             \mathbf{i+k} & \mathbf{j+k} \\
            \mathbf{0} & \mathbf{0}
           \end{bmatrix},\\
 \mathcal{Z}_{1}(:,:,2,2)=\begin{bmatrix}
             \mathbf{0} & \mathbf{0} \\
            \mathbf{2+k} & \mathbf{0}
           \end{bmatrix},\
 \mathcal{X}_{1}(:,:,2,2)=\begin{bmatrix}
             \mathbf{0} & \mathbf{0} \\
            \mathbf{0} & \mathbf{2j}
           \end{bmatrix},\
\mathcal{X}_{2}(:,:,2,2)=\begin{bmatrix}
             \mathbf{i} & \mathbf{0} \\
            \mathbf{4} & \mathbf{-k}
           \end{bmatrix},\
\mathcal{Y}_{3}(:,:,2,2)=\begin{bmatrix}
             \mathbf{0} & \mathbf{2i} \\
            \mathbf{0} & \mathbf{3j}
           \end{bmatrix}.
\end{gathered}              
\end{align*}
\end{small}
\end{example}
\begin{remark}
\label{system 22t1}
If we set $\mathcal{C}_{i}=\mathcal{B}_{i}=\mathcal{I}$ in \eqref{1.4aa} where $i=\overline{1,3},$ we obtain the Sylvester-like quaternion system of tensor equations \eqref{1.7aa}.
\end{remark}
\begin{remark}
\label{system 22t1R}
If we set $\mathcal{A}_{i}=\mathcal{D}_{i}=0$ in \eqref{1.7aa} where $i=\overline{1,3},$ we derive the  Sylvester-like quaternion system of tensor equations \eqref{1.5aa}.
\end{remark}
\begin{remark}
\label{system 22t1Rt}
If we set $\mathcal{G}_{i}=\mathcal{F}_{i}^{\eta^{*}}=0$,  $\mathcal{J}_{i}=\mathcal{H}_{i}^{\eta^{*}}=0$ and $\mathcal{E}_{i}=\mathcal{E}_{i}^{\eta^{*}}=0$  in \eqref{1.5aa} where $i=\overline{1,3},$ we  investigate $\eta$-Hermitian solution for \eqref{1.6aa}.
\end{remark}
In the following $Section$, we establish  the consistency conditions and the  general solution to \eqref{1.5aa}. In a direct implementation, we investigate some necessary and sufficient conditions for the existence of a common $\eta$-Hermitian solution of \eqref{1.6aa}.
\section{\textbf{Some implementations of the central system  $\mathbf{(1.4)}$}}
\begin{theorem}\label{system eeee1}
Consider the quaternion system of tensor equations  \eqref{1.5aa}, where
\begin{align*}
&\mathcal{F}_{4} \in \mathbb{H}^{I(N)\times J(N)},\ \mathcal{G}_{4} \in \mathbb{H}^{L(M)\times K(M)},\
 \mathcal{H}_{4} \in \mathbb{H}^{I(N)\times Q(N)},\ \mathcal{J}_{4} \in \mathbb{H}^{S(M)\times K(M)},\\
&\mathcal{E}_{4} \in \mathbb{H}^{I(N)\times K(M)},\ \mathcal{E}_{5} \in \mathbb{H}^{I(N)\times K(M)},\
\mathcal{F}_{i} \in \mathbb{H}^{A(N)\times J(N)},\ \mathcal{G}_{i} \in \mathbb{H}^{L(M)\times F(M)},\\
&\mathcal{H}_{i} \in \mathbb{H}^{A(N)\times P(N)},\  \mathcal{J}_{i} \in \mathbb{H}^{S(M)\times F(M)},\
\mathcal{E}_{i} \in \mathbb{H}^{A(N)\times F(M)}\ (i=\overline{1,3})
\end{align*} are given tensors over $\mathbb{H}$. Set
\begin{subequations}
\begin{align}
&\mathcal{\widehat{M}}_{i}=\mathcal{R}_{\mathcal{F}_{i}}*_{N}\mathcal{H}_{i},\
\mathcal{\widehat{N}}_{i}=\mathcal{J}_{i}*_{M}\mathcal{L}_{\mathcal{G}_{i}},\  \mathcal{\widehat{S}}_{i}=\mathcal{H}_{i}*_{N}\mathcal{L}_{\mathcal{\widehat{M}}_{i}},\ (i=\overline{1,3}),\\
\label{abc55}
&\mathcal{A}_{11}=\begin{bmatrix}\mathcal{L}_{\mathcal{F}_{4}} & -\mathcal{L}_{\mathcal{F}_{1}}\end{bmatrix},\
\mathcal{D}_{11}=\begin{bmatrix}\mathcal{R}_{\mathcal{G}_{4}} \\ -\mathcal{R}_{\mathcal{G}_{1}}\end{bmatrix},\
\mathcal{\widehat{A}}_{11}=\mathcal{F}_{1}^{\dagger}*_{N}\mathcal{\widehat{S}}_{1},\
\mathcal{\widehat{B}}_{11}=R_{\mathcal{\widehat{N}}_{1}}*_{M}\mathcal{J}_{1}*_{M}\mathcal{G}_{1}^{\dagger},\\
\label{abc66}
&\begin{array}{l}
 \mathcal{E}_{11}=\mathcal{F}_{1}^{\dagger}*_{N}\mathcal{E}_{1}*_{M}\mathcal{G}_{1}^{\dagger}
-\mathcal{F}_{1}^{\dagger}*_{N}\mathcal{H}_{1}*_{N}\mathcal{\widehat{M}}_{1}^{\dagger}*_{N}\mathcal{E}_{1}
*_{M}\mathcal{G}_{1}^{\dagger}
-\mathcal{F}_{1}^{\dagger}*_{N}\mathcal{\widehat{S}}_{1}*_{N}\mathcal{H}_{1}^{\dagger}*_{N}\mathcal{E}_{1}\\
\ \ \ \ \ *_{M}\mathcal{\widehat{N}}_{1}^{\dagger}*_{M}\mathcal{J}_{1}*_{M}\mathcal{G}_{1}^{\dagger}
-\mathcal{F}_{4}^{\dagger}*_{N}\mathcal{E}_{4}*_{M}\mathcal{G}_{4}^{\dagger},
 \end{array}\\
\label{abc77}
&\mathcal{A}_{22}=\begin{bmatrix}\mathcal{L}_{\mathcal{H}_{4}} & -\mathcal{L}_{\mathcal{\widehat{M}}_{3}}*_{N}\mathcal{L}_{\mathcal{\widehat{S}}_{3}}\end{bmatrix},\
\mathcal{D}_{22}=\begin{bmatrix}\mathcal{R}_{\mathcal{J}_{4}} \\ -\mathcal{R}_{\mathcal{J}_{3}}\end{bmatrix},\
\mathcal{\widehat{A}}_{22}=\mathcal{L}_{\mathcal{\widehat{M}}_{3}},\
\mathcal{\widehat{B}}_{22}=R_{\mathcal{\widehat{N}}_{3}},\\
\label{abc88}
&\begin{array}{l}
 \mathcal{E}_{22}=\mathcal{\widehat{M}}_{3}^{\dagger}*_{N}\mathcal{E}_{3}*_{M}\mathcal{J}_{3}^{\dagger}
 +\mathcal{\widehat{S}}_{3}^{\dagger}*_{N}\mathcal{\widehat{S}}_{3}*_{N}\mathcal{H}_{3}^{\dagger}
*_{N}\mathcal{E}_{3}*_{M}\mathcal{\widehat{N}}_{3}^{\dagger}
-\mathcal{H}_{4}^{\dagger}*_{N}\mathcal{E}_{5}*_{M}\mathcal{J}_{4}^{\dagger},
 \end{array}\\
 \label{abc99}
&\mathcal{\widehat{\widehat{A}}}_{ii}=\mathcal{R}_{\mathcal{A}_{ii}}*_{N}\mathcal{\widehat{A}}_{ii},\
\mathcal{\widehat{\widehat{B}}}_{ii}=\mathcal{\widehat{B}}_{ii}*_{M}\mathcal{L}_{\mathcal{D}_{ii}},\
\mathcal{\widehat{\widehat{E}}}_{ii}=\mathcal{R}_{\mathcal{A}_{ii}}*_{N}\mathcal{E}_{ii}*_{M}\mathcal{L}_{\mathcal{D}_{ii}},\ (i=1,2),\\
 \label{abc9911}
&\mathcal{\overline{A}}_{1}=\begin{bmatrix}-\mathcal{L}_{\mathcal{\widehat{M}}_{1}}*_{N}\mathcal{L}_{\mathcal{\widehat{S}}_{1}}& \mathcal{L}_{\mathcal{F}_{2}}\end{bmatrix},\ \mathcal{\overline{A}}_{2}=\begin{bmatrix}-\mathcal{L}_{\mathcal{\widehat{M}}_{2}}*_{N}\mathcal{L}_{\mathcal{\widehat{S}}_{2}}& \mathcal{L}_{\mathcal{F}_{3}}\end{bmatrix},\ \mathcal{\overline{F}}_{1}=\mathcal{F}^{\dagger}_{2}*_{N}\mathcal{\widehat{S}}_{2},\\
\label{abc9912}
&\mathcal{\overline{B}}_{1}=\begin{bmatrix}-\mathcal{R}_{\mathcal{J}_{1}}\\ \mathcal{R}_{\mathcal{G}_{2}}\end{bmatrix},\ \mathcal{\overline{B}}_{2}=\begin{bmatrix}-\mathcal{R}_{\mathcal{J}_{2}}\\ \mathcal{R}_{\mathcal{G}_{3}}\end{bmatrix},
\mathcal{\overline{F}}_{2}=\mathcal{F}^{\dagger}_{3}*_{N}\mathcal{\widehat{S}}_{3},\
\mathcal{\overline{G}}_{1}=\mathcal{J}_{2}*_{N}\mathcal{G}^{\dagger}_{2},\
\mathcal{\overline{G}}_{2}=\mathcal{J}_{3}*_{N}\mathcal{G}^{\dagger}_{3},\\
\label{abc99131}
&\mathcal{\overline{H}}_{1}=\mathcal{L}_{\mathcal{\widehat{M}}_{1}},\ \mathcal{\overline{J}}_{1}=\mathcal{R}_{\mathcal{\widehat{N}}_{1}},\
\mathcal{\overline{H}}_{2}=\mathcal{L}_{\mathcal{\widehat{M}}_{2}},\ \mathcal{\overline{J}}_{2}=\mathcal{R}_{\mathcal{\widehat{N}}_{2}},\\
\label{abc99132}
&\begin{array}{l}
 \mathcal{\overline{E}}_{1}=-\mathcal{\widehat{M}}_{1}^{\dagger}*_{N}\mathcal{E}_{1}*_{M}\mathcal{J}_{1}^{\dagger}
 -\mathcal{\widehat{S}}_{1}^{\dagger}*_{N}\mathcal{\widehat{S}}_{1}*_{N}\mathcal{H}_{1}^{\dagger}
*_{N}\mathcal{E}_{1}*_{M}\mathcal{\widehat{N}}_{1}^{\dagger}
+\mathcal{F}_{2}^{\dagger}*_{N}\mathcal{E}_{2}*_{M}\mathcal{G}_{2}^{\dagger}
-\mathcal{F}_{2}^{\dagger}\\
\ \ \ \ \ \ \ \ \ \ \ *_{N}\mathcal{H}_{2}*_{N}\mathcal{\widehat{M}}_{2}^{\dagger}*_{N}\mathcal{E}_{2}
*_{M}\mathcal{G}_{2}^{\dagger}
-\mathcal{F}_{2}^{\dagger}*_{N}\mathcal{\widehat{S}}_{2}*_{N}\mathcal{H}_{2}^{\dagger}*_{N}\mathcal{E}_{2}
*_{M}\mathcal{\widehat{N}}_{2}^{\dagger}*_{M}\mathcal{J}_{2}*_{M}\mathcal{G}_{2}^{\dagger},
 \end{array}\\
\label{abc991322}
&\begin{array}{l}
 \mathcal{\overline{E}}_{2}=-\mathcal{\widehat{M}}_{2}^{\dagger}*_{N}\mathcal{E}_{2}*_{M}\mathcal{J}_{2}^{\dagger}
 -\mathcal{\widehat{S}}_{2}^{\dagger}*_{N}\mathcal{\widehat{S}}_{2}*_{N}\mathcal{H}_{2}^{\dagger}
*_{N}\mathcal{E}_{2}*_{M}\mathcal{\widehat{N}}_{2}^{\dagger}
+\mathcal{F}_{2}^{\dagger}*_{N}\mathcal{E}_{2}*_{M}\mathcal{G}_{2}^{\dagger}
-\mathcal{F}_{3}^{\dagger}\\
\ \ \ \ \ \ \ \ \ \ \ *_{N}\mathcal{H}_{3}*_{N}\mathcal{\widehat{M}}_{3}^{\dagger}*_{N}\mathcal{E}_{3}
*_{M}\mathcal{G}_{3}^{\dagger}
-\mathcal{F}_{3}^{\dagger}*_{N}\mathcal{\widehat{S}}_{3}*_{N}\mathcal{H}_{3}^{\dagger}*_{N}\mathcal{E}_{3}
*_{M}\mathcal{\widehat{N}}_{3}^{\dagger}*_{M}\mathcal{J}_{3}*_{M}\mathcal{G}_{3}^{\dagger},
 \end{array}\\
\label{abxc22}
&\mathcal{\overline{F}}_{ii}=\mathcal{R}_{\mathcal{\overline{A}}_{i}}*_{N}\mathcal{\overline{F}}_{i},\ \mathcal{\overline{G}}_{ii}=\mathcal{\overline{G}}_{i}*_{M}\mathcal{L}_{\mathcal{\overline{B}}_{i}},\ \mathcal{\overline{H}}_{ii}=\mathcal{R}_{\mathcal{\overline{A}}_{i}}*_{N}\mathcal{\overline{H}}_{i},\ \mathcal{\overline{J}}_{ii}=\mathcal{\overline{J}}_{i}*_{M}\mathcal{L}_{\mathcal{\overline{B}}_{i}},\\
\label{abxc33}
&\mathcal{\overline{E}}_{ii}=\mathcal{R}_{\mathcal{\overline{A}}_{i}}*_{N}\mathcal{\overline{E}}_{i}*_{M}\mathcal{L}_{\mathcal{\overline{B}}_{i}},\ \mathcal{\overline{M}}_{ii}=\mathcal{R}_{\mathcal{\overline{F}}_{ii}}*_{N}\mathcal{\overline{H}}_{ii},\ \mathcal{\overline{N}}_{ii}=\mathcal{\overline{J}}_{ii}*_{M}\mathcal{L}_{\mathcal{\overline{G}}_{ii}},\ \mathcal{\overline{S}}_{ii}=\mathcal{\overline{H}}_{ii}*_{N}\mathcal{L}_{\mathcal{\overline{M}}_{ii}},\\
\label{abcD55}
&\mathcal{\overline{\overline{A}}}_{1}=\begin{bmatrix}\mathcal{L}_{\mathcal{\overline{F}}_{11}} &
 -\mathcal{L}_{\mathcal{\overline{M}}_{22}}*_{N}\mathcal{L}_{\mathcal{\overline{S}}_{22}}\end{bmatrix},\
\mathcal{\overline{\overline{B}}}_{1}=\begin{bmatrix}\mathcal{R}_{\mathcal{\overline{G}}_{11}} \\ -\mathcal{R}_{\mathcal{\overline{J}}_{11}}\end{bmatrix},\
\mathcal{\overline{\overline{F}}}_{1}=\mathcal{\overline{F}}_{11}^{\dagger}*_{N}\mathcal{\overline{S}}_{11},\\
\label{abcDE55}
&\mathcal{\overline{\overline{G}}}_{1}=R_{\mathcal{\overline{N}}_{11}}*_{M}\mathcal{\overline{J}}_{11}*_{M}\mathcal{\overline{G}}_{11}^{\dagger},\
\mathcal{\overline{\overline{H}}}_{1}=\mathcal{L}_{\mathcal{\overline{M}}_{22}},\
\mathcal{\overline{\overline{J}}}_{1}=\mathcal{R}_{\mathcal{\overline{N}}_{22}},\\
\label{abcD66}
&\begin{array}{l}
 \mathcal{\overline{\overline{E}}}_{1}=\mathcal{\overline{F}}_{11}^{\dagger}*_{N}\mathcal{\overline{E}}_{11}*_{M}\mathcal{\overline{G}}_{11}^{\dagger}
-\mathcal{\overline{F}}_{11}^{\dagger}*_{N}\mathcal{\overline{H}}_{11}*_{N}\mathcal{\overline{M}}_{11}^{\dagger}*_{N}\mathcal{\overline{E}}_{11}
*_{M}\mathcal{\overline{G}}_{11}^{\dagger}
-\mathcal{\overline{F}}_{11}^{\dagger}*_{N}\mathcal{\overline{S}}_{11}*_{N}\mathcal{\overline{H}}_{11}^{\dagger}\\ \ \ \ \ \ \ \ \ *_{N}\mathcal{\overline{E}}_{11}*_{M}\mathcal{\overline{N}}_{11}^{\dagger}*_{M}\mathcal{\overline{J}}_{11}*_{M}\mathcal{\overline{G}}_{11}^{\dagger}-\mathcal{\overline{M}}_{22}^{\dagger}*_{N}\mathcal{\overline{E}}_{22}*_{M}\mathcal{\overline{J}}_{22}^{\dagger}
 -\mathcal{\overline{S}}_{22}^{\dagger}*_{N}\mathcal{\overline{S}}_{22}*_{N}\mathcal{\overline{H}}_{22}^{\dagger}\\
\ \  \ \ \ \ \ \ \ *_{N}\mathcal{\overline{E}}_{22}*_{M}\mathcal{\overline{N}}_{22}^{\dagger}
 \end{array}\\
\label{abKxc22}
&\mathcal{\overline{\overline{F}}}_{11}=\mathcal{R}_{\mathcal{\overline{\overline{A}}}_{1}}*_{N}\mathcal{\overline{\overline{F}}}_{1},\ \mathcal{\overline{\overline{G}}}_{11}=\mathcal{\overline{\overline{G}}}_{1}*_{M}\mathcal{L}_{\mathcal{\overline{\overline{B}}}_{1}},\ \mathcal{\overline{\overline{H}}}_{11}=\mathcal{R}_{\mathcal{\overline{\overline{A}}}_{1}}*_{N}\mathcal{\overline{\overline{H}}}_{1},\ \mathcal{\overline{\overline{J}}}_{11}=\mathcal{\overline{\overline{J}}}_{1}*_{M}\mathcal{L}_{\mathcal{\overline{\overline{B}}}_{1}},
\end{align}
\begin{align}
\label{abKKxAc33}
&\mathcal{\overline{\overline{E}}}_{11}=\mathcal{R}_{\mathcal{\overline{\overline{A}}}_{1}}*_{N}\mathcal{\overline{\overline{E}}}_{1}
*_{M}\mathcal{L}_{\mathcal{\overline{\overline{B}}}_{1}},\ \mathcal{\overline{\overline{M}}}_{11}=\mathcal{R}_{\mathcal{\overline{\overline{F}}}_{11}}*_{N}\mathcal{\overline{\overline{H}}}_{11},\ \mathcal{\overline{\overline{N}}}_{11}=\mathcal{\overline{\overline{J}}}_{11}*_{M}\mathcal{L}_{\mathcal{\overline{\overline{G}}}_{11}},\\
\label{abc99911}
&\mathcal{\overline{\overline{S}}}_{11}=\mathcal{\overline{\overline{H}}}_{11}*_{N}\mathcal{L}_{\mathcal{\overline{\overline{M}}}_{11}},\
\mathcal{\widetilde{A}}_{1}=\begin{bmatrix}\mathcal{L}_{\mathcal{\overline{M}}_{11}}*_{N}\mathcal{L}_{\mathcal{\overline{S}}_{11}}& -\mathcal{L}_{\mathcal{\widehat{\widehat{A}}}_{11}}\end{bmatrix},\ \mathcal{\widetilde{A}}_{2}=\begin{bmatrix}\mathcal{L}_{\mathcal{\widehat{\widehat{A}}}_{22}}& -\mathcal{L}_{\mathcal{\overline{F}}_{22}}\end{bmatrix},\\
\label{abc99912}
&\mathcal{\widetilde{B}}_{1}=\begin{bmatrix}\mathcal{R}_{\mathcal{\overline{J}}_{11}}\\-\mathcal{R}_{\mathcal{\widehat{\widehat{B}}}_{11}}\end{bmatrix},\ \mathcal{\widetilde{B}}_{2}=\begin{bmatrix}\mathcal{R}_{\mathcal{\widehat{\widehat{B}}}_{22}}\\ -\mathcal{R}_{\mathcal{\overline{G}}_{22}}\end{bmatrix},
\mathcal{\widetilde{C}}_{1}=\mathcal{L}_{\mathcal{\overline{M}}_{11}}
\mathcal{\widetilde{D}}_{1}=\mathcal{R}_{\mathcal{\overline{N}}_{11}},\\
\label{abc999131}
&\mathcal{\widetilde{C}}_{2}=\mathcal{\overline{F}}^{\dagger}_{22}*_{N}\mathcal{\overline{S}}_{22},\
\mathcal{\widetilde{D}}_{2}=\mathcal{R}_{\mathcal{\overline{N}}_{22}}*_{M}\mathcal{\overline{J}}_{22}*_{M}\mathcal{\overline{G}}^{\dagger}_{22},\\
\label{abcF131}
&\mathcal{\widetilde{E}}_{1}=\widehat{\widehat{\mathcal{A}}}_{11}^{\dagger}*_{N}\widehat{\widehat{\mathcal{E}}}_{11}
*_{M}\widehat{\widehat{\mathcal{B}}}_{11}^{\dagger}-\mathcal{\overline{M}}_{11}^{\dagger}*_{N}\mathcal{\overline{E}}_{11}*_{M}\mathcal{\overline{J}}_{11}^{\dagger}
 -\mathcal{\overline{S}}_{11}^{\dagger}*_{N}\mathcal{\overline{S}}_{11}*_{N}\mathcal{\overline{H}}_{11}^{\dagger}
*_{N}\mathcal{\overline{E}}_{11}*_{M}\mathcal{\overline{N}}_{11}^{\dagger},\\
\label{abcD66}
&\begin{array}{l}
\mathcal{\widetilde{E}}_{2}=\mathcal{\overline{F}}_{22}^{\dagger}*_{N}\mathcal{\overline{E}}_{22}*_{M}\mathcal{\overline{G}}_{22}^{\dagger}
-\mathcal{\overline{F}}_{22}^{\dagger}*_{N}\mathcal{\overline{H}}_{22}*_{N}\mathcal{\overline{M}}_{22}^{\dagger}*_{N}\mathcal{\overline{E}}_{22}
*_{M}\mathcal{\overline{G}}_{22}^{\dagger}
-\mathcal{\overline{F}}_{22}^{\dagger}*_{N}\mathcal{\overline{S}}_{22}*_{N}\mathcal{\overline{H}}_{22}^{\dagger}\\
\ \ \ \ \ *_{N}\mathcal{\overline{E}}_{22}*_{M}\mathcal{\overline{N}}_{22}^{\dagger}*_{M}\mathcal{\overline{J}}_{22}*_{M}\mathcal{\overline{G}}_{22}^{\dagger}
-\widehat{\widehat{\mathcal{A}}}_{22}^{\dagger}*_{N}\widehat{\widehat{\mathcal{E}}}_{22}
*_{M}\widehat{\widehat{\mathcal{B}}}_{22}^{\dagger},\\
 \end{array}\\
\label{abc999WW11}
& \mathcal{\widetilde{F}}_{1}=\begin{bmatrix}\mathcal{L}_{\mathcal{\widetilde{C}}_{11}}& -\mathcal{L}_{\mathcal{\overline{\overline{F}}}_{11}}\end{bmatrix},\
\mathcal{\widetilde{F}}_{2}=\begin{bmatrix}\mathcal{L}_{\mathcal{\overline{\overline{M}}}_{11}}*_{N}\mathcal{L}_{\mathcal{\overline{\overline{S}}}_{11}}& -\mathcal{L}_{\mathcal{\widetilde{C}}_{22}}\end{bmatrix},\
\mathcal{\widetilde{H}}_{1}=\mathcal{\overline{\overline{F}}}_{11}^{\dagger}*_{N}\mathcal{\overline{\overline{S}}}_{11},\\
\label{abc999WWW12}
&\mathcal{\widetilde{J}}_{1}=\mathcal{R}_{\mathcal{\overline{\overline{N}}}_{11}}
*_{M}\mathcal{\overline{\overline{J}}}_{11}*_{M}\mathcal{\overline{\overline{G}}}_{11}^{\dagger}\
\mathcal{\widetilde{G}}_{1}=\begin{bmatrix}\mathcal{R}_{\mathcal{\widetilde{D}}_{11}}\\-\mathcal{R}_{\mathcal{\overline{\overline{G}}}_{11}}\end{bmatrix},\ \mathcal{\widetilde{G}}_{2}=\begin{bmatrix}\mathcal{R}_{\mathcal{\overline{\overline{J}}}_{11}}\\ -\mathcal{R}_{\mathcal{\widetilde{D}}_{22}}\end{bmatrix},
\mathcal{\widetilde{H}}_{2}=\mathcal{L}_{\mathcal{\overline{\overline{M}}}_{11}},\
\mathcal{\widetilde{J}}_{2}=\mathcal{R}_{\mathcal{\overline{\overline{N}}}_{11}},\\
\label{abcDS66}
&\begin{array}{l}
 \mathcal{\widetilde{\widetilde{E}}}_{1}=\mathcal{\overline{\overline{F}}}_{11}^{\dagger}*_{N}\mathcal{\overline{\overline{E}}}_{11}*_{M}\mathcal{\overline{\overline{G}}}_{11}^{\dagger}
-\mathcal{\overline{\overline{F}}}_{11}^{\dagger}*_{N}\mathcal{\overline{\overline{H}}}_{11}*_{N}\mathcal{\overline{\overline{M}}}_{11}^{\dagger}
*_{N}\mathcal{\overline{\overline{E}}}_{11}
*_{M}\mathcal{\overline{\overline{G}}}_{11}^{\dagger}
-\mathcal{\overline{\overline{F}}}_{11}^{\dagger}*_{N}\mathcal{\overline{\overline{S}}}_{11}\\
\ \ \ \ \ *_{N}\mathcal{\overline{\overline{H}}}_{11}^{\dagger}*_{N}\mathcal{\overline{\overline{E}}}_{11}*_{M}\mathcal{\overline{\overline{N}}}_{11}^{\dagger}*_{M}\mathcal{\overline{\overline{J}}}_{11}
*_{M}\mathcal{\overline{\overline{G}}}_{11}^{\dagger}-\mathcal{\widetilde{C}}_{11}^{\dagger}*_{N}\mathcal{\widetilde{E}}_{11}
*_{M}\mathcal{\widetilde{D}}_{11}^{\dagger},
 \end{array}
\end{align}
\end{subequations}
\begin{subequations}
\begin{align}
\label{abcAAD66}
&\mathcal{\widetilde{\widetilde{E}}}_{2}= \mathcal{\widetilde{C}}_{22}^{\dagger}*_{N}\mathcal{\widetilde{E}}_{22}
*_{M}\mathcal{\widetilde{D}}_{22}^{\dagger}
-\mathcal{\overline{\overline{M}}}_{11}^{\dagger}*_{N}\mathcal{\overline{\overline{E}}}_{11}*_{M}\mathcal{\overline{\overline{J}}}_{11}^{\dagger}
 -\mathcal{\overline{\overline{S}}}_{11}^{\dagger}*_{N}\mathcal{\overline{\overline{S}}}_{11}*_{N}\mathcal{\overline{\overline{H}}}_{11}^{\dagger}
*_{N}\mathcal{\overline{\overline{E}}}_{11}*_{M}\mathcal{\overline{\overline{N}}}_{11}^{\dagger},\\
\label{abcDGG66}
&\mathcal{\widetilde{H}}_{11}=\mathcal{R}_{\mathcal{\widetilde{F}}_{1}}*_{N}\mathcal{\widetilde{H}}_{1},\
\mathcal{\widetilde{H}}_{22}=\mathcal{R}_{\mathcal{\widetilde{F}}_{2}}*_{N}\mathcal{\widetilde{H}}_{2},\
\mathcal{\widetilde{J}}_{11}=\mathcal{\widetilde{J}}_{1}*_{M}\mathcal{L}_{\mathcal{\widetilde{G}}_{1}},\
\mathcal{\widetilde{J}}_{22}=\mathcal{\widetilde{J}}_{2}*_{M}\mathcal{L}_{\mathcal{\widetilde{G}}_{2}},\\
\label{abcHHD66}
&\mathcal{\widetilde{\widetilde{E}}}_{11}=\mathcal{R}_{\mathcal{\widetilde{F}}_{1}}*_{N}\mathcal{\widetilde{\widetilde{E}}}_{1}
*_{M}\mathcal{L}_{\mathcal{\widetilde{G}}_{1}},\
\mathcal{\widetilde{\widetilde{E}}}_{22}=\mathcal{R}_{\mathcal{\widetilde{F}}_{2}}*_{N}\mathcal{\widetilde{\widetilde{E}}}_{2}
*_{M}\mathcal{L}_{\mathcal{\widetilde{G}}_{2}},\
\mathcal{\widetilde{A}}=\begin{bmatrix}\mathcal{L}_{\mathcal{\widetilde{H}}_{11}}& -\mathcal{L}_{\mathcal{\widetilde{H}}_{22}}\end{bmatrix},\\
\label{abcHHRD66}
&\mathcal{\widetilde{B}}=\begin{bmatrix}\mathcal{R}_{\mathcal{\widetilde{J}}_{11}}& -\mathcal{R}_{\mathcal{\widetilde{J}}_{22}}\end{bmatrix},\
\mathcal{\widetilde{E}}=\mathcal{\widetilde{H}}_{22}^{\dagger}*_{N}\mathcal{\widetilde{\widetilde{E}}}_{22}
*_{M}\mathcal{\widetilde{J}}_{22}^{\dagger}-
\mathcal{\widetilde{H}}_{11}^{\dagger}*_{N}\mathcal{\widetilde{\widetilde{E}}}_{11}
*_{M}\mathcal{\widetilde{J}}_{11}^{\dagger}.
\end{align}
\end{subequations}
Then the system \eqref{1.4aa} is consistent if and only if
\begin{align}
\label{b22}
&\mathcal{R}_{\mathcal{\widehat{M}}_{i}}*_{N}\mathcal{R}_{\mathcal{F}_{i}}*_{N}\mathcal{E}_{i}=0,\ \mathcal{E}_{i}*_{M}\mathcal{L}_{\mathcal{G}_{i}}*_{M}\mathcal{L}_{\mathcal{\widehat{N}}_{i}}=0,\\
\label{b33}
&\mathcal{R}_{\mathcal{F}_{i}}*_{N}\mathcal{E}_{i}*_{M}\mathcal{L}_{\mathcal{J}_{i}}=0,\
\mathcal{R}_{\mathcal{H}_{i}}*_{N}\mathcal{E}_{i}*_{M}\mathcal{L}_{\mathcal{G}_{i}}=0,\ (i=\overline{1,3}),\\
\label{b44}
&\mathcal{R}_{\mathcal{F}_{4}}*_{N}\mathcal{E}_{4}=0,\ \mathcal{E}_{4}*_{M}\mathcal{L}_{\mathcal{G}_{4}}=0,\
\mathcal{R}_{\mathcal{H}_{4}}*_{N}\mathcal{E}_{5}=0,\ \mathcal{E}_{5}*_{M}\mathcal{L}_{\mathcal{J}_{4}}=0,\\
\label{b55}
&\mathcal{R}_{\mathcal{\widehat{\widehat{A}}}_{kk}}*_{N}\mathcal{\widehat{\widehat{E}}}_{kk}=0,\
\mathcal{\widehat{\widehat{E}}}_{kk}*_{M}\mathcal{L}_{\mathcal{\widehat{\widehat{B}}}_{kk}}=0,\\
\label{b66}
&\mathcal{R}_{\mathcal{\overline{M}}_{kk}}*_{N}\mathcal{R}_{\mathcal{\overline{F}}_{kk}}*_{N}\mathcal{\overline{E}}_{kk}=0,\ \mathcal{\overline{E}}_{kk}*_{M}\mathcal{L}_{\mathcal{\overline{G}}_{kk}}*_{M}\mathcal{L}_{\mathcal{\overline{N}}_{kk}}=0,\\
\label{b77}
&\mathcal{R}_{\mathcal{\overline{F}}_{kk}}*_{N}\mathcal{\overline{E}}_{kk}*_{M}\mathcal{L}_{\mathcal{\overline{J}}_{kk}}=0,\
\mathcal{R}_{\mathcal{\overline{H}}_{kk}}*_{N}\mathcal{\overline{E}}_{kk}*_{M}\mathcal{L}_{\mathcal{\overline{G}}_{kk}}=0,\ (k=1,2),\\
\label{b88}
&\mathcal{R}_{\mathcal{\overline{\overline{M}}}_{11}}*_{N}\mathcal{R}_{\mathcal{\overline{\overline{F}}}_{11}}*_{N}\mathcal{\overline{\overline{E}}}_{11}=0,\ \mathcal{\overline{\overline{E}}}_{11}*_{M}\mathcal{L}_{\mathcal{\overline{\overline{G}}}_{11}}*_{M}\mathcal{L}_{\mathcal{\overline{\overline{N}}}_{11}}=0,\\
\label{b99}
&\mathcal{R}_{\mathcal{\overline{\overline{F}}}_{11}}*_{N}\mathcal{\overline{\overline{E}}}_{11}*_{M}\mathcal{L}_{\mathcal{\overline{\overline{J}}}_{11}}=0,\
\mathcal{R}_{\mathcal{\overline{\overline{H}}}_{11}}*_{N}\mathcal{\overline{\overline{E}}}_{11}*_{M}\mathcal{L}_{\mathcal{\overline{\overline{G}}}_{11}}=0,\\
\label{bb111}
&\mathcal{R}_{\mathcal{\widetilde{C}}_{jj}}*_{N}\mathcal{\widetilde{E}}_{jj}=0,\
\mathcal{\widetilde{E}}_{jj}*_{M}\mathcal{L}_{\mathcal{\widetilde{D}}_{jj}}=0\ (j=1,2),\\
\label{bb1121}
&\mathcal{R}_{\mathcal{\widetilde{H}}_{ll}}*_{N}\mathcal{\widetilde{\widetilde{E}}}_{ll}=0,\
\mathcal{\widetilde{\widetilde{E}}}_{ll}*_{M}\mathcal{L}_{\mathcal{\widetilde{G}}_{ll}}=0,\
\mathcal{R}_{\mathcal{\widetilde{A}}}*_{N}\mathcal{\widetilde{E}}*_{M}\mathcal{L}_{\mathcal{\widetilde{B}}},\ (l=1,2),\\
\label{bb122121}
&\mathcal{R}_{\mathcal{\widetilde{A}}}*_{N}\mathcal{\widetilde{E}}*_{M}\mathcal{L}_{\mathcal{\widetilde{B}}}=0.
\end{align}
Under these conditions, the general solution to system \eqref{1.4aa} can be expressed as follows:
\begin{align}
\label{sssa}
&\mathcal{Z}_{1}=\mathcal{F}_{4}^{\dagger}*_{N}\mathcal{E}_{4}*_{M}\mathcal{G}_{4}^{\dagger}
+\mathcal{L}_{\mathcal{F}_{4}}*_{N}\mathcal{W}_{1}
+\mathcal{W}_{2}*_{M}\mathcal{R}_{\mathcal{G}_{4}},\\
&\mathcal{Z}_{4}=\mathcal{H}_{4}^{\dagger}*_{N}\mathcal{E}_{5}*_{M}\mathcal{J}_{4}^{\dagger}
+\mathcal{L}_{\mathcal{H}_{4}}*_{N}\mathcal{\acute{W}}_{1}
+\mathcal{W}_{3}*_{M}\mathcal{R}_{\mathcal{J}_{4}}, \\
&\begin{array}{l}
 \mathcal{Z}_{2}=\mathcal{\widehat{M}}_{1}^{\dagger}*_{N}\mathcal{\widehat{E}}_{1}*_{M}\mathcal{J}_{1}^{\dagger}
 +\mathcal{\widehat{S}}_{1}^{\dagger}*_{N}\mathcal{\widehat{S}}_{1}*_{N}\mathcal{H}_{1}^{\dagger}
*_{N}\mathcal{E}_{1}*_{M}\mathcal{\widehat{N}}_{1}^{\dagger}+\mathcal{L}_{\mathcal{\widehat{M}}_{1}}
*_{N}\mathcal{L}_{\mathcal{\widehat{S}}_{1}}*_{N}\mathcal{\widehat{U}}_{1}\\
\ \ \ \ \ \ +\mathcal{L}_{\mathcal{\widehat{M}}_{1}}*_{N}\mathcal{\widehat{U}}_{2}*_{M}\mathcal{R}_{\mathcal{\widehat{N}}_{1}}
+\mathcal{\widehat{U}}_{3}*_{M}\mathcal{R}_{\mathcal{J}_{1}},
 \end{array}\\
&\begin{array}{l}
or\ \mathcal{Z}_{2}=\mathcal{F}_{2}^{\dagger}*_{N}\mathcal{E}_{2}*_{M}\mathcal{G}_{2}^{\dagger}
-\mathcal{F}_{2}^{\dagger}*_{N}\mathcal{H}_{2}*_{N}\mathcal{\widehat{M}}_{2}^{\dagger}*_{N}\mathcal{E}_{2}
*_{M}\mathcal{G}_{2}^{\dagger}
-\mathcal{F}_{2}^{\dagger}*_{N}\mathcal{\widehat{S}}_{2}*_{N}\mathcal{H}_{2}^{\dagger}*_{N}\mathcal{E}_{2}\\
\ \ \ \ \ *_{M}\mathcal{\widehat{N}}_{2}^{\dagger}*_{M}\mathcal{J}_{2}*_{M}\mathcal{G}_{2}^{\dagger}-\mathcal{F}_{2}^{\dagger}
*_{N}\mathcal{\widehat{S}}_{2}
*_{N}\mathcal{\widehat{V}}_{2}*_{M}\mathcal{R}_{\mathcal{\widehat{N}}_{2}}*_{M}\mathcal{J}_{2}*_{M}\mathcal{G}_{2}^{\dagger}
+\mathcal{L}_{\mathcal{F}_{2}}*_{N}\mathcal{\widehat{V}}_{4}\\
\ \ \ \ \ +\mathcal{\widehat{V}}_{5}*_{M}\mathcal{R}_{\mathcal{F}_{2}},
 \end{array}\\
&\begin{array}{l}
 \mathcal{Z}_{3}=\mathcal{\widehat{M}}_{2}^{\dagger}*_{N}\mathcal{E}_{2}*_{M}\mathcal{J}_{2}^{\dagger}
 +\mathcal{\widehat{S}}_{2}^{\dagger}*_{N}\mathcal{\widehat{S}}_{2}*_{N}\mathcal{H}_{2}^{\dagger}
*_{N}\mathcal{E}_{2}*_{M}\mathcal{\widehat{N}}_{2}^{\dagger}+\mathcal{L}_{\mathcal{\widehat{M}}_{2}}
*_{N}\mathcal{L}_{\mathcal{\widehat{S}}_{2}}*_{N}\mathcal{\widehat{V}}_{1}\\
\ \ \ \ \ \ +\mathcal{L}_{\mathcal{\widehat{M}}_{2}}*_{N}\mathcal{\widehat{V}}_{2}*_{M}\mathcal{R}_{\mathcal{\widehat{N}}_{2}}
+\mathcal{\widehat{V}}_{3}*_{M}\mathcal{R}_{\mathcal{J}_{2}},
 \end{array}\\
&\begin{array}{l}
or\  \mathcal{Z}_{3}=\mathcal{F}_{3}^{\dagger}*_{N}\mathcal{E}_{3}*_{M}\mathcal{G}_{3}^{\dagger}
-\mathcal{F}_{3}^{\dagger}*_{N}\mathcal{H}_{3}*_{N}\mathcal{\widehat{M}}_{3}^{\dagger}*_{N}\mathcal{E}_{3}
*_{M}\mathcal{G}_{3}^{\dagger}
-\mathcal{F}_{3}^{\dagger}*_{N}\mathcal{\widehat{S}}_{3}*_{N}\mathcal{H}_{3}^{\dagger}*_{N}\mathcal{E}_{3}\\
\ \ \ \ \ *_{M}\mathcal{\widehat{N}}_{3}^{\dagger}*_{M}\mathcal{J}_{3}*_{M}\mathcal{G}_{3}^{\dagger}-\mathcal{F}_{3}^{\dagger}
*_{N}\mathcal{\widehat{S}}_{3}
*_{N}\mathcal{\widehat{K}}_{2}*_{M}\mathcal{R}_{\mathcal{\widehat{N}}_{3}}*_{M}\mathcal{J}_{3}*_{M}\mathcal{G}_{3}^{\dagger}
+\mathcal{L}_{\mathcal{F}_{3}}*_{N}\mathcal{\widehat{K}}_{4}\\
\ \ \ \ \ +\mathcal{\widehat{K}}_{5}*_{M}\mathcal{R}_{\mathcal{G}_{3}},\ (i=\overline{1,3}).
 \end{array}
\end{align}
Where the arbitrary tensors $\mathcal{W}_{j}$, $\mathcal{\widehat{V}}_{i}$, $\mathcal{\widehat{U}}_{j}$, $\mathcal{\widehat{K}}_{k}$ and $\mathcal{\acute{W}}_{1}$ $(j=\overline{1,3},\ i=\overline{1,5},\ k\in \{2,4,5\})$ can be reduced by \eqref{Xyqs1}-\eqref{XZyw996}.
\begin{proof} See $Remark$ (\ref{system 22t1})-$Remark$ (\ref{system 22t1R}).
\end{proof}
\end{theorem}
\begin{corollary}
\label{system eee1}
Consider the quaternion system of tensor equations  \eqref{1.5aa}, where
\begin{align*}
&\mathcal{F}_{4} \in \mathbb{H}^{I(N)\times J(N)},\
 \mathcal{H}_{4} \in \mathbb{H}^{I(N)\times Q(N)},\
\mathcal{E}_{4} \in \mathbb{H}^{I(N)\times I(N)},\ \mathcal{E}_{5} \in \mathbb{H}^{I(N)\times I(N)},\\
&\mathcal{F}_{i} \in \mathbb{H}^{A(N)\times J(N)},\ \mathcal{H}_{i} \in \mathbb{H}^{A(N)\times I(N)},\
\mathcal{E}_{i} \in \mathbb{H}^{A(N)\times A(N)}\ (i=\overline{1,3})
\end{align*} are given tensors over $\mathbb{H}$. Set
\begin{subequations}
\begin{align}
\label{abc5522}
&\mathcal{\widehat{M}}_{i}=\mathcal{R}_{\mathcal{F}_{i}}*_{N}\mathcal{H}_{i},\
\mathcal{\widehat{N}}_{i}=(\mathcal{\widehat{M}}_{i})^{\eta^{*}},\  \mathcal{\widehat{S}}_{i}=\mathcal{H}_{i}*_{N}\mathcal{L}_{\mathcal{\widehat{M}}_{i}},\ (i=\overline{1,3})\ \mathcal{A}_{11}=\begin{bmatrix}\mathcal{L}_{\mathcal{F}_{4}} & -\mathcal{L}_{\mathcal{F}_{1}}\end{bmatrix},\\
\label{abc55}
&\mathcal{D}_{11}=\begin{bmatrix}\mathcal{R}_{\mathcal{F}_{4}^{\eta^{*}}} \\ -\mathcal{R}_{\mathcal{F}_{1}^{\eta^{*}}}\end{bmatrix},\
\mathcal{\widehat{A}}_{11}=\mathcal{F}_{1}^{\dagger}*_{N}\mathcal{\widehat{S}}_{1},\
\mathcal{\widehat{B}}_{11}=R_{\mathcal{\widehat{N}}_{1}}*_{N}\mathcal{H}_{1}^{\eta^{*}}*_{N}(\mathcal{F}_{1}^{\eta^{*}})^{\dagger},\\
\label{abc66}
&\begin{array}{l}
 \mathcal{E}_{11}=\mathcal{F}_{1}^{\dagger}*_{N}\mathcal{E}_{1}*_{N}(\mathcal{F}_{1}^{\eta^{*}})^{\dagger}
-\mathcal{F}_{1}^{\dagger}*_{N}\mathcal{H}_{1}*_{N}\mathcal{\widehat{M}}_{1}^{\dagger}*_{N}\mathcal{E}_{1}
*_{N}(\mathcal{F}_{1}^{\eta^{*}})^{\dagger}
-\mathcal{F}_{1}^{\dagger}*_{N}\mathcal{\widehat{S}}_{1}\\
\ \ \ \ \ \ \ \ \ \ \ \ \ \ \ *_{N}\mathcal{H}_{1}^{\dagger}*_{N}\mathcal{E}_{1}*_{N}\mathcal{\widehat{N}}_{1}^{\dagger}*_{M}\mathcal{H}_{1}^{\eta^{*}}*_{N}(\mathcal{F}_{1}^{\eta^{*}})^{\dagger}
-\mathcal{F}_{4}^{\dagger}*_{N}\mathcal{E}_{4}*_{N}(\mathcal{F}_{4}^{\eta^{*}})^{\dagger},
 \end{array}\\
\label{abc77}
&\mathcal{A}_{22}=\begin{bmatrix}\mathcal{L}_{\mathcal{H}_{4}} & -\mathcal{L}_{\mathcal{\widehat{M}}_{3}}*_{N}\mathcal{L}_{\mathcal{\widehat{S}}_{3}}\end{bmatrix},\
\mathcal{D}_{22}=\begin{bmatrix}\mathcal{R}_{\mathcal{H}_{4}^{\eta^{*}}} \\ -\mathcal{R}_{\mathcal{H}_{3}^{\eta^{*}}}\end{bmatrix},\
\mathcal{\widehat{A}}_{22}=\mathcal{L}_{\mathcal{\widehat{M}}_{3}},\
\mathcal{\widehat{B}}_{22}=R_{\mathcal{\widehat{N}}_{3}},\\
\label{abc88}
&\begin{array}{l}
 \mathcal{E}_{22}=\mathcal{\widehat{M}}_{3}^{\dagger}*_{N}\mathcal{E}_{3}*_{N}(\mathcal{H}_{3}^{\eta^{*}})^{\dagger}
 +\mathcal{\widehat{S}}_{3}^{\dagger}*_{N}\mathcal{\widehat{S}}_{3}*_{N}\mathcal{H}_{3}^{\dagger}
*_{N}\mathcal{E}_{3}*_{N}\mathcal{\widehat{N}}_{3}^{\dagger}
-\mathcal{H}_{4}^{\dagger}*_{N}\mathcal{E}_{5}*_{N}(\mathcal{J}_{4}^{\eta^{*}})^{\dagger},
 \end{array}\\
 \label{abc99}
&\mathcal{\widehat{\widehat{A}}}_{ii}=\mathcal{R}_{\mathcal{A}_{ii}}*_{N}\mathcal{\widehat{A}}_{ii},\
\mathcal{\widehat{\widehat{B}}}_{ii}=\mathcal{\widehat{B}}_{ii}*_{N}\mathcal{L}_{\mathcal{D}_{ii}},\
\mathcal{\widehat{\widehat{E}}}_{ii}=\mathcal{R}_{\mathcal{A}_{ii}}*_{N}\mathcal{E}_{ii}*_{N}\mathcal{L}_{\mathcal{D}_{ii}},\ (i=1,2),\\
\label{abc9911}
&\mathcal{\overline{A}}_{1}=\begin{bmatrix}-\mathcal{L}_{\mathcal{\widehat{M}}_{1}}*_{N}\mathcal{L}_{\mathcal{\widehat{S}}_{1}}& \mathcal{L}_{\mathcal{F}_{2}}\end{bmatrix},\ \mathcal{\overline{A}}_{2}=\begin{bmatrix}-\mathcal{L}_{\mathcal{\widehat{M}}_{2}}*_{N}\mathcal{L}_{\mathcal{\widehat{S}}_{2}}& \mathcal{L}_{\mathcal{F}_{3}}\end{bmatrix},\ \mathcal{\overline{F}}_{1}=\mathcal{F}^{\dagger}_{2}*_{N}\mathcal{\widehat{S}}_{2},\\
\label{abc9912}
&\mathcal{\overline{B}}_{1}=\begin{bmatrix}-\mathcal{R}_{\mathcal{H}_{1}^{\eta^{*}}}\\ \mathcal{R}_{\mathcal{F}_{2}^{\eta^{*}}}\end{bmatrix},\ \mathcal{\overline{B}}_{2}=\begin{bmatrix}-\mathcal{R}_{\mathcal{H}_{2}^{\eta^{*}}}\\ \mathcal{R}_{\mathcal{F}_{3}^{\eta^{*}}}\end{bmatrix},
\mathcal{\overline{F}}_{2}=\mathcal{F}^{\dagger}_{3}*_{N}\mathcal{\widehat{S}}_{3},\
\mathcal{\overline{G}}_{1}=\mathcal{J}_{2}*_{N}(\mathcal{G}_{2}^{\eta^{*}})^{\dagger},\\
\label{abc99131}
&\mathcal{\overline{G}}_{2}=\mathcal{J}_{3}*_{N}(\mathcal{G}_{3}^{\eta^{*}})^{\dagger},\ \mathcal{\overline{H}}_{1}=\mathcal{L}_{\mathcal{\widehat{M}}_{1}},\ \mathcal{\overline{J}}_{1}=\mathcal{R}_{\mathcal{\widehat{N}}_{1}},\
\mathcal{\overline{H}}_{2}=\mathcal{L}_{\mathcal{\widehat{M}}_{2}},\ \mathcal{\overline{J}}_{2}=\mathcal{R}_{\mathcal{\widehat{N}}_{2}},
 \end{align}
\begin{align} 
\label{abc99132}
&\begin{array}{l}
 \mathcal{\overline{E}}_{1}=-\mathcal{\widehat{M}}_{1}^{\dagger}*_{N}\mathcal{E}_{1}*_{M}(\mathcal{H}_{1}^{\eta^{*}})^{\dagger}
 -\mathcal{\widehat{S}}_{1}^{\dagger}*_{N}\mathcal{\widehat{S}}_{1}*_{N}\mathcal{H}_{1}^{\dagger}
*_{N}\mathcal{E}_{1}*_{N}\mathcal{\widehat{N}}_{1}^{\dagger}
+\mathcal{F}_{2}^{\dagger}*_{N}\mathcal{E}_{2}*_{N}
\\
\ \ \ \ \ \ \ \ \ \ \ (\mathcal{F}_{2}^{\eta^{*}})^{\dagger}-\mathcal{F}_{2}^{\dagger}*_{N}\mathcal{H}_{2}*_{N}\mathcal{\widehat{M}}_{2}^{\dagger}*_{N}\mathcal{E}_{2}
*_{N}(\mathcal{F}_{2}^{\eta^{*}})^{\dagger}
-\mathcal{F}_{2}^{\dagger}*_{N}\mathcal{\widehat{S}}_{2}*_{N}\mathcal{H}_{2}^{\dagger}*_{N}\mathcal{E}_{2}
\\ \ \ \ \ \ \ \ \ \ \ \ \ \ *_{N}\mathcal{\widehat{N}}_{2}^{\dagger}*_{N}\mathcal{J}_{2}^{\eta^{*}}*_{N}(\mathcal{F}_{2}^{\eta^{*}})^{\dagger},
 \end{array}\\
\label{abc991322}
&\begin{array}{l}
 \mathcal{\overline{E}}_{2}=-\mathcal{\widehat{M}}_{2}^{\dagger}*_{N}\mathcal{E}_{2}*_{N}(\mathcal{H}_{2}^{\eta^{*}})^{\dagger}
 -\mathcal{\widehat{S}}_{2}^{\dagger}*_{N}\mathcal{\widehat{S}}_{2}*_{N}\mathcal{H}_{2}^{\dagger}
*_{N}\mathcal{E}_{2}*_{N}\mathcal{\widehat{N}}_{2}^{\dagger}
+\mathcal{F}_{2}^{\dagger}*_{N}\mathcal{E}_{2}*_{N}\\
\ \ \ \ \ \ \ \ \ \ \ (\mathcal{F}_{2}^{\eta^{*}})^{\dagger}
-\mathcal{F}_{3}^{\dagger}*_{N}\mathcal{H}_{3}*_{N}\mathcal{\widehat{M}}_{3}^{\dagger}*_{N}\mathcal{E}_{3}
*_{N}(\mathcal{F}_{3}^{\eta^{*}})^{\dagger}
-\mathcal{F}_{3}^{\dagger}*_{N}\mathcal{\widehat{S}}_{3}*_{N}\mathcal{H}_{3}^{\dagger}*_{N}\mathcal{E}_{3}
\\ \  \ \ \ \ \ \ \ \ \ *_{N}\mathcal{\widehat{N}}_{3}^{\dagger}*_{M}\mathcal{H}_{3}^{\eta^{*}}*_{N}(\mathcal{F}_{3}^{\eta^{*}})^{\dagger},
 \end{array}\\
\label{abxc22}
&\mathcal{\overline{F}}_{ii}=\mathcal{R}_{\mathcal{\overline{A}}_{i}}*_{N}\mathcal{\overline{F}}_{i},\ \mathcal{\overline{G}}_{ii}=\mathcal{\overline{G}}_{i}*_{N}\mathcal{L}_{\mathcal{\overline{B}}_{i}},\ \mathcal{\overline{H}}_{ii}=\mathcal{R}_{\mathcal{\overline{A}}_{i}}*_{N}\mathcal{\overline{H}}_{i},\ \mathcal{\overline{J}}_{ii}=\mathcal{\overline{J}}_{i}*_{N}\mathcal{L}_{\mathcal{\overline{B}}_{i}},\\
\label{abxc33}
&\mathcal{\overline{E}}_{ii}=\mathcal{R}_{\mathcal{\overline{A}}_{i}}*_{N}\mathcal{\overline{E}}_{i}*_{N}\mathcal{L}_{\mathcal{\overline{B}}_{i}},\ \mathcal{\overline{M}}_{ii}=\mathcal{R}_{\mathcal{\overline{F}}_{ii}}*_{N}\mathcal{\overline{H}}_{ii},\ \mathcal{\overline{N}}_{ii}=\mathcal{\overline{J}}_{ii}*_{N}\mathcal{L}_{\mathcal{\overline{G}}_{ii}},\ \mathcal{\overline{S}}_{ii}=\mathcal{\overline{H}}_{ii}*_{N}\mathcal{L}_{\mathcal{\overline{M}}_{ii}},\\
\label{abcD55}
&\mathcal{\overline{\overline{A}}}_{1}=\begin{bmatrix}\mathcal{L}_{\mathcal{\overline{F}}_{11}} &
 -\mathcal{L}_{\mathcal{\overline{M}}_{22}}*_{N}\mathcal{L}_{\mathcal{\overline{S}}_{22}}\end{bmatrix},\
\mathcal{\overline{\overline{B}}}_{1}=\begin{bmatrix}\mathcal{R}_{\mathcal{\overline{G}}_{11}} \\ -\mathcal{R}_{\mathcal{\overline{J}}_{11}}\end{bmatrix},\
\mathcal{\overline{\overline{F}}}_{1}=\mathcal{\overline{F}}_{11}^{\dagger}*_{N}\mathcal{\overline{S}}_{11},\\
\label{abcDE55}
&\mathcal{\overline{\overline{G}}}_{1}=R_{\mathcal{\overline{N}}_{11}}*_{N}\mathcal{\overline{J}}_{11}*_{M}\mathcal{\overline{G}}_{11}^{\dagger},\
\mathcal{\overline{\overline{H}}}_{1}=\mathcal{L}_{\mathcal{\overline{M}}_{22}},\
\mathcal{\overline{\overline{J}}}_{1}=\mathcal{R}_{\mathcal{\overline{N}}_{22}},\\
\label{abcD66}
&\begin{array}{l}
 \mathcal{\overline{\overline{E}}}_{1}=\mathcal{\overline{F}}_{11}^{\dagger}*_{N}\mathcal{\overline{E}}_{11}*_{N}\mathcal{\overline{G}}_{11}^{\dagger}
-\mathcal{\overline{F}}_{11}^{\dagger}*_{N}\mathcal{\overline{H}}_{11}*_{N}\mathcal{\overline{M}}_{11}^{\dagger}*_{N}\mathcal{\overline{E}}_{11}
*_{M}\mathcal{\overline{G}}_{11}^{\dagger}
-\mathcal{\overline{F}}_{11}^{\dagger}*_{N}\mathcal{\overline{S}}_{11}*_{N}\mathcal{\overline{H}}_{11}^{\dagger}\\ \ \ \ \ \ \ \ \ *_{N}\mathcal{\overline{E}}_{11}*_{N}\mathcal{\overline{N}}_{11}^{\dagger}*_{N}\mathcal{\overline{J}}_{11}*_{N}\mathcal{\overline{G}}_{11}^{\dagger}-\mathcal{\overline{M}}_{22}^{\dagger}*_{N}\mathcal{\overline{E}}_{22}*_{M}\mathcal{\overline{J}}_{22}^{\dagger}
 -\mathcal{\overline{S}}_{22}^{\dagger}*_{N}\mathcal{\overline{S}}_{22}*_{N}\mathcal{\overline{H}}_{22}^{\dagger}\\
\ \  \ \ \ \ \ \ \ *_{N}\mathcal{\overline{E}}_{22}*_{N}\mathcal{\overline{N}}_{22}^{\dagger}
 \end{array}\\
\label{abKxc22}
&\mathcal{\overline{\overline{F}}}_{11}=\mathcal{R}_{\mathcal{\overline{\overline{A}}}_{1}}*_{N}\mathcal{\overline{\overline{F}}}_{1},\ \mathcal{\overline{\overline{G}}}_{11}=\mathcal{\overline{\overline{G}}}_{1}*_{N}\mathcal{L}_{\mathcal{\overline{\overline{B}}}_{1}},\ \mathcal{\overline{\overline{H}}}_{11}=\mathcal{R}_{\mathcal{\overline{\overline{A}}}_{1}}*_{N}\mathcal{\overline{\overline{H}}}_{1},\ \mathcal{\overline{\overline{J}}}_{11}=\mathcal{\overline{\overline{J}}}_{1}*_{N}\mathcal{L}_{\mathcal{\overline{\overline{B}}}_{1}},\\
\label{abKKxAc33}
&\mathcal{\overline{\overline{E}}}_{11}=\mathcal{R}_{\mathcal{\overline{\overline{A}}}_{1}}*_{N}\mathcal{\overline{\overline{E}}}_{1}
*_{N}\mathcal{L}_{\mathcal{\overline{\overline{B}}}_{1}},\ \mathcal{\overline{\overline{M}}}_{11}=\mathcal{R}_{\mathcal{\overline{\overline{F}}}_{11}}*_{N}\mathcal{\overline{\overline{H}}}_{11},\ \mathcal{\overline{\overline{N}}}_{11}=\mathcal{\overline{\overline{J}}}_{11}*_{N}\mathcal{L}_{\mathcal{\overline{\overline{G}}}_{11}},\\
\label{abc99911}
&\mathcal{\overline{\overline{S}}}_{11}=\mathcal{\overline{\overline{H}}}_{11}*_{N}\mathcal{L}_{\mathcal{\overline{\overline{M}}}_{11}},\
\mathcal{\widetilde{A}}_{1}=\begin{bmatrix}\mathcal{L}_{\mathcal{\overline{M}}_{11}}*_{N}\mathcal{L}_{\mathcal{\overline{S}}_{11}}& -\mathcal{L}_{\mathcal{\widehat{\widehat{A}}}_{11}}\end{bmatrix},\ \mathcal{\widetilde{A}}_{2}=\begin{bmatrix}\mathcal{L}_{\mathcal{\widehat{\widehat{A}}}_{22}}& -\mathcal{L}_{\mathcal{\overline{F}}_{22}}\end{bmatrix},\\
\label{abc99912}
&\mathcal{\widetilde{B}}_{1}=\begin{bmatrix}\mathcal{R}_{\mathcal{\overline{J}}_{11}}\\-\mathcal{R}_{\mathcal{\widehat{\widehat{B}}}_{11}}\end{bmatrix},\ \mathcal{\widetilde{B}}_{2}=\begin{bmatrix}\mathcal{R}_{\mathcal{\widehat{\widehat{B}}}_{22}}\\ -\mathcal{R}_{\mathcal{\overline{G}}_{22}}\end{bmatrix},
\mathcal{\widetilde{C}}_{1}=\mathcal{L}_{\mathcal{\overline{M}}_{11}}
\mathcal{\widetilde{D}}_{1}=\mathcal{R}_{\mathcal{\overline{N}}_{11}},\\
\label{abc999131}
&\mathcal{\widetilde{C}}_{2}=\mathcal{\overline{F}}^{\dagger}_{22}*_{N}\mathcal{\overline{S}}_{22},\
\mathcal{\widetilde{D}}_{2}=\mathcal{R}_{\mathcal{\overline{N}}_{22}}*_{N}\mathcal{\overline{J}}_{22}*_{N}\mathcal{\overline{G}}^{\dagger}_{22},\\
\label{abcF131}
&\mathcal{\widetilde{E}}_{1}=\widehat{\widehat{\mathcal{A}}}_{11}^{\dagger}*_{N}\widehat{\widehat{\mathcal{E}}}_{11}
*_{M}\widehat{\widehat{\mathcal{B}}}_{11}^{\dagger}-\mathcal{\overline{M}}_{11}^{\dagger}*_{N}\mathcal{\overline{E}}_{11}
*_{N}\mathcal{\overline{J}}_{11}^{\dagger}
 -\mathcal{\overline{S}}_{11}^{\dagger}*_{N}\mathcal{\overline{S}}_{11}*_{N}\mathcal{\overline{H}}_{11}^{\dagger}
*_{N}\mathcal{\overline{E}}_{11}*_{M}\mathcal{\overline{N}}_{11}^{\dagger},\\
\label{abcD66}
&\begin{array}{l}
\mathcal{\widetilde{E}}_{2}=\mathcal{\overline{F}}_{22}^{\dagger}*_{N}\mathcal{\overline{E}}_{22}*_{N}\mathcal{\overline{G}}_{22}^{\dagger}
-\mathcal{\overline{F}}_{22}^{\dagger}*_{N}\mathcal{\overline{H}}_{22}*_{N}\mathcal{\overline{M}}_{22}^{\dagger}*_{N}\mathcal{\overline{E}}_{22}
*_{M}\mathcal{\overline{G}}_{22}^{\dagger}
-\mathcal{\overline{F}}_{22}^{\dagger}*_{N}\\
\ \ \ \ \ \ \ \mathcal{\overline{S}}_{22}*_{N}\mathcal{\overline{H}}_{22}^{\dagger}*_{N}\mathcal{\overline{E}}_{22}*_{N}\mathcal{\overline{N}}_{22}^{\dagger}*_{N}\mathcal{\overline{J}}_{22}*_{M}\mathcal{\overline{G}}_{22}^{\dagger}
-\widehat{\widehat{\mathcal{A}}}_{22}^{\dagger}*_{N}\widehat{\widehat{\mathcal{E}}}_{22}
*_{N}\widehat{\widehat{\mathcal{B}}}_{22}^{\dagger},
 \end{array}\\
\label{abc999WW11}
& \mathcal{\widetilde{F}}_{1}=\begin{bmatrix}\mathcal{L}_{\mathcal{\widetilde{C}}_{11}}& -\mathcal{L}_{\mathcal{\overline{\overline{F}}}_{11}}\end{bmatrix},\
\mathcal{\widetilde{F}}_{2}=\begin{bmatrix}\mathcal{L}_{\mathcal{\overline{\overline{M}}}_{11}}*_{N}\mathcal{L}_{\mathcal{\overline{\overline{S}}}_{11}}& -\mathcal{L}_{\mathcal{\widetilde{C}}_{22}}\end{bmatrix},\
\mathcal{\widetilde{H}}_{1}=\mathcal{\overline{\overline{F}}}_{11}^{\dagger}*_{N}\mathcal{\overline{\overline{S}}}_{11},\\
\label{abc999WWW12}
&\mathcal{\widetilde{J}}_{1}=\mathcal{R}_{\mathcal{\overline{\overline{N}}}_{11}}
*_{N}\mathcal{\overline{\overline{J}}}_{11}*_{N}\mathcal{\overline{\overline{G}}}_{11}^{\dagger}\
\mathcal{\widetilde{G}}_{1}=\begin{bmatrix}\mathcal{R}_{\mathcal{\widetilde{D}}_{11}}\\-\mathcal{R}_{\mathcal{\overline{\overline{G}}}_{11}}\end{bmatrix},\ \mathcal{\widetilde{G}}_{2}=\begin{bmatrix}\mathcal{R}_{\mathcal{\overline{\overline{J}}}_{11}}\\ -\mathcal{R}_{\mathcal{\widetilde{D}}_{22}}\end{bmatrix},
\mathcal{\widetilde{H}}_{2}=\mathcal{L}_{\mathcal{\overline{\overline{M}}}_{11}},\
\mathcal{\widetilde{J}}_{2}=\mathcal{R}_{\mathcal{\overline{\overline{N}}}_{11}},\\
\label{abcDS66}
&\begin{array}{l}
 \mathcal{\widetilde{\widetilde{E}}}_{1}=\mathcal{\overline{\overline{F}}}_{11}^{\dagger}*_{N}\mathcal{\overline{\overline{E}}}_{11}*_{N}\mathcal{\overline{\overline{G}}}_{11}^{\dagger}
-\mathcal{\overline{\overline{F}}}_{11}^{\dagger}*_{N}\mathcal{\overline{\overline{H}}}_{11}*_{N}\mathcal{\overline{\overline{M}}}_{11}^{\dagger}
*_{N}\mathcal{\overline{\overline{E}}}_{11}
*_{N}\mathcal{\overline{\overline{G}}}_{11}^{\dagger}
-\mathcal{\overline{\overline{F}}}_{11}^{\dagger}*_{N}\mathcal{\overline{\overline{S}}}_{11}\\
\ \ \ \ \ *_{N}\mathcal{\overline{\overline{H}}}_{11}^{\dagger}*_{N}\mathcal{\overline{\overline{E}}}_{11}*_{N}\mathcal{\overline{\overline{N}}}_{11}^{\dagger}*_{N}\mathcal{\overline{\overline{J}}}_{11}
*_{M}\mathcal{\overline{\overline{G}}}_{11}^{\dagger}-\mathcal{\widetilde{C}}_{11}^{\dagger}*_{N}\mathcal{\widetilde{E}}_{11}
*_{M}\mathcal{\widetilde{D}}_{11}^{\dagger},
 \end{array}
\end{align}
 \end{subequations}
 \vspace*{-\baselineskip}
 \begin{subequations}
\begin{align}
\label{abcAAD66}
&\mathcal{\widetilde{\widetilde{E}}}_{2}= \mathcal{\widetilde{C}}_{22}^{\dagger}*_{N}\mathcal{\widetilde{E}}_{22}
*_{M}\mathcal{\widetilde{D}}_{22}^{\dagger}
-\mathcal{\overline{\overline{M}}}_{11}^{\dagger}*_{N}\mathcal{\overline{\overline{E}}}_{11}*_{N}\mathcal{\overline{\overline{J}}}_{11}^{\dagger}
 -\mathcal{\overline{\overline{S}}}_{11}^{\dagger}*_{N}\mathcal{\overline{\overline{S}}}_{11}*_{N}\mathcal{\overline{\overline{H}}}_{11}^{\dagger}
*_{N}\mathcal{\overline{\overline{E}}}_{11}*_{M}\mathcal{\overline{\overline{N}}}_{11}^{\dagger},\\
\label{abcDGG66}
&\mathcal{\widetilde{H}}_{11}=\mathcal{R}_{\mathcal{\widetilde{F}}_{1}}*_{N}\mathcal{\widetilde{H}}_{1},\
\mathcal{\widetilde{H}}_{22}=\mathcal{R}_{\mathcal{\widetilde{F}}_{2}}*_{N}\mathcal{\widetilde{H}}_{2},\
\mathcal{\widetilde{J}}_{11}=\mathcal{\widetilde{J}}_{1}*_{N}\mathcal{L}_{\mathcal{\widetilde{G}}_{1}},\
\mathcal{\widetilde{J}}_{22}=\mathcal{\widetilde{J}}_{2}*_{N}\mathcal{L}_{\mathcal{\widetilde{G}}_{2}},\\
\label{abcHHD66}
&\mathcal{\widetilde{\widetilde{E}}}_{11}=\mathcal{R}_{\mathcal{\widetilde{F}}_{1}}*_{N}\mathcal{\widetilde{\widetilde{E}}}_{1}
*_{N}\mathcal{L}_{\mathcal{\widetilde{G}}_{1}},\
\mathcal{\widetilde{\widetilde{E}}}_{22}=\mathcal{R}_{\mathcal{\widetilde{F}}_{2}}*_{N}\mathcal{\widetilde{\widetilde{E}}}_{2}
*_{N}\mathcal{L}_{\mathcal{\widetilde{G}}_{2}},\
\mathcal{\widetilde{A}}=\begin{bmatrix}\mathcal{L}_{\mathcal{\widetilde{H}}_{11}}& -\mathcal{L}_{\mathcal{\widetilde{H}}_{22}}\end{bmatrix},\\
\label{abcH44HRD66}
&\mathcal{\widetilde{B}}=\begin{bmatrix}\mathcal{R}_{\mathcal{\widetilde{J}}_{11}}& -\mathcal{R}_{\mathcal{\widetilde{J}}_{22}}\end{bmatrix},\
\mathcal{\widetilde{E}}=\mathcal{\widetilde{H}}_{22}^{\dagger}*_{N}\mathcal{\widetilde{\widetilde{E}}}_{22}
*_{M}\mathcal{\widetilde{J}}_{22}^{\dagger}-
\mathcal{\widetilde{H}}_{11}^{\dagger}*_{N}\mathcal{\widetilde{\widetilde{E}}}_{11}
*_{N}\mathcal{\widetilde{J}}_{11}^{\dagger}.
\end{align}
\end{subequations}
Then the system \eqref{1.4aa} is consistent if and only if
\begin{align}
\label{b22}
&\mathcal{R}_{\mathcal{\widehat{M}}_{i}}*_{N}\mathcal{R}_{\mathcal{F}_{i}}*_{N}\mathcal{E}_{i}=0,\
\mathcal{R}_{\mathcal{F}_{i}}*_{N}\mathcal{E}_{i}*_{M}\mathcal{L}_{\mathcal{H}_{i}^{\eta^{*}}}=0,\
(i=\overline{1,3}),\\
\label{b44}
&\mathcal{R}_{\mathcal{F}_{4}}*_{N}\mathcal{E}_{4}=0,\
\mathcal{R}_{\mathcal{H}_{4}}*_{N}\mathcal{E}_{5}=0,\
\mathcal{R}_{\mathcal{\widehat{\widehat{A}}}_{kk}}*_{N}\mathcal{\widehat{\widehat{E}}}_{kk}=0,\
\mathcal{\widehat{\widehat{E}}}_{kk}*_{M}\mathcal{L}_{\mathcal{\widehat{\widehat{B}}}_{kk}}=0,\ \\
\label{b66}
&\mathcal{R}_{\mathcal{\overline{M}}_{kk}}*_{N}\mathcal{R}_{\mathcal{\overline{F}}_{kk}}*_{N}\mathcal{\overline{E}}_{kk}=0,\ \mathcal{\overline{E}}_{kk}*_{M}\mathcal{L}_{\mathcal{\overline{G}}_{kk}}*_{M}\mathcal{L}_{\mathcal{\overline{N}}_{kk}}=0,\\
\label{b77}
&\mathcal{R}_{\mathcal{\overline{F}}_{kk}}*_{N}\mathcal{\overline{E}}_{kk}*_{M}\mathcal{L}_{\mathcal{\overline{J}}_{kk}}=0,\
\mathcal{R}_{\mathcal{\overline{H}}_{kk}}*_{N}\mathcal{\overline{E}}_{kk}*_{M}\mathcal{L}_{\mathcal{\overline{G}}_{kk}}=0,\ (k=1,2),\\
\label{b88}
&\mathcal{R}_{\mathcal{\overline{\overline{M}}}_{11}}*_{N}\mathcal{R}_{\mathcal{\overline{\overline{F}}}_{11}}*_{N}\mathcal{\overline{\overline{E}}}_{11}=0,\ \mathcal{\overline{\overline{E}}}_{11}*_{M}\mathcal{L}_{\mathcal{\overline{\overline{G}}}_{11}}*_{M}\mathcal{L}_{\mathcal{\overline{\overline{N}}}_{11}}=0,\\
\label{b99}
&\mathcal{R}_{\mathcal{\overline{\overline{F}}}_{11}}*_{N}\mathcal{\overline{\overline{E}}}_{11}*_{M}\mathcal{L}_{\mathcal{\overline{\overline{J}}}_{11}}=0,\
\mathcal{R}_{\mathcal{\overline{\overline{H}}}_{11}}*_{N}\mathcal{\overline{\overline{E}}}_{11}*_{M}\mathcal{L}_{\mathcal{\overline{\overline{G}}}_{11}}=0,\\
\label{bb111}
&\mathcal{R}_{\mathcal{\widetilde{C}}_{jj}}*_{N}\mathcal{\widetilde{E}}_{jj}=0,\
\mathcal{\widetilde{E}}_{jj}*_{M}\mathcal{L}_{\mathcal{\widetilde{D}}_{jj}}=0\ (j=1,2),\ \mathcal{R}_{\mathcal{\widetilde{H}}_{ll}}*_{N}\mathcal{\widetilde{\widetilde{E}}}_{ll}=0,\\
\label{bb1121}
&\mathcal{\widetilde{\widetilde{E}}}_{ll}*_{M}\mathcal{L}_{\mathcal{\widetilde{G}}_{ll}}=0,\
\mathcal{R}_{\mathcal{\widetilde{A}}}*_{N}\mathcal{\widetilde{E}}*_{M}\mathcal{L}_{\mathcal{\widetilde{B}}},\ (l=1,2),\
\mathcal{R}_{\mathcal{\widetilde{A}}}*_{N}\mathcal{\widetilde{E}}*_{M}\mathcal{L}_{\mathcal{\widetilde{B}}}=0.
\end{align}
Under these conditions, the general solution to system \eqref{1.4aa} can be expressed as follows:
\begin{align}
\label{1.6ssffss}
&\mathcal{Z}_{k}=
\frac{\mathcal{\acute{Z}}_{k}+\mathcal{\acute{Z}}_{k}^{\eta^{*}}}{2},\ \ (k=\overline{1,4}),
\end{align}
where
\begin{align}
\label{sssa}
&\mathcal{\acute{Z}}_{1}=\mathcal{F}_{4}^{\dagger}*_{N}\mathcal{E}_{4}*_{M}(\mathcal{F}_{4}^{\eta^{*}})^{\dagger}
+\mathcal{L}_{\mathcal{F}_{4}}*_{N}\mathcal{W}_{1}
+\mathcal{W}_{2}*_{M}\mathcal{R}_{\mathcal{F}_{4}^{\eta^{*}}},\\
&\mathcal{\acute{Z}}_{4}=\mathcal{H}_{4}^{\dagger}*_{N}\mathcal{E}_{5}*_{M}(\mathcal{H}_{4}^{\eta^{*}})^{\dagger}
+\mathcal{L}_{\mathcal{H}_{4}}*_{N}\mathcal{\acute{W}}_{1}
+\mathcal{W}_{3}*_{M}\mathcal{R}_{\mathcal{H}_{4}^{\eta^{*}}}, \\
&\begin{array}{l}
 \mathcal{\acute{Z}}_{2}=\mathcal{\widehat{M}}_{1}^{\dagger}*_{N}\mathcal{\widehat{E}}_{1}*_{M}(\mathcal{H}_{1}^{\eta^{*}})^{\dagger}
 +\mathcal{\widehat{S}}_{1}^{\dagger}*_{N}\mathcal{\widehat{S}}_{1}*_{N}\mathcal{H}_{1}^{\dagger}
*_{N}\mathcal{E}_{1}*_{N}\mathcal{\widehat{N}}_{1}^{\dagger}+\mathcal{L}_{\mathcal{\widehat{M}}_{1}}
*_{N}\mathcal{L}_{\mathcal{\widehat{S}}_{1}}*_{N}\mathcal{\widehat{U}}_{1}\\
\ \ \ \ \ \ +\mathcal{L}_{\mathcal{\widehat{M}}_{1}}*_{N}\mathcal{\widehat{U}}_{2}*_{M}\mathcal{R}_{\mathcal{\widehat{N}}_{1}}
+\mathcal{\widehat{U}}_{3}*_{M}\mathcal{R}_{\mathcal{H}_{1}^{\eta^{*}}},
 \end{array}\\
&\begin{array}{l}
or\ \mathcal{\acute{Z}}_{2}=\mathcal{F}_{2}^{\dagger}*_{N}\mathcal{E}_{2}*_{M}(\mathcal{F}_{2}^{\eta^{*}})^{\dagger}
-\mathcal{F}_{2}^{\dagger}*_{N}\mathcal{H}_{2}*_{N}\mathcal{\widehat{M}}_{2}^{\dagger}*_{N}\mathcal{E}_{2}
*_{N}(\mathcal{F}_{2}^{\eta^{*}})^{\dagger}
-\mathcal{F}_{2}^{\dagger}*_{N}\mathcal{\widehat{S}}_{2}*_{N}\mathcal{H}_{2}^{\dagger}\\
\ \ \ \ \ *_{N}\mathcal{E}_{2}*_{N}\mathcal{\widehat{N}}_{2}^{\dagger}*_{M}\mathcal{H}_{2}^{\eta^{*}}*_{N}(\mathcal{F}_{2}^{\eta^{*}})^{\dagger}-\mathcal{F}_{2}^{\dagger}
*_{N}\mathcal{\widehat{S}}_{2}
*_{N}\mathcal{\widehat{V}}_{2}*_{M}\mathcal{R}_{\mathcal{\widehat{N}}_{2}}*_{M}\mathcal{H}_{2}^{\eta^{*}}*_{M}(\mathcal{G}_{2}^{\eta^{*}})^{\dagger}
\\
\ \ \ \ \ +\mathcal{L}_{\mathcal{F}_{2}}*_{N}\mathcal{\widehat{V}}_{4}+\mathcal{\widehat{V}}_{5}*_{M}\mathcal{R}_{\mathcal{F}_{2}},
 \end{array}\\
&\begin{array}{l}
 \mathcal{\acute{Z}}_{3}=\mathcal{\widehat{M}}_{2}^{\dagger}*_{N}\mathcal{E}_{2}*_{N}(\mathcal{H}_{2}^{\eta^{*}})^{\dagger}
 +\mathcal{\widehat{S}}_{2}^{\dagger}*_{N}\mathcal{\widehat{S}}_{2}*_{N}\mathcal{H}_{2}^{\dagger}
*_{N}\mathcal{E}_{2}*_{M}\mathcal{\widehat{N}}_{2}^{\dagger}+\mathcal{L}_{\mathcal{\widehat{M}}_{2}}
*_{N}\mathcal{L}_{\mathcal{\widehat{S}}_{2}}*_{N}\mathcal{\widehat{V}}_{1}\\
\ \ \ \ \ \ +\mathcal{L}_{\mathcal{\widehat{M}}_{2}}*_{N}\mathcal{\widehat{V}}_{2}*_{M}\mathcal{R}_{\mathcal{\widehat{N}}_{2}}
+\mathcal{\widehat{V}}_{3}*_{M}\mathcal{R}_{\mathcal{H}_{2}^{\eta^{*}}},
 \end{array}\\
&\begin{array}{l}
or\  \mathcal{\acute{Z}}_{3}=\mathcal{F}_{3}^{\dagger}*_{N}\mathcal{E}_{3}*_{M}(\mathcal{F}_{3}^{\eta^{*}})^{\dagger}
-\mathcal{F}_{3}^{\dagger}*_{N}\mathcal{H}_{3}*_{N}\mathcal{\widehat{M}}_{3}^{\dagger}*_{N}\mathcal{E}_{3}
*_{N}(\mathcal{F}_{3}^{\eta^{*}})^{\dagger}
-\mathcal{F}_{3}^{\dagger}*_{N}\mathcal{\widehat{S}}_{3}*_{N}\mathcal{H}_{3}^{\dagger}\\
\ \ \ \ \ *_{N}\mathcal{E}_{3}*_{N}\mathcal{\widehat{N}}_{3}^{\dagger}*_{N}(\mathcal{H}_{3}^{\eta^{*}})*_{N}(\mathcal{F}_{3}^{\eta^{*}})^{\dagger}-\mathcal{F}_{3}^{\dagger}
*_{N}\mathcal{\widehat{S}}_{3}
*_{N}\mathcal{\widehat{K}}_{2}*_{M}\mathcal{R}_{\mathcal{\widehat{N}}_{3}}*_{M}\mathcal{H}_{3}^{\eta^{*}}*_{N}(\mathcal{F}_{3}^{\eta^{*}})^{\dagger}
\\
\ \ \ \ \ +\mathcal{L}_{\mathcal{F}_{3}}*_{N}\mathcal{\widehat{K}}_{4}+\mathcal{\widehat{K}}_{5}*_{M}\mathcal{R}_{\mathcal{F}_{3}^{\eta^{*}}},\ (i=\overline{1,3}).
 \end{array}
\end{align}
Where the arbitrary tensors $\mathcal{W}_{j}$, $\mathcal{\widehat{V}}_{i}$, $\mathcal{\widehat{U}}_{j}$, $\mathcal{\widehat{K}}_{k}$ and $\mathcal{\acute{W}}_{1}$ $(j=\overline{1,3},\ i=\overline{1,5},\ k\in \{2,4,5\})$ can be reduced by \eqref{Xyqs1}-\eqref{XZyw996} under the definitions \eqref{abc5522}-\eqref{abcH44HRD66}.
\begin{proof} Consider the following quaternion system of tensor equations:
\begin{align}
\label{1.6aaa}
\left\{
\begin{array}{rll}
\begin{gathered}
\mathcal{F}_{4}*_{N}\mathcal{\acute{Z}}_{1}*_{N}\mathcal{F}_{4}^{\eta^{*}}=\mathcal{E}_{4},\\
\mathcal{F}_{i}*_{N}\mathcal{\acute{Z}}_{i}*_{N}\mathcal{F}_{i}^{\eta^{*}}
+\mathcal{H}_{i}*_{N}\mathcal{\acute{Z}}_{i+1}*_{N}\mathcal{H}_{i}^{\eta^{*}}=\mathcal{E}_{i},\\
 \mathcal{H}_{4}*_{N}\mathcal{\acute{Z}}_{4}*_{N}\mathcal{H}_{4}^{\eta^{*}}=\mathcal{E}_{5},
\end{gathered}
\end{array}
  \right.
\end{align}
where $(i=\overline{1,3})$. Suppose  that the system \eqref{1.6aa} is consistent. Claim that $(\mathcal{Z}_{1}, \mathcal{Z}_{2}, \mathcal{Z}_{3}, \mathcal{Z}_{4})$ is a solution to the quaternion system of tensor equations \eqref{1.6aa}, then it is evident that  $(\mathcal{\acute{Z}}_{1}, \mathcal{\acute{Z}}_{2}, \mathcal{\acute{Z}}_{3}, \mathcal{\acute{Z}}_{4})$ $=(\mathcal{Z}_{1}, \mathcal{Z}_{2}, \mathcal{Z}_{3}, \mathcal{Z}_{4})$ is a solution to the system \eqref{1.6aaa}. Conversely, if the system \eqref{1.6aaa} has a solution $(\mathcal{\acute{Z}}_{1}, \mathcal{\acute{Z}}_{2}, \mathcal{\acute{Z}}_{3}, \mathcal{\acute{Z}}_{4})$. It is sufficient to show that
\begin{align}
\label{1.6ssss}
&(\mathcal{Z}_{1},\mathcal{Z}_{2},\mathcal{Z}_{3}, \mathcal{Z}_{4})=
\left(\frac{\mathcal{\acute{Z}}_{1}+\mathcal{\acute{Z}}_{1}^{\eta^{*}}}{2},
\frac{\mathcal{\acute{Z}}_{2}+\mathcal{\acute{Z}}_{2}^{\eta^{*}}}{2},
\frac{\mathcal{\acute{Z}}_{3}+\mathcal{\acute{Z}}_{3}^{\eta^{*}}}{2},
\frac{\mathcal{\acute{Z}}_{4}+\mathcal{\acute{Z}}_{4}^{\eta^{*}}}{2}\right),
\end{align}
is a solution to system \eqref{1.6aa}. Clearly, the quaternion tensors $\mathcal{Z}_{i},$  $(i=\overline{1,4})$ are  $\eta$-Hermitian tensors. By Applying \eqref{1.6ssss} on the system \eqref{1.6aa} yields:
\begin{align*}
&\mathcal{F}_{4}*_{N}\mathcal{Z}_{1}*_{N}\mathcal{F}_{4}^{\eta^{*}}
=\mathcal{F}_{4}*_{N}\left(\frac{\mathcal{\acute{Z}}_{1}+\mathcal{\acute{Z}}_{1}^{\eta^{*}}}{2}\right)*_{N}F_{4}^{\eta^{*}}\\
&=\frac{1}{2}\mathcal{F}_{4}*_{N}\mathcal{\acute{Z}}_{1}*_{N}\mathcal{F}_{4}^{\eta^{*}}
+\frac{1}{2}\left(\mathcal{F}_{4}*_{N}\mathcal{\acute{Z}}_{1}*_{N}\mathcal{F}_{4}^{\eta^{*}}\right)^{\eta^{*}}
=\mathcal{E}_{4}.
\end{align*}
Similarly, it can be shown that
\begin{align*}
&\mathcal{H}_{4}*_{N}\mathcal{Z}_{4}*_{N}\mathcal{H}_{4}^{\eta^{*}}=\mathcal{E}_{5}.
 \end{align*}
Moreover,
\begin{align*}
&\mathcal{F}_{i}*_{N}\mathcal{Z}_{i}*_{N}\mathcal{F}_{i}^{\eta^{*}}
+\mathcal{H}_{i}*_{N}\mathcal{Z}_{i+1}*_{N}\mathcal{H}_{i}^{\eta^{*}}\\
&=\mathcal{F}_{i}*_{N}\left(\frac{\mathcal{\acute{Z}}_{i}+\mathcal{\acute{Z}}_{i}^{\eta^{*}}}{2}\right)*_{N}\mathcal{F}_{i}^{\eta^{*}}
+\mathcal{H}_{i}*_{N}\left(\frac{\mathcal{\acute{Z}}_{i+1}+\mathcal{\acute{Z}}_{i+1}^{\eta^{*}}}{2}\right)*_{N}\mathcal{H}_{i}^{\eta^{*}}\\
&=\frac{1}{2}\left[ \mathcal{F}_{i}*_{N}\mathcal{\acute{Z}}_{i}*_{N}\mathcal{F}_{i}^{\eta^{*}}
+\mathcal{H}_{i}*_{N}\mathcal{\acute{Z}}_{i+1}*_{N}\mathcal{H}_{i}^{\eta^{*}}\right]\\
&+\frac{1}{2}\left[ \mathcal{F}_{i}*_{N}\mathcal{\acute{Z}}_{i}*_{N}\mathcal{F}_{i}^{\eta^{*}}
+\mathcal{H}_{i}*_{N}\mathcal{\acute{Z}}_{i+1}*_{N}\mathcal{H}_{i}^{\eta^{*}}\right]^{\eta^{*}}=\mathcal{E}_{i},\ (i=\overline{1,3}).
\end{align*}
 Therefore, \eqref{1.6ssss} is a solution to the system \eqref{1.6aa}. Consequently, apply $Theorem$ $\ref{system eeee1},$ on the system \eqref{1.6aaa}, we can establish the solvability conditions and the general solution to the quaternion system \eqref{1.6aa}.
\end{proof}
\end{corollary}

\section{\textbf{Conclusion}}
Having first established the necessary and sufficient conditions  for the presence of a solution to \eqref{1.4aa}, we, therefore, manifest an expression of its general solution. If $\mathcal{A}_{i}=\mathcal{D}_{i}=0$ in \eqref{1.7aa}, where $(i=\overline{1,3}),$ we obtain the Sylvester-like quaternion system of tensor equations \eqref{1.5aa}. As an application of system \eqref{1.5aa}, we investigate an $\eta$-Hermitian solution to system \eqref{1.6aa}. We also construct a numerical example to validate the system \eqref{1.4aa}. It is notable that the primary conclusions of this study are particularly beneficial for the corresponding systems over the real and complex number fields. These conclusions can also obtain the corresponding matrix equation systems to \eqref{1.4aa}-\eqref{1.6aa}.

All results are valid over an arbitrary division ring. As a  direct consequence, the corresponding systems of quaternion matrix equations to the systems \eqref{1.4aa}, \eqref{1.5aa}, \eqref{1.6aa}, and \eqref{1.7aa} can be described by rank equalities and Moore-Penrose inverses of matrices whenever $N=M=1$.

In further future work, we infer the solvability constraints to the $n$-system of the matrix equations
\begin{small}
\begin{align*}
\left\{
\begin{gathered}
A_{1}X_{1}B_{1}+C_{1}Y_{1}D_{1}+C_{1}(G_{1}Z_{1}F_{1}+H_{1}Z_{2}J_{1})B_{1}=E_{1}\\
A_{2}X_{2}B_{2}+C_{2}Y_{2}D_{2}+C_{2}(G_{2}Z_{2}F_{2}+H_{2}Z_{3}J_{2})B_{2}=E_{2}\\
\vdots\\
A_{n}X_{n}B_{n}+C_{n}Y_{n}D_{n}+C_{n}(G_{n}Z_{n}F_{n}+H_{n}Z_{n+1}J_{n})B_{n}=E_{n}
\end{gathered}
\right.
\end{align*}
\end{small}
can be characterized by rank equalities and  Moore-Penrose inverses of some known matrices and hence, we can derive a formula of its general solution. Moreover, we intend to study that system over an arbitrary regular ring.

\end{document}